\documentclass[notitlepage]{report}
\usepackage{tikzexternal}

\def\thisChapterRef{}

\usepackage{amsmath}
\usepackage{amsthm}
\usepackage{pxfonts}
\usepackage{hyphenat}
\usepackage{breqn}
\usepackage{listings}
\usepackage{tensor}
\usepackage{mathtools}
\usepackage{fancyvrb}
\usepackage{xfrac}
\usepackage{hyperref}
\usepackage{cleveref}

\makeatletter
\AtBeginDocument{%
  \let\oldlabel=\label
  \let\oldref=\ref
  \let\oldcref=\cref
  \let\oldCref=\Cref
  \renewcommand\ref[1]{%
    \oldref{\thisChapterRef#1}}%
  \renewcommand\cref[1]{%
    \oldcref{\thisChapterRef#1}}%
  \renewcommand\Cref[1]{%
    \oldCref{\thisChapterRef#1}}%
  \renewcommand\label[1]{%
    \oldlabel{\thisChapterRef\@firstofone#1}}%
  \let\ltx@label=\label
}
\makeatother

\theoremstyle{plain}
\newtheorem{theorem}{Theorem}[section]
\newtheorem{proposition}[theorem]{Proposition}
\theoremstyle{definition}
\newtheorem{defn}[theorem]{Definition}
\theoremstyle{remark}
\newtheorem{question}[theorem]{Question}

\newcommand\m[1]{\mathcal{#1}}
\newcommand\R{\mathbb{R}}
\newcommand\Z{\mathbb{Z}}

\makeatletter
\providecommand{\BreakableEnDash}{\leavevmode%
    \prw@zbreak--\discretionary{}{}{}\prw@zbreak}
  \DeclareRobustCommand{\enhyp}{%
    \ifmmode--\else\BreakableEnDash\fi}
  \newcommand{\BreakableEmDash}{\leavevmode%
    \prw@zbreak---\discretionary{}{}{}\prw@zbreak}
  \DeclareRobustCommand{\emhyp}{%
    \ifmmode---\else\BreakableEmDash\fi}
\makeatother
\newcommand{\noproof}{%
    \ifmmode
    \pushQED{\qed}\qedhere
    \else
    {\hspace*{\fill}\qed}%
    \fi}

\let\centring=\centering
\newcommand{\ssetminus}{\!\smallsetminus\!}

\makeatletter
\lst@CCPutMacro
    \lst@ProcessOther {"2D}{\lst@ttfamily{-{}}{--}}
    \@empty\z@\@empty
\makeatother

\tikzsetexternalprefix{figures/envelope-}
\tikzsetfigurename{figure}
\let\origtikzsetnextfilename\tikzsetnextfilename
\def\tikzsetnextfilename#1{%
  \origtikzsetnextfilename{#1}
  \mysetlabel{#1}
}

\makeatletter
\long\def\chappendix#1\end#2{%
  \g@addto@macro\appendices{#1}
  \end{#2}
}

\gdef\appendices{}
\newcommand\inputchapter[1]{%
  \def\thisChapterRef{#1-}%
  \tikzsetexternalprefix{figures/envelope-#1-}
  \let\appendix=\chappendix
  \g@addto@macro\appendices{\def\thisChapterRef{#1-}}%
  \input{#1}%
}
\makeatother

\newcommand{\mysetlabel}[1]{%
  \gdef\mynextlabel{#1}}

\newcommand{\autolabel}{%
  \label{fig:\mynextlabel}
  \global\let\mynextlabel\relax
}

\newcommand{\kb}[2]{\langle \csname kb#1\endcsname(#2) \rangle}

\title{Investigations of Higher Order Links}

\author{Nils A.\ Baas \and Andrew Stacey}

\begin{document}
\hypersetup{pageanchor=false}
\maketitle
\begin{abstract}
The following is an amalgamation of four preprints and some computer programs which together represent the current state of our investigations of higher order links.
This investigation was motivated by questions discussed and raised in \cite{nbNewSta}.
An important motivation has been to suggest the synthesis of new types of molecules (see \cite{nbNewStr,nbNewSta,nbStrOr,nbHigOrd,nbnsChemSynth,nbnsas} in the bibliography).
This discussion is not final, but we think that the results and methods are worth making public and would be useful for other investigators.
\end{abstract}

\tableofcontents

\begingroup
\long\def\documentclass#1\myabstract#2\begin#3{
\begingroup
\begin{abstract}
#2
\end{abstract}
}
\let\enddocument=\endinput
\let\mymaketitle=\relax

\chapter{Operads}

\mymaketitle

\section{Introduction}

\section{The Framed Link Operad}

In this section we define the framed link operad and show how to use it to define the notion of the level of a link.

\begin{defn}
We define the \emph{framed link operad}, which we shall denote by \(\m{L}\), in the following way.
The objects of \(\m{L}\) are the natural numbers.
The morphisms \(\m{L}(n,1)\) are smooth embeddings of \(n\) copies of the thickened torus, \(S^1 \times D^2\), in itself with the property that the boundary of the source is taken into the interior of the target.
Composition of morphisms corresponds to composition of embeddings.

We consider this as an operad enriched in the category of smooth spaces.
\end{defn}

When we restrict to the endomorphism monoid \(\m{L}(1,1)\), we recover the notion of a satelite knot (modulo the fact that our knots are in the torus rather than \(\R^3\)).

Given a morphism \(f \in \m{L}(n,1)\), we can try to decompose it as an operadic composition.
That is, we look for \(f_k\) and \(g_1, \dotsc, g_k\) such that
\[
  f = f_k \circ (g_1, \dotsc, g_k)
\]
and then we repeat this process on \(f_k\) and on the \(g_j\).
When considering these decompositions, we wish to disallow ``trivial'' decompositions.

\begin{defn}
A \emph{trivial decomposition} is one of the form
\[
  f = h \circ g, \quad f = g \circ (h_1, \dotsc, h_k), \quad f = f_k \circ (h_1 \circ g_1, \dotsc, h_k \circ g_k)
\]
with \(h, h_i \in \m{L}(1,1)\) where removing the, as appropriate, \(h\) or all of the \(h_i\) does not change the isotopy class of the morphism.

A \emph{globally trivial decomposition} is one of the form
\[
  f = g \circ (h_1, \dotsc, h_k), \quad f = f_k \circ (h_1 \circ g_1, \dotsc, h_k \circ g_k)
\]
with \(h_i \in \m{L}(1,1)\) where removing \emph{one} of the \(h_i\) does not change the isotopy class of the morphism.
\end{defn}

That is to say, the morphisms are in the same path component of \(\m{L}(n,1)\).
Using the notation of the definition, for a trivial decomposition we test if \(f \simeq g\) in the first two cases, and \(f \simeq f_k \circ (g_1, \dotsc, g_k)\) in the third.
For a globally trivial decomposition, we test this with only removing one of the \(h_i\).

The reason for the name ``globally'' is related to how we intend to use these notions.
We shall start with a morphism and keep decomposing it until the only way to decompose it further is by adding trivial decompositions, or globally trivial decompositions.
If we use trivial decompositions, we can continue decomposing so long as one input of the morphism admits further decomposition.
Thus we can continue to decompose the morphism provided we can locally do so.
When using globally trivial decompositions, every input has to admit further decomposition and so we decompose the morphism only so long as we can do so globally.

There is an obvious ordering on the family of decompositions (omitting decompositions with (globally) trivial subdecompositions) of a given morphism given by refinement.
We then consider maximal elements in this family of decompositions.
Such a decomposition can be represented as a tree, where the nodes correspond to morphisms.

\begin{defn}
We define the \emph{level} of a morphism \(f \in \m{L}(n,1)\) to be the maximum of the heights of maximal elements in its family of decompositions, when we disregard decompositions that have trivial subdecompositions.

We define the \emph{global level} of a morphism \(f \in \m{L}(n,1)\) to be the maximum of the heights of maximal elements in its family of decompositions, when we disregard decompositions that have globally trivial subdecompositions.
\end{defn}

The basic idea behind these definitions is to view a link through a series of filters by truncating its decomposition at a certain height.
At each stage, we ignore the lower level structure and think of it as occuring on a scale too small to see.
As we move through the levels of filtration, we see finer and finer detail on each component.
In the first stages, we see finer detail on every component.
There comes a point where some components have revealled all of their secrets and increasing the filtration power has no effect on them.
However, it may still have an effect on other components.
The inner level is this point at which some components have been fully refined.
The outer level is the point at which all the components have been fully refined.

As an illustration, consider the Whitehead link as in Figure~\ref{fig:whitehead}.
By drawing it as in Figure~\ref{fig:whiteheadv2}, we see that this can be factored through the Hopf link by replacing one of the circles in the Hopf link by a morphism in \(\m{L}(1,1)\).
Neither component admits further decomposition and thus it has level \(2\) and global level \(1\).

\tikzsetnextfilename{whitehead}
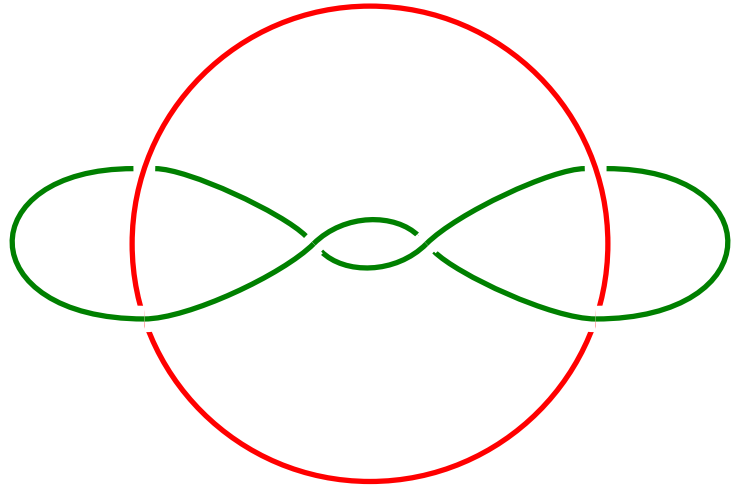
\begin{figure}
\centering
\begin{tikzpicture}[every path/.style={thick knot,double=ring2},every node/.style={transform shape,knot crossing,inner sep=1.5pt}]
\pgfmathsetmacro{\whx}{3}
\pgfmathsetmacro{\why}{1}
\pgfmathsetmacro{\radius}{sqrt(\whx^2 + \why^2)}
\draw[double=ring1] (0,0) circle (\radius);
\node (l) at (-.75,0) {};
\node (r) at (.75,0) {};
\node (ul) at (-\whx,\why) {};
\node (ur) at (\whx,\why) {};
\node[fill=white] (ll) at (-\whx,-\why) {};
\node[fill=white] (lr) at (\whx,-\why) {};
\draw (l.center) .. controls (l.4 north east) and (r.4 north west) .. (r);
\draw (l) .. controls (l.4 south east) and (r.4 south west) .. (r.center);
\draw (r) .. controls (r.4 south east) and (lr.4 west) .. (lr.center);
\draw (r.center) .. controls (r.4 north east) and (ur.4 west) .. (ur);
\draw (l.center) .. controls (l.4 south west) and (ll.4 east) .. (ll.center);
\draw (l) .. controls (l.4 north west) and (ul.4 east) .. (ul);
\draw (ul) .. controls (ul.16 west) and (ll.16 west) .. (ll.center);
\draw (ur) .. controls (ur.16 east) and (lr.16 east) .. (lr.center);
\end{tikzpicture}
\caption{The Whitehead Link}
\autolabel
\end{figure}

\tikzsetnextfilename{whiteheadv2}
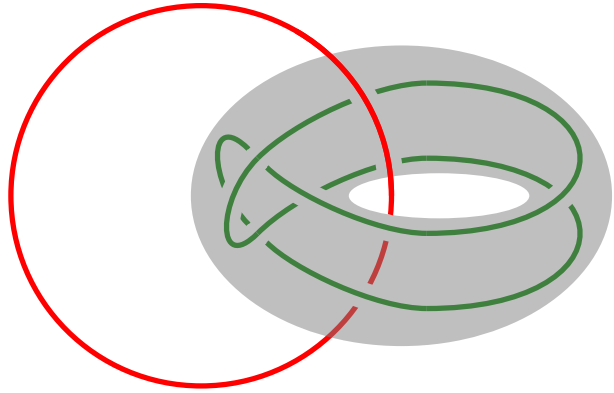
\begin{figure}
\centering
\begin{tikzpicture}[every path/.style={thick knot,double=ring2},every node/.style={transform shape,knot crossing,inner sep=2pt}]
\pgfmathsetmacro{\whx}{3}
\pgfmathsetmacro{\why}{1}
\pgfmathsetmacro{\radius}{sqrt(\whx^2 + \why^2)}
\pgfmathsetmacro{\whyy}{2}
\pgfmathsetmacro{\whxx}{sqrt(\whx^2 + \why^2 - \whyy^2)}
\node (l) at (.75,1) {};
\node (r) at (.75,0) {};
\node (ul) at (\whx,\why + 1) {};
\node (ur) at (\whx,\why) {};
\node (ll) at (\whx,0) {};
\node (lr) at (\whx,-\why) {};
%\draw (l.center) .. controls (l.4 north) and (r.4 north west) .. (r);
%\draw (l) .. controls (l.4 south east) and (r.4 south west) .. (r.center);
\draw (r.center) .. controls (r.4 north east) and (ur.4 west) .. (ur.center);
\draw (l.center) .. controls (l.4 north east) and (ul.4 west) .. (ul.center);
\draw (ur.center) .. controls (ur.16 east) and (lr.16 east) .. (lr.center);
\draw (ul.center) .. controls (ul.16 east) and (ll.16 east) .. (ll.center);
\draw[double=ring1] (0,.5) circle (.8*\radius);
\draw (l) .. controls (l.4 south east) and (ll.4 west) .. (ll.center);
\draw (r) .. controls (r.4 south east) and (lr.4 west) .. (lr.center);
\draw (r) to[out=135,in=135,looseness=2] (l);
\draw (r.center) to[out=-135,in=-135,looseness=2] (l.center);
\fill[opacity=.5,gray,even odd rule] (\radius-.5,.5) ellipse (2.8 and 2) (\radius,.5) ellipse (1.2 and .3);
\draw[double=none,ring1,line width=2pt] (.8*\radius,.5) arc(0:60:.8*\radius);
\end{tikzpicture}
\caption{The Deformed Whitehead Link}
\autolabel
\end{figure}

\begin{question}
How unique is the decomposition?
We would like to think that it was unique up to some reasonable idea of equivalence.
\end{question}

Given two components in a link, we can consider their relative level.
That is, we consider the tree associated to a maximal decomposition (either with the local or global notion).
The components of the link correspond to leaves on the tree so our two components specify two leaves.
We then consider the node at which the branches ending at those two leaves diverged.
The height of the subtree with this node as root gives a measure of the disconnectedness of the two components.
Note that the two components can themselves be linked.
In the doubled Hopf link (with the shown decomposition), as in Figure~\ref{fig:dblehopf}, the red and blue components have relative level \(1\) even though they are linked.

\tikzsetnextfilename{dblehopf}

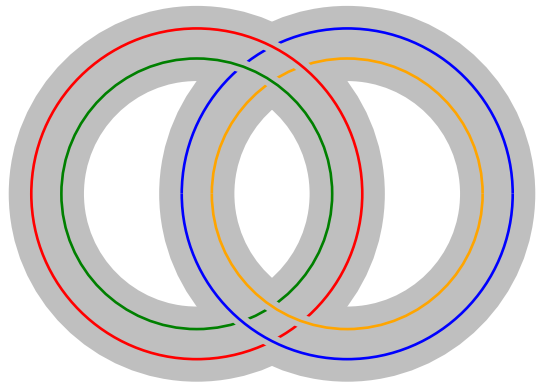
\begin{figure}
\centering
\begin{tikzpicture}
\resetbgcolour{torus}
\draw[torus,line width=1cm] (-1,0) circle (2);
\draw[torus,line width=1cm] (1,0) circle (2);
\draw[knot,double=ring3] (1,0) ++(2.2,0) arc(0:180:2.2);
\draw[knot,double=ring4] (1,0) ++(1.8,0) arc(0:180:1.8);
\draw[knot,double=ring1] (-1,0) circle (2.2);
\draw[knot,double=ring2] (-1,0) circle (1.8);
\draw[knot,double=ring3] (1,0) ++(2.2,0) arc(0:-180:2.2);
\draw[knot,double=ring4] (1,0) ++(1.8,0) arc(0:-180:1.8);
\end{tikzpicture}
\caption{Doubled Hopf link}
\autolabel
\end{figure}

\section{Pure Links}

The framed link operad covers all knots and links.
As such, it is too big a home for the links that we are particularly interested in.
Our purpose is to study the composition structure of this operad and see what happens as we build higher and higher order links.
In order to better examine what happens under composition, we wish to remove other sources of complexity.
One of the main such sources is the complexity in \(\m{L}(1,1)\).
Although there is a rich composition structure here (studied as the subject of satelite knots), it is intertwined with the question of classification of knots which is somewhat orthogonal to our quest.
Thus we wish to remove the knottedness of our links and concentrate on them purely as links.

One approach to this is due to Milnor~\cite{MR0071020} where the components of links are allowed to pass through themselves, but not through each other.
This defines a quotient of each morphism space.
However, the resulting spaces do not carry a well\hyp{}defined composition.
To see this, we compare the Whitehead link from Figure~\ref{fig:whitehead} and the result obtained by doubling the  components.
This is no longer equivalent to the unlink under Milnor's relation since we cannot pass the green and blue components through each other.
See \cite[Section 3]{MR0071020}.

\tikzsetnextfilename{whiteheadrepl}
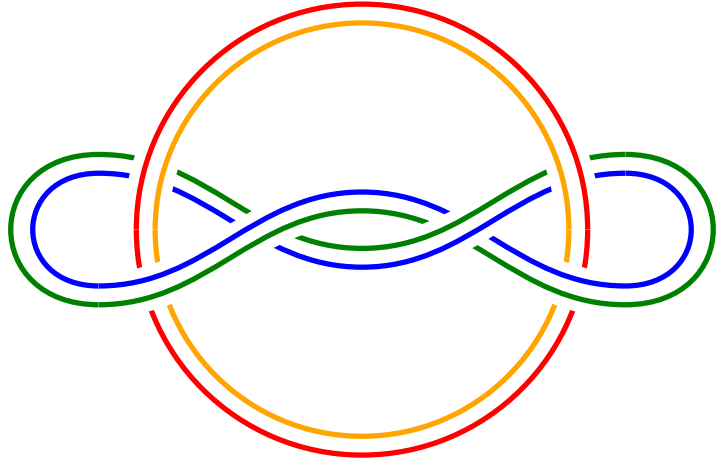
\begin{figure}
\centering
\begin{tikzpicture}[every path/.style={thick knot,double=ring2},every node/.style={transform shape,knot crossing,inner sep=1.5pt}]
\pgfmathsetmacro{\radius}{3}
\draw[double=ring1] (-\radius,0) arc(-180:0:\radius);
\draw[double=ring4] (-\radius+\brstep,0) arc(-180:0:\radius-\brstep);
\foreach \brj in {1,-1} {
\begin{scope}[scale=\brj]
\foreach \brk in {2,3} {
\draw[double=ring\brk] (0,\brk * \brstep - \brstep) to[out=0,in=180] (3.5,-1.5+\brk * \brstep) to[out=0,in=0,looseness=2] (3.5,1.5-\brk*\brstep);
}
\foreach \brk in {2,3} {
\draw[double=ring\brk] (0,-\brk * \brstep + \brstep) to[out=0,in=180] (3.5,1.5-\brk * \brstep); 
}
\end{scope}
}
\draw[double=ring1] (-\radius,0) arc(180:0:\radius);
\draw[double=ring4] (-\radius+\brstep,0) arc(180:0:\radius-\brstep);
\end{tikzpicture}
\caption{The doubled Whitehead link}
\autolabel
\end{figure}

Our approach is to take a suboperad.

\begin{defn}
The \emph{pure link operad}, denoted \(\m{P L}\), is the suboperad of \(\m{L}\) wherein for \(n \ge 2\) the morphisms in \(\m{P L}(n,1)\) are those smooth embeddings \(\coprod_n (S^1 \times D^2) \to S^1 \times D^2\) with the property that the restriction to any component is isotopic to the unknot, whilst in \(\m{P L}(1,1)\) we take morphisms isotopic to the identity.
\end{defn}

To illustrate the difference between the suboperad and Milnor's quotient, consider the link in Figure~\ref{fig:brring2cptsdef}.
As we shall see later, this is the Brunnian ring with two components.
In Milnor's classification, this is the unlink with two components since the green component can pass through itself and unlink from the red component.
In our classification, this move is not allowed and it is not equivalent to an unlink.
To show that it is an object in the pure link operad, let us redraw it as in Figure~\ref{fig:brring2cpts} with the torus shaded in.
In the torus, the unlinking in Milnor's classification still works: bring the upper ends of the green component round the torus to meet the lower ends and pass them through.
This unhooks the green component from the red.

\tikzsetnextfilename{brring2cptsdef}

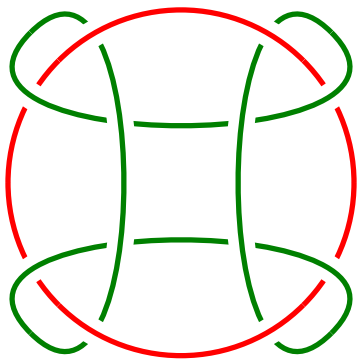
\begin{figure}
\centering
\begin{tikzpicture}[every path/.style={thick knot},every to/.style={looseness=1.5}]
\draw[double=ring1] (2.3,0) arc (0:45:2.3) (2.3,0) arc (0:-45:2.3);
\draw[double=ring1] (-2.3,0) arc (180:135:2.3) (-2.3,0) arc (-180:-135:2.3);
\draw[double=ring2] (-2,2) to[out=-135,in=-45] (2,2);
\draw[double=ring2] (2,-2) to[out=45,in=135] (-2,-2);
\draw[double=ring2] (2,2) to[out=135,in=-135] (2,-2);
\draw[double=ring2] (-2,-2) to[out=-45,in=45] (-2,2);
\draw[double=ring1] (0,2.3) arc (90:45:2.3) (0,2.3) arc (90:135:2.3);
\draw[double=ring1] (0,-2.3) arc (-90:-45:2.3) (0,-2.3) arc (-90:-135:2.3);
\end{tikzpicture}
\caption{Deformed Brunnian ring with 2 components}
\autolabel
\end{figure}

\tikzsetnextfilename{brring2cpts}

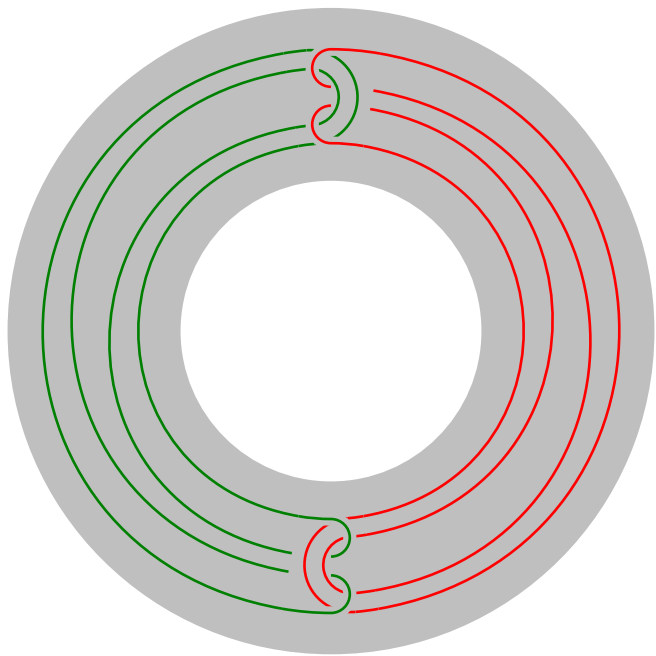
\begin{figure}
\centering
\begin{tikzpicture}
\resetbgcolour{torus}
\fill[torus,even odd rule] (0,0) circle (4.3) circle (2);
\foreach \brk in {1,2} {
\begin{scope}[rotate=\brk * 180]
\colorlet{chain}{ring\brk}
\brunnianlink{2.5}{180}
\end{scope}
}
\end{tikzpicture}
\caption{Brunnian ring with 2 components embedded in a torus}
\autolabel
\end{figure}

A variant of this suboperad is to allow twistings.
That is, we add in morphisms which differ from one already there by a twist.
To see what difference this would make, consider the Hopf link.
With twists allowed, this factors as the trivial link with \(2\) components followed by a twist, as shown in Figure~\ref{fig:twisthopf}.
Without twists, this does not factor.

As we only allow things in \(\mathcal{P L}(1,1)\) that are isotopic to the identity, in a factorisation these play no part.
Therefore when we consider a factorisation of a pure link each non\hyp{}trivial term in the factorisation must involve two or more of the inputs.
Thus when we take the graph of the factorisation, it is a rooted tree in which each node is at least \(2\)\enhyp{}valent.
At each level of the factorisation, therefore, the number of inputs must strictly decrease.
Hence the level of a pure link can be at most one less than the number of components.

\bigskip

This suboperad is still not quite what we want.
It contains more links than we would like.
Specifically, we would like our links to notice that they are in a torus.
If it is possible to cut the torus and not disturb the link, then this link is really in \(S^3\) or \(\R^3\) and not in the torus.
We would like to exclude this possibility.
Thus we study links that actually circumnavigate the torus.
Indeed, this is the primary property that we would like our links to have.
As such, we wish to ensure not just that we have this property, but that this property is somehow a central property of the link in question.
More concretely, we want to know that every component of the link contributes to the link having this circumnavigation property.

\begin{defn}
We say that a link in the solid torus is \emph{atomic} if it circumnavigates the torus and if the removal of any component means that it no longer does so.

That is to say, the meridian of the torus is non\hyp{}trivial in the fundamental group of the complement of the link, but if any component is removed then it becomes trivial.

The \emph{pure atomic link operad}, \(\m{P A L}\), is the suboperad of the pure link operad wherein for \(n \ge 2\) the morphisms consist of atomic links.
\end{defn}

\tikzsetnextfilename{twisthopf}

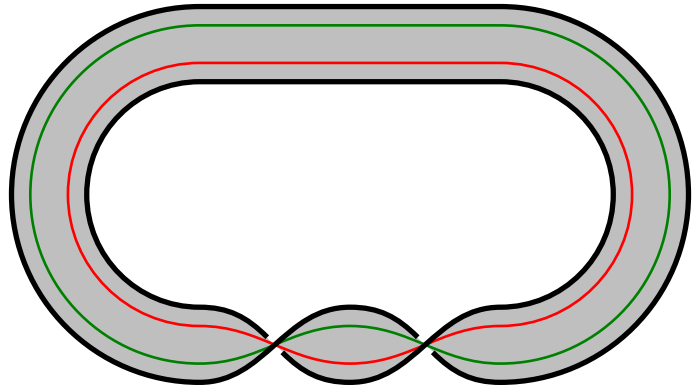
\begin{figure}
\centering
\begin{tikzpicture}
\draw[torus, line width=1cm] (0,-2) arc (270:90:2) -- ++(4,0) arc(90:-90:2);
\fill[torus,even odd rule] (0,-2) ++(0,.5) to[out=0,in=180] ++(2,-1) to[out=0,in=180] ++(2,1) arc (-90:90:1.5) -- ++(-4,0) arc(90:270:1.5) ++(0,-1) to[out=0,in=180] ++(2,1) to[out=0,in=180] ++(2,-1) arc(-90:90:2.5) -- ++(-4,0) arc(90:270:2.5);
\draw[string,ring1] (0,-2) ++(0,.25) to[out=0,in=180] ++(2,-.5) to[out=0,in=180] ++(2,.5) arc (-90:90:1.75) -- ++(-4,0) arc(90:270:1.75);
\draw[string,ring2] (0,-2) ++(0,-.25) to[out=0,in=180] ++(2,.5) to[out=0,in=180] ++(2,-.5) arc(-90:90:2.25) -- ++(-4,0) arc(90:270:2.25);
\draw[line width=2pt] (0,-2) ++(0,.5) to[out=0,in=135] ++(.9,-.4) ++(.2,-.2) to[out=-45,in=180] ++(.9,-.4) to[out=0,in=180] ++(2,1) arc (-90:90:1.5) -- ++(-4,0) arc(90:270:1.5) ++(0,-1) to[out=0,in=180] ++(2,1) to[out=0,in=135] ++(.9,-.4) ++(.2,-.2) to[out=-45,in=180] ++(.9,-.4) arc(-90:90:2.5) -- ++(-4,0) arc(90:270:2.5);
\end{tikzpicture}
\caption{The Hopf link inside a twisted torus}
\autolabel
\end{figure}

\section{Invariants}

In this section we shall investigate which of the currently available knot and link invariants detect the possible factorisations of a knot or link.

\subsection{Skein Relations}

There are various knot and link invariants that are based on the Skein relations, or variants thereof.
Examples include the Homfly\hyp{}PT polynomial, the Jones polynomial, and the Kauffman bracket.
The basic idea behind these invariants is to take a crossing, cut the strands involved in the crossing, and splice them together again in various ways.
If the strands are oriented then there are three possible outcomes (one of which is the original crossing), shown in Figure~\ref{fig:orskein}.
If the strands are not oriented then there are four possible outcomes, shown in Figure~\ref{fig:uorskein}.

\tikzsetnextfilename{orskein}
\begin{figure}
\centering
\begin{tikzpicture}[every path/.style={>=triangle 60,thick knot,double=Red,
    decoration={markings,mark=at position 1 with {\arrow[double=none,Red,line width=1pt]{>}}},
    shorten >=1pt}]
\draw[postaction={decorate}] (2,0) -- (0,2);
\draw[postaction={decorate}] (0,0) -- (2,2);
\begin{scope}[xshift=3cm]
\draw[postaction={decorate}] (0,0) -- (2,2);
\draw[postaction={decorate}] (2,0) -- (0,2);
\begin{scope}[xshift=3cm]
\draw[postaction={decorate}] (2,0) to[bend left] (2,2);
\draw[postaction={decorate}] (0,0) to[bend right] (0,2);
\end{scope}
\end{scope}
\end{tikzpicture}
\caption{The Oriented Skeins}
\autolabel
\end{figure}
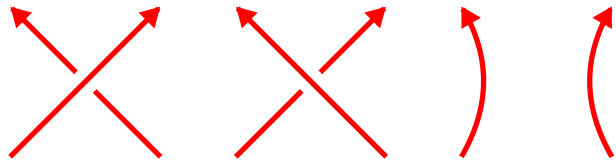

\tikzsetnextfilename{uorskein}
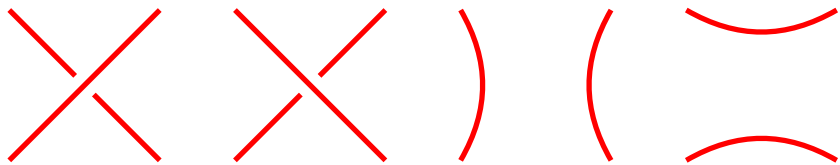
\begin{figure}
\centering
\begin{tikzpicture}[every path/.style={thick knot,double=Red}]
\draw (2,0) -- (0,2);
\draw (0,0) -- (2,2);
\begin{scope}[xshift=3cm]
\draw (0,0) -- (2,2);
\draw (2,0) -- (0,2);
\begin{scope}[xshift=3cm]
\draw (2,0) to[bend left] (2,2);
\draw (0,0) to[bend right] (0,2);
\begin{scope}[xshift=3cm]
\draw (0,0) to[bend left] (2,0);
\draw (0,2) to[bend right] (2,2);
\end{scope}
\end{scope}
\end{scope}
\end{tikzpicture}
\caption{The Unoriented Skeins}
\autolabel
\end{figure}

One purpose of defining the level of a link as we have done is to be able to study the link at a particular level.
In so doing, we ignore all structure coming from more refined levels.
In particular, when examining a ``strand'' at one level, we cannot know its finer structure.
Thus we cannot splice it to anything other than itself, which means that in the skein relations, we cannot allow contributions from diagrams where a strand is spliced to a different strand.
This simply leaves us with exchanging crossings as the only allowed operation.
From this, it is a reasonable conjecture that the various invariants that using the skein relations will not detect the level of a link as they must work always with the full link and cannot truncate to a particular level.

To see this concretely, let us consider the Hopf ring of length \(2\) as in the left in Figure~\ref{fig:hopfring2}.
Let us also compose this with itself in both components, producing the Hopf ring of level \(2\) on the right in Figure~\ref{fig:hopfring2}.

\tikzsetnextfilename{hopfring2}

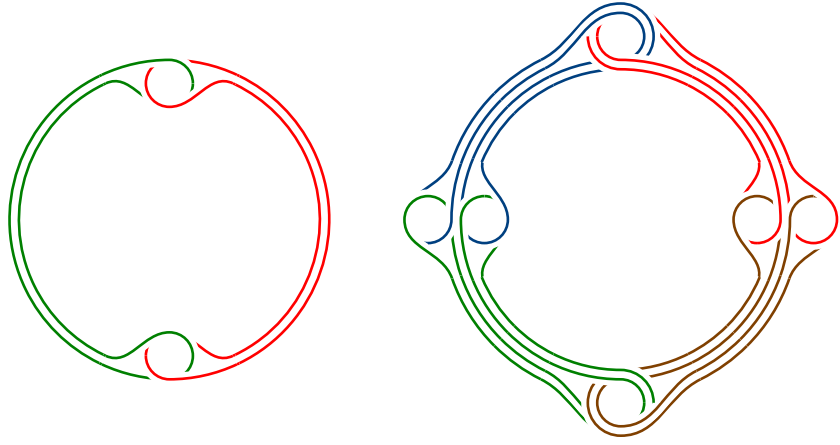
\begin{figure}
\centering
\begin{tikzpicture}
\setbrstep{.125}
\hopfring{2}{2}
\begin{scope}[xshift=6cm]
\hopftwo{2}{2}{2}
\end{scope}
\end{tikzpicture}
\caption{Hopf ring with two components and the same composed with itself}
\autolabel
\end{figure}

Computing the HOMFLY\emhyp{}PT polynomials of these, using the \texttt{homfly} program, we obtain for the Hopf ring:
\[
- m^{-1} l^{-5} - m^{-1} l^{-3} + m l^{-3} - m l^{-1}
\]
and for the level two Hopf ring:
\begin{dmath*}
- m^{-3} l^{-11} - 3 m^{-3} l^{-9} - 3 m^{-3} l^{-7} - m^{-3} l^{-5} + 2 m^{-1} l^{-9} + 2 m^{-1} l^{-7} - 2 m^{-1} l^{-5} - 2 m^{-1} l^{-3} - m l^{-7} + m l^{-5} + m l^{-3} - m l^{-1} + 6 m^{3} l^{-7} + 16 m^{3} l^{-5} + 12 m^{3} l^{-3} - 2 m^{3} l - 11 m^{5} l^{-7} - 35 m^{5} l^{-5} - 26 m^{5} l^{-3} + m^{5} l - m^{5} l^{3} + 6 m^{7} l^{-7} + 28 m^{7} l^{-5} + 22 m^{7} l^{-3} - m^{9} l^{-7} - 9 m^{9} l^{-5} - 8 m^{9} l^{-3} + m^{11} l^{-5} + m^{11} l^{-3}
\end{dmath*}

There is almost too much information there to see if there is a pattern.
So let us consider instead their Kauffman brackets:
\[
- A^{10} + A^{6} - A^{2} - A^{-6}
\]
and
\begin{dmath*}
- A^{46} + 5 A^{42} - 11 A^{38} + 14 A^{34} - 12 A^{30} + 9 A^{26} - 7 A^{22} + 2 A^{18} + 2 A^{14} - 7 A^{10} + 8 A^{6} - 9 A^{2} + 6 A^{-2} - 9 A^{-6} + 7 A^{-10} - 8 A^{-14} + 4 A^{-18} - A^{-22}.
\end{dmath*}

Although it is nigh\hyp{}on impossible to show that there is no relationship between the two, as evidence for that fact let us examine what happens to the second when we quotient by the ideal generated by the first.
We can find a representative of the equivalence class in the quotient with the smallest span, smallest degree, and only positive degree terms.
This is:
\[
- A^{12} - 5 A^{8} + 2 A^{4} + 2.
\]

\subsection{Complements}

Unlike invariants defined using the skein relations, invariants defined using the link complement stand a better chance of detecting or being compatible with the notion of level.
This is because a factorisation of the link defines a filtration on the complement.
Let \(L\) be a link represented by an embedding \(f \colon \coprod S^1 \times D^2 \to S^1 \times D^2\) with a factorisation
\begin{equation}
\label{eq:factor}
  f = f_k \circ (g_1, \dotsc, g_k).
\end{equation}
Then \(f_k\) is also a link and the image of \(f_k\) contains the image of \(f\).
Hence the complement of \(f\) contains the complement of \(f_k\).

A full factorisation in the operad will not, in general, be a linear factorisation.
If we continue the factorisation in \eqref{eq:factor} then we next look at the \(g_j\) and the factorisations of these are independent.
Thus a full factorisation will look like a tree, and for each rooted subtree we obtain a link complement.
Inclusion of subtrees maps to inclusion of spaces.

Thus if we are given a functor on the category of link complements (in the solid torus) with inclusions as morphisms, the factorisation of the link defines a diagram in the target category.

We can use this in two ways: to detect factorisations and to study them.

\subsection{The Fundamental Group}
\label{sec:fungrp}

As the fundamental group is a covariant functor, when applied to the factorisation of a link we obtain a homomorphism from the fundamental group of the complement of a subtree to that of the whole link.

Let us consider a factorisation as in Equation~\eqref{eq:factor}:
\[
  f = f_k \circ (g_1, \dotsc, g_k).
\]
Let us write \(C_f\), \(C_{f_k}\), and \(C_{g_i}\) for the obvious complements.
As our link complements are path connected, the exact location of the basepoint is not important, so let us choose it to be on the exterior of the torus.
We have an inclusion \(C_{f_k} \subseteq C_f\) and therefore a group homomorphism \(\pi_1(C_{f_k}) \to \pi_1(C_f)\).
In this case the ambient tori are the same so there is no basepoint ambiguity.

Now consider the operation of composition with \(g_1\).
This pastes in a new torus (with some bits missing) into the gap left by the complement of the first component in \(f_k\).
Let us write \(f_{k,1}\) for this composition.
That is,
\[
  f_{k,1} = f_k \circ (g_1, 1, \dotsc, 1).
\]
The complement of \(f_{k,1}\) is then the union of the complement of \(f_k\) and the complement of \(g_1\), with the torus surrounding \(g_1\) glued in to the gap left by first component of \(f_k\).

As we have assumed that our links do not touch the boundary of the torus we can adjust the pieces by homotopies in this union so that each part is an open set in \(C_{f_{k,1}}\).
The intersection of the two pieces is homotopy equivalent to an unfilled torus, \(S^1 \times S^1\).
This, then, fulfils the criteria for the Van Kampen Theorem and we have that the fundamental group of \(C_{f_{k,1}}\) is the free product of the fundamental groups of \(C_{f_k}\) and \(C_{g_1}\) modulo the identification of the meridian and longitude of the surrounding torus in the two groups.

We proceed by induction and conclude that the following is true.

\begin{theorem}
The fundamental group of \(C_f\) is formed by taking the free product of the fundamental groups of \(C_{f_k}\) and the \(C_{g_i}\) and identifying the meridians and longitudes of the surrounding tori.
\noproof
\end{theorem}

If we take the Wirtinger presentation of the fundamental groups then in \(\pi_1(C_{f_k})\) the meridians are actually generators.
On the other hand, we have assumed that our factorisation is in the pure link operad and so in \(\pi_1(C_{g_i})\) the meridians are non\hyp{}trivial.
We deduce, therefore, that the meridians are non\hyp{}trivial in \(\pi_1(C_f)\).

Now the meridian corresponding to, say, \(g_1\) has the property that removing any of its subcomponents trivialises it, but removing any other components does not.
Therefore if we know the meridian corresponding to \(g_1\) but do not know which components are involved in \(g_1\), we can determine this by looking for those components which have the property that upon their removal the meridian collapses.

This suggests a possible strategy for finding a factorisation: look in \(\pi_1(C)\) for elements that could be meridians of embedded tori and then determine those components that collapse it.
In a genuine factorisation these meridians would form a tree structure, where the ordering was that \(\alpha \preceq \beta\) if every component that collapses \(\alpha\) also collapses \(\beta\).
Note that the meridian of the outermost torus is the topmost element of this ordering, and if we include the elements that correspond to components these are the minimal elements.

Of interest also are the longitudes.
In a torus then the longitude and meridian commute.
So when we have a potential meridian we need to look for its possible longitudes by looking in its centraliser subgroup.

\subsection{Homology and Cohomology}

In a similar fashion, the homology and cohomology functors applied to the link complement will be filtered according to the factorisation structure of the link.
For homology, this will be a filtration by subgroups whilst for cohomology this will be a projective filtration, though we could convert this to a question of ideals in the cohomology ring.

To use homology and cohomology to detect the filtration we will need to use more than just the groups (rings) themselves.
The homology groups of a link complement are determined solely by the number of components via a Mayer\enhyp{}Vietoris argument.
Indeed, let \(f \colon \coprod_n S^1 \times D^2 \to S^1 \times D^2\) be an \(n\)\enhyp{}component link and let \(C_f\) be its complement.
Then \(S^1 \times D^2\) is the union of \(C_f\) and \(\coprod_n S^1 \times D^2\).
The intersection is \(\coprod_n S^1 \times S^1\).
The Mayer\enhyp{}Vietoris sequence is thus:
\begin{multline*}
  \to H_{k+1}(S^1 \times D^2) \to H_k(\coprod_n S^1 \times S^1) \to H_k(\coprod_n S^1 \times D^2) \oplus H_k(C_f) \\
 \to H_k(S^1 \times D^2) \to H_{k-1}(\coprod_n S^1 \times S^1)
\end{multline*}
The first non\hyp{}trivial term is with \(k = 2\) where we have:
\[
  0 \to \Z^n \to H_2(C_f) \to 0
\]
and thus \(H_2(C_f) \cong \Z^n\).
At the bottom end we note that the inclusion \(S^1 \times S^1 \to S^1 \times D^2\) induces an isomorphism on \(H_0\), whence the corresponding term in the Mayer\enhyp{}Vietoris sequence is an injection.
Thus for \(k=1\) we have:
\[
  0 \to \Z^{2n} \oplus \Z^n \oplus H_1(C_f) \to \Z \to 0.
\]
The generators of \(\Z^{2n}\) are the meridians and longitudes.
The map to the \(\Z^n\) takes care of the longitudes and the meridians map to the \(H_1(C_f)\) factor.
The final factor of \(\Z\) is the longitude of the outer torus.
This \(H_1(C_f) \cong \Z^{n+1}\).
This fits with \(H_1(C_f)\) being the abelianisation of the fundamental group.

This, therefore, contains no information about their filtration.
Moreover the key elements in the fundamental group, the meridians, are all zero in homology since they are always commutators (this is because the link is pure).

Turning to cohomology we can use the product as a more refined tool.
However, this \emhyp{} and the Massey products \emhyp{} only detects the type of link between two components, not whether or not there is a containing torus.

We would therefore need to use some further structure on the chains or cochains where it was possible to detect the factorisation structure prior to taking homology.

\subsection{The Relative Link Group}

One avenue for detecting factorisations is to look for appropriate subgroups of the fundamental group of the link complement.
Although this is a very natural place to look, the complexity of the presentations works against it.

The natural way to look for a factorisation is to consider a family of components and ask if that family can be put inside a torus without disturbing the rest of the link.
When doing this, the intricacies of the rest of the link are irrelevant, all that matters is how the other components wrap around the family.
Thus after deciding on a family to consider, we regard the other components simply as elements in the fundamental group of the complement of the sublink.
If there is a surrounding torus, these elements will be in the subgroup generated by its meridian.
Therefore we look for an element in the fundamental group of the complement of the sublink that generates a subgroup containing the elements defined by the other components of the main link.

There are some additional properties that this element will have to satisfy which come from the restrictions that we have placed on our factorisations.
These properties can be summarised as saying that the sublink is a minimal sublink which carries that particular element.
That is to say, in the fundamental group associated to any smaller sublink, the element that we are looking at must vanish.

It is also important to note that as the sublink cannot circumnavigate the outer torus, we can view the link as being in \(\R^3\) instead of the torus.
That is, we can remove the meridian from the outermost torus from the fundamental group.

\medskip

As an example, we consider two links both with three components.
The first has the following fundamental group.
The generators are labelled such that those with main symbol \(a\) correspond to the background torus whilst the other generators are labelled according to the corresponding component of the link.
The generators are:
\begin{dmath}
\left\{ a_1, a_2, a_3, a_4,  b_1, b_2, b_3, b_4,  c_1, c_2, c_3, c_4,  d_1, d_2, d_3, d_4,  d_5, d_6, d_7, d_8 \right\}
\end{dmath}
The relations are:
\begin{dmath}[style={\eqmargin=0pt}]
d_1 c_4 d_1^{-1} c_1^{-1},\;\allowbreak
d_1 c_3^{-1} d_2^{-1} c_3,\;\allowbreak
d_2 a_2 d_2^{-1} a_3^{-1},\;\allowbreak
d_2 a_3^{-1} d_3^{-1} a_3,\;\allowbreak
d_3 c_2 d_3^{-1} c_3^{-1},\;\allowbreak
d_3 c_1^{-1} d_4^{-1} c_1,\;\allowbreak
d_4 b_3 d_4^{-1} b_4^{-1},\;\allowbreak
d_4 b_3^{-1} d_5^{-1} b_3,\;\allowbreak
d_5 c_1 d_6^{-1} c_1^{-1},\;\allowbreak
d_6 c_2^{-1} d_6^{-1} c_1,\;\allowbreak
d_6 a_3 d_7^{-1} a_3^{-1},\;\allowbreak
d_7 a_2^{-1} d_7^{-1} a_1,\;\allowbreak
d_7 c_3 d_8^{-1} c_3^{-1},\;\allowbreak
d_8 c_4^{-1} d_8^{-1} c_3,\;\allowbreak
d_8 b_1 d_8^{-1} b_2^{-1},\;\allowbreak
d_8 b_1^{-1} d_1^{-1} b_1,\;\allowbreak
b_2 a_4 b_2^{-1} a_1^{-1},\;\allowbreak
b_2 a_3^{-1} b_3^{-1} a_3,\;\allowbreak
b_4 a_3 b_1^{-1} a_3^{-1},\;\allowbreak
b_1 a_4^{-1} b_1^{-1} a_3.
\end{dmath}
We also record the elements in the fundamental group corresponding to each component (and the background torus), where we start each longitudinal loop at the strand labelled \({}_1\):
\begin{align*}
a &: d_7 d_2^{-1} b_1 b_2^{-1}, \\
b &: d_8^{-1} a_3^{-1} d_4^{-1} a_3, \\
c &: d_6 d_3^{-1} d_8 d_1^{-1}, \\
d &: c_3^{-1} a_3^{-1} c_1^{-1} b_3^{-1} c_1 a_3 c_3 b_1^{-1}.
\end{align*}

We now consider sublinks.
As our sublinks must be proper sublinks, and cannot comprise just one component, in this case there are three to consider and each is given by deleting one component.

We can immediately see that there is only one sublink to consider.
If we delete the \(b\)\enhyp{}component, then the element of the fundamental group of the \(c{-}d\)\enhyp{}sublink represented by \(b\) is not trivial when we pass further to the sublink containing just the \(d\)\enhyp{}component: this can be seen by looking at the exponents.
A similar thing happens when the \(d\) component is removed.
We are thus left with removing the \(c\)\enhyp{}component.

From the relations this produces the following equalities: \(d_1 = d_2\), \(d_3 = d_4\), \(d_5 = d_6\), \(d_7 = d_8\).
The longitude along \(c\) is \(d_6 d_3^{-1} d_8 d_1^{-1}\).
Clearly, if we remove \(d\) then this collapses.
We require it also to collapse if we remove \(b\).
Removing \(b\) identifies \(d_4 = d_5\) whence \(d_3 = d_6\) and the first pair in \(c\) cancel.
Removing \(b\) also identifies \(d_8 = d_1\), whence \(c\) collapses.

Thus we have a candidate for a sublink involving just \(b\) and \(d\). To find the corresponding meridian we need to look at \(d_6 d_3^{-1} d_8 d_1^{-1}\) and see if we can simplify it in the fundamental group of the complement of \(b\) and \(d\).
Since we are now in a sublink, we can identify the \(a_i\) with the identity.
This identifies \(d_2 = d_3\), \(d_6 = d_7\), \(b_2 = b_3\), and \(b_1 = b_4\).
We therefore have \(d_1 = d_2 = d_3 = d_4\) and \(d_5 = d_6 = d_7 = d_8\).
The longitude of \(c\) is thus \(d_5 d_1^{-1} d_5 d_1^{-1}\) which is the square of \(d_5 d_1^{-1}\).
Note that removing \(b\) still collapses this element so this still satisfies the requirement for a meridian.
It is clearly minimal.

Thus our suggestion for a sublink is the link comprising components \(b\) and \(d\), with the meridian of the surrounding torus being \(d_5 d_1^{-1}\) (or equivalent thereof).

Now the element \(d_5 d_1^{-1}\) here is in the complement of the components \(b\) and \(d\).
In the complement of the full link we might not want to use this precise element but one that becomes it upon removal of \(c\).
The key property that we need is that the meridian collapses upon removal of just \(d\) or \(b\).
Now it is clear that \(d_5 d_1^{-1}\) collapses upon removal of \(d\), but not clear what happens if we remove just \(b\).
Using the relations we see that:
\begin{dmath*}
  d_5 d_1^{-1} = (c_1 d_6 c_1^{-1}) d_1^{-1} = (c_1 a_3 d_7 a_3^{-1} c_1^{-1}) d_1^{-1} = (c_1 a_3 c_3 d_8 c_3^{-1} a_3^{-1} c_1^{-1}) d_1^{-1} = (c_1 a_3 c_3 b_1^{-1} d_1 b_1 c_3^{-1} a_3^{-1} c_1^{-1}) d_1^{-1}
\end{dmath*}
And this does not collapse if we remove just \(b\) since that does not allow us to bring \(d_1\) and \(d_1^{-1}\) together.
But in place of \(d_5\) we could use one of \(d_6\), \(d_7\), or \(d_8\) and it is clear from the above that \(d_8\) would work since \(d_8 d_1^{-1} = b_1^{-1} d_1 b_1 d_1^{-1}\).
This, then, is our proposal for a meridian.

Now if we take just the complement of the \(b\) and \(d\) components then we have generators \(b_1\), \(b_2\), \(d_1\), and \(d_8\).
The relations that just involve \(b\) and \(d\) are (adjusted for the identified generators):
\begin{dmath}[style={\eqmargin=0pt}]
d_1 b_2 d_1^{-1} b_1^{-1},\;\allowbreak
d_1 b_2^{-1} d_8^{-1} b_2,\;\allowbreak
d_8 b_1 d_8^{-1} b_2^{-1},\;\allowbreak
d_8 b_1^{-1} d_1^{-1} b_1.
\end{dmath}
However, this is simply the complement of \(b\) and \(d\) in \(\R^3\).
We need to put them in a torus with meridian \(d_8 d_1^{-1}\).
This means that we need to introduce a new component, say \(e\), which goes around \(d_8\) and \(d_1\) (in the appropriate directions).
This splits \(d_8\) and \(d_1\) so that we have new generators \(d_8'\), \(d_1'\), \(e_1\), and \(e_2\) with relations:
\begin{dmath}[style={\eqmargin=0pt}]
e_1^{-1} d_8 e_2 d_8^{-1},\;\allowbreak
e_2^{-1} d_1^{-1} e_1 d_1,\;\allowbreak
e_1 {d_8'}^{-1} e_1^{-1} d_8,\;\allowbreak
e_1 d_1' e_1^{-1} d_1^{-1}.
\end{dmath}
Comparing these with the original relations, it would make sense to identify \(d_8' = d_7\) and \(d_1' = d_2\).
Then \(e\) looks a little like part of \(c\) (as it should since \(c\) should entwine around \(b\) and \(d\) only as multiplies of \(e\)).
The longitude in this torus is simply \(e_1\).

Back in the larger group, we want to replace the \(b\) and \(d\) components by a single new component, say \(f\), such that \(d_8 d_1^{-1}\) is one of the new generators, say \(f_1\).
This will interact with the \(c\) and \(a\) components in some fashion.
Looking at the relations, we see that:
\[
  c_3 d_8 d_1^{-1} c_3^{-1} = d_7 d_2^{-1}
\]
so if we write \(f_2 = d_7 d_2^{-1}\) then we have \(c_3 f_1 c_3^{-1} f_2^{-1} = 1\).
Continuing in this vein we see that we have generators \(f_1, f_2, f_3, f_4\).

\medskip

Now let us consider the second example, also with three components.
The generators for this link are (with the same convention as above):
\[
  \left\{a_1, a_2, b_1, b_2, b_3, b_4, c_1, c_2, d_1, d_2\right\}
\]
The relations are:
\begin{dmath}[style={\eqmargin=0pt}]
a_1^{-1} b_1^{-1} a_2 b_1, \;\allowbreak
a_1 b_4^{-1} a_1^{-1} b_1, \;\allowbreak
a_1 b_3 a_1^{-1} b_2^{-1}, \;\allowbreak
a_2^{-1} b_2 a_1 b_2^{-1}, \;\allowbreak
c_2^{-1} b_1 c_1 b_1^{-1}, \;\allowbreak
b_2 c_2^{-1} b_1^{-1} c_1, \;\allowbreak
d_2^{-1} c_1 d_1 c_1^{-1}, \;\allowbreak
c_2 d_2^{-1} c_1^{-1} d_1, \;\allowbreak
d_1^{-1} b_3 d_2 b_3^{-1}, \;\allowbreak
b_4 d_1^{-1} b_3^{-1} d_2
\end{dmath}
The longitudinal paths are:
\begin{align*}
a &: b_1^{-1} b_2, \\
b &: a_1^{-1} d_1 a_1 c_2, \\
c &: b_1 d_2, \\
d &: c_1 b_3.
\end{align*}

In this case then there are no sublinks.
As before, the potential sublinks are those where we remove one component.
Let us remove the \(d\) component.
Then the corresponding longitude is \(c_1 b_3\) in the complement of \(b\) and \(c\).
But we want this to be trivial if we remove a further component, and this is not the case.
A similar thing happens for the other components.
Therefore this link has no non\hyp{}trivial factorisation. 

\bigskip

The two links under consideration were the two in Figures~\ref{fig:fundgrphopf} and~\ref{fig:fundgrphopfring}.

\tikzsetnextfilename{fundgrphopf}
%\tikzset{external/export next=false}
\begin{figure}
\centering
\begin{tikzpicture}[every path/.style={knot=Red},every node/.style={black}]
\pgfmathsetmacro{\brover}{.25}
\pgfmathsetmacro{\bradj}{.1}
\pgfmathsetmacro{\brlen}{1}
\pgfmathsetmacro{\brsep}{.25}
\foreach \bri/\brdir/\brnud/\brv in {
  1/1/1/-1,
  2/1/0/-1,
  3/-1/2/-1,
  4/1/0/1,
  5/1/1/1,
  6/-1/2/1
} {
\coordinate (u\bri) at (\brover * \brdir * \brnud,-\bri * \brsep + \brv * \bradj);
\draw (u\bri) -- ++(0,-\brv * \bradj) -- +(-\brdir * \brlen - \brover * \brdir * \brnud,0) coordinate (u\bri e);
}
\draw (u1) -- (u5);
\draw (u2) -- (u4);
\draw (u3) -- (u6);

\begin{scope}[yshift=-2cm]
\foreach \bri/\brdir/\brnud/\brv in {
  1/1/1/-1,
  2/1/0/-1,
  3/-1/2/-1,
  4/1/0/1,
  5/1/1/1,
  6/-1/2/1
} {
\coordinate (l\bri) at (\brover * \brdir * \brnud,-\bri * \brsep + \brv * \bradj);
\draw (l\bri) -- ++(0,-\brv * \bradj) -- +(-\brdir * \brlen - \brover * \brdir * \brnud,0) coordinate (l\bri e);
}
\draw (l1) -- (l5);
\draw (l2) -- (l4);
\draw (l3) -- (l6);
\end{scope}
\draw (u3e) -- +(4*\brsep,0) |- (l6e);
\draw (u6e) -- +(\brsep,0) |- (l3e);
\draw ($(u5e)!.5!(l1e)$) -- ++(0,\brsep) -- ++(-9*\brsep,0) -- ++(0,-\brsep);
\draw (u5e) -- +(-\brsep,0) |- (l1e);
\draw (u4e) -- +(-2*\brsep,0) |- (l2e);
\draw (u1e) -- ++(-2*\brsep,0) -- ++(0,-.5*\brsep) -- ++(-2.5*\brsep,0);
\draw (u1e) ++(-8*\brsep,-\brsep) coordinate (ul1e) -- ++(0,\brsep) -- ++(5*\brsep,0) -- ++(0,-\brsep) -- ++(-4*\brsep,0) -- ++(0,-\brsep) coordinate (ul2e);
\draw (u2e) -- ++(-\brsep,0) -- ++(0,-.5*\brsep) -- ++(-4*\brsep,0) -- ++(0,\brsep) -- ++(.5*\brsep,0);
\draw (l4e) -- ++(-2*\brsep,0) -- ++(0,-.5*\brsep) -- ++(-2.5*\brsep,0);
\draw (l4e) ++(-7*\brsep,\brsep) coordinate (ll4e) -- ++(0,-\brsep) -- ++(4*\brsep,0) -- ++(0,-\brsep) -- ++(-5*\brsep,0) -- ++(0,\brsep) coordinate (ll5e);
\draw (l5e) -- ++(-\brsep,0) -- ++(0,-.5*\brsep) -- ++(-4*\brsep,0) -- ++(0,\brsep) -- ++(.5*\brsep,0);
\draw (ul1e) -- (ll5e);
\draw (ul2e) -- (ll4e);
\draw ($(u5e)!.5!(l1e)$) -- ++(0,-\brsep) -- ++(-9*\brsep,0) -- ++(0,\brsep);
%\fill[black] (0,0) circle (2pt);
%\draw[opacity=.5,gray,double=none,step=2*\brsep,thin] (current bounding box.south west) grid (current bounding box.north east);
\node at (-12*\brsep,0) {\(b_1\)};
\node at (-10*\brsep,-3*\brsep) {\(b_2\)};
\node at (-10*\brsep,-11*\brsep) {\(b_3\)};
\node at (-12*\brsep,-14*\brsep) {\(b_4\)};
\node at (-9*\brsep,-5*\brsep) {\(a_1\)};
\node at (-7*\brsep,-7*\brsep) {\(a_2\)};
\node at (-9*\brsep,-9*\brsep) {\(a_3\)};
\node at (-13*\brsep,-5*\brsep) {\(a_4\)};
\node at (6*\brsep,-2*\brsep) {\(c_1\)};
\node at (2*\brsep,-10*\brsep) {\(c_2\)};
\node at (4*\brsep,-8*\brsep) {\(c_3\)};
\node at (2*\brsep,-4*\brsep) {\(c_4\)};
\node at (-4*\brsep,0) {\(d_1\)};
\node at (-3*\brsep,-6*\brsep) {\(d_2\)};
\node at (-\brsep,-8*\brsep) {\(d_3\)};
\node at (-4*\brsep,-14*\brsep) {\(d_4\)};
\node at (-4*\brsep,-11*\brsep) {\(d_5\)};
\node at (-7*\brsep,-10*\brsep) {\(d_6\)};
\node at (-7*\brsep,-4*\brsep) {\(d_7\)};
\node at (-3*\brsep,-3*\brsep) {\(d_8\)};
\end{tikzpicture}
\caption{Hopf ring with sub\hyp{}Hopf ring}
\autolabel
\end{figure}
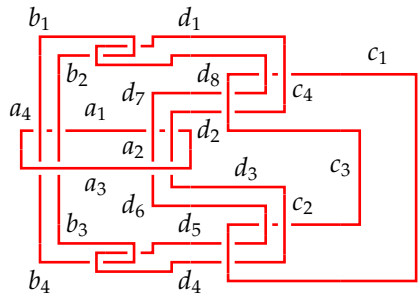

\tikzsetnextfilename{fundgrphopfring}
%\tikzset{external/export next=false}
\begin{figure}
\centering
\begin{tikzpicture}[every path/.style={knot=Red},every node/.style={black}]
\pgfmathsetmacro{\brover}{.25}
\pgfmathsetmacro{\bradj}{.1}
\pgfmathsetmacro{\brlen}{1}
\pgfmathsetmacro{\brsep}{.25}
\def\brprefix{ul}
\foreach \bri/\brdir/\brnud/\brv in {
  1/1/1/-1,
  2/-1/1/-1,
  3/1/1/1,
  4/-1/1/1%
} {
\coordinate (\brprefix\bri) at (\brover * \brdir * \brnud,-\bri * \brsep + \brv * \bradj);
\draw (\brprefix\bri) -- ++(0,-\brv * \bradj) -- +(-\brdir * \brlen - \brover * \brdir * \brnud,0) coordinate (\brprefix\bri e);
}
\draw (\brprefix1) -- (\brprefix3);
\draw (\brprefix2) -- (\brprefix4);

\begin{scope}[yshift=-2cm,xshift=1cm]
\def\brprefix{l}
\foreach \bri/\brdir/\brnud/\brv in {
  1/1/1/-1,
  2/-1/1/-1,
  3/1/1/1,
  4/-1/1/1%
} {
\coordinate (\brprefix\bri) at (\brover * \brdir * \brnud,-\bri * \brsep + \brv * \bradj);
\draw (\brprefix\bri) -- ++(0,-\brv * \bradj) -- +(-\brdir * \brlen - \brover * \brdir * \brnud,0) coordinate (\brprefix\bri e);
}
\draw (\brprefix1) -- (\brprefix3);
\draw (\brprefix2) -- (\brprefix4);
\end{scope}

\begin{scope}[xshift=2cm,yshift=-\brsep cm]
\def\brprefix{ur}
\foreach \bri/\brdir/\brnud/\brv in {
  1/1/1/-1,
  2/-1/1/-1,
  3/1/1/1,
  4/-1/1/1%
} {
\coordinate (\brprefix\bri) at (\brover * \brdir * \brnud,-\bri * \brsep + \brv * \bradj);
\draw (\brprefix\bri) -- ++(0,-\brv * \bradj) -- +(-\brdir * \brlen - \brover * \brdir * \brnud,0) coordinate (\brprefix\bri e);
}
\draw (\brprefix1) -- (\brprefix3);
\draw (\brprefix2) -- (\brprefix4);
\end{scope}
\draw (ul3e) ++(\brsep,-1cm) ++(0,\brsep) -- ++(0,\brsep) -- ++(-6*\brsep,0) -- ++(0,-\brsep);
\draw (ul1e) -- +(-3*\brsep,0) |- (l3e);
\draw (ul3e) -- +(-\brsep,0) |- (l1e);
\draw (ul2e) -- (ur1e);
\draw (ul4e) -- (ur3e);
\draw (ur2e) -- +(3*\brsep,0) |- (l4e);
\draw (ur4e) -- +(\brsep,0) |- (l2e);
\draw (ul3e) ++(\brsep,-1cm) ++(0,\brsep) -- ++(0,-\brsep) -- ++(-6*\brsep,0) -- ++(0,\brsep);
%\fill[black] (0,0) circle (2pt);
%\draw[opacity=.5,gray,double=none,step=2*\brsep,thin] (current bounding box.south west) grid (current bounding box.north east);
\node at (-4*\brsep,0) {\(b_1\)};
\node at (-4*\brsep,-4*\brsep) {\(b_2\)};
\node at (0,-8*\brsep) {\(b_3\)};
\node at (0,-12*\brsep) {\(b_4\)};
\node at (-10*\brsep,-6*\brsep) {\(a_1\)};
\node at (-6*\brsep,-6*\brsep) {\(a_2\)};
%\node at (-9*\brsep,-9*\brsep) {\(a_3\)};
%\node at (-13*\brsep,-5*\brsep) {\(a_4\)};
\node at (5*\brsep,-\brsep) {\(c_1\)};
\node at (5*\brsep,-5*\brsep) {\(c_2\)};
%\node at (4*\brsep,-8*\brsep) {\(c_3\)};
%\node at (2*\brsep,-4*\brsep) {\(c_4\)};
\node at (12*\brsep,-7*\brsep) {\(d_2\)};
\node at (16*\brsep,-7*\brsep) {\(d_1\)};
%\node at (-\brsep,-8*\brsep) {\(d_3\)};
%\node at (-4*\brsep,-14*\brsep) {\(d_4\)};
%\node at (-4*\brsep,-11*\brsep) {\(d_5\)};
%\node at (-7*\brsep,-10*\brsep) {\(d_6\)};
%\node at (-7*\brsep,-4*\brsep) {\(d_7\)};
%\node at (-3*\brsep,-3*\brsep) {\(d_8\)};
\end{tikzpicture}
\caption{Hopf ring with three components}
\autolabel
\end{figure}
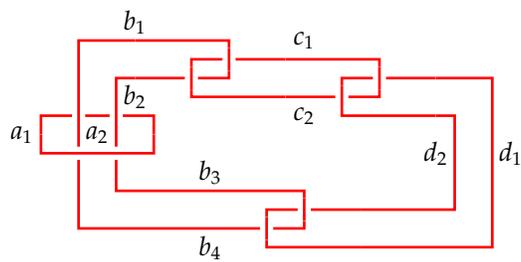

\chapter{Kauffman Computations}

\def\pgfmathsetmacro#1#2{}

\mymaketitle

\section{The Rules of the Kauffman Bracket}

The rules of the Kauffman bracket are simple.

\begin{enumerate}
%\tikzset{external/export=false}
\item \(\langle \tikz[baseline=-.8ex] \draw[knot,double=Red] (0,0) circle (1ex); \rangle = 1\)
\item \(\langle \tikz[baseline=-.8ex] \node[knot under cross,knot,draw,double=Red] {}; \rangle = A \langle \tikz[baseline=-.8ex] \node[knot vert,knot,draw,double=Red] {}; \rangle + A^{-1} \langle \tikz[baseline=-.8ex] \node[knot horiz,knot,draw,double=Red] {}; \rangle\)
\item \(\langle  \tikz[baseline=-.8ex] \draw[knot,double=Red] (0,0) circle (1ex); \cup L \rangle = (-A^2 - A^{-2}) \langle L \rangle\)
\end{enumerate}

\section{The Simplification Technique}
\label{sec:simple}

To compute the Kauffman bracket of our links, we shall use a simplification technique.
The identity for the Kauffman bracket replaces a diagram with a crossing by two diagrams without that crossing.
It therefore reduces the number of crossings by \(1\) whilst doubling the number of diagrams that need to be considered.
A na\"ive algorithm for computing the Kauffman bracket is therefore to replace all the crossings by their uncrossed versions and add together the contributions from each of the resulting diagrams (which will consist of a number of unlinked, unknotted loops).

The basic idea of the simplification technique is to ``gather terms'' at various stages in this algorithm.
To do that, we need to identify regions of the original diagram where it is likely that the results from applying the crossing\hyp{}replacements will contain repetitions.
By identifying these repetitions, we reduce the number of diagrams that need to be further considered.

There is a simple rule to identify potential regions where this simplification may work.
We consider a region of the diagram and we assume that we can isolate this region by drawing a simple closed curve in the plane which is nowhere tangential to the link.
This will contain a certain number of crossings, say \(k\).
Various strands of the diagram will enter and exit this region (but not simply touch it), say \(n\) distinct strands (thus \(2 n\) distinct entry\hyp{}exit points on the boundary of the region).
At the end of the replacement process, there will be the same number of entry\hyp{}exit points on the boundary and each will still be connected to another such point.
However, the connecting lines can no longer cross.
This puts an upper bound on the number of configurations.
If this upper bound is less than \(2^k\), it is worth ``gathering terms'' before proceeding to another region.

The upper bound satisfies a simple recursion formula:
\[
  a_n = \sum_{k=0}^{n-1} a_k a_{n - 1 - k}
\]
and the first few terms of the sequence are:
\begin{centre}
\begin{tabular}{c|ccccccccc}
\(n\) & \(0\) & \(1\) & \(2\) & \(3\) & \(4\) & \(5\) & \(6\) & \(7\) & \(8\) \\
\hline
\(a_n\) & \(1\) & \(1\) & \(2\) & \(5\) & \(14\) & \(42\) & \(132\) & \(429\) & \(1430\) \\
\(\lceil \log_2(a_n)\rceil\) & \(0\) & \(0\) & \(1\) & \(3\) & \(4\) & \(6\) & \(8\) & \(9\) & \(10\)
\end{tabular}
\end{centre}
The third line in the table is the number of links that need to be in the region in order for the simplification to be worth doing.
We want to find \(k\) such that \(2^k > a_n\), so we want \(k > \log_2 a_n\).

The sequence \((a_n)\) is the sequence of Catalan numbers\footnote{A000108 in The On-Line Encyclopedia of Integer Sequences} which has general formula:
\[
  \frac{1}{n+1}\binom{2n}{n}.
\]

\section{Simplifying Hopf Links}

We can draw the Hopf link prettily as in Figure~\ref{fig:hopf}.
We can also represent the linking part more schematically as in Figure~\ref{fig:schhopflink}.

\tikzsetnextfilename{hopf}

\begin{figure}
\centering
\begin{tikzpicture}[scale=.5]
\draw[knot,double=ring1] (0,0) arc(180:0:3);
\draw[knot,double=ring2] (6,0) circle (3);
\draw[knot,double=ring1] (0,0) arc(-180:0:3);
\end{tikzpicture}
\caption{The Hopf link}
\autolabel
\end{figure}
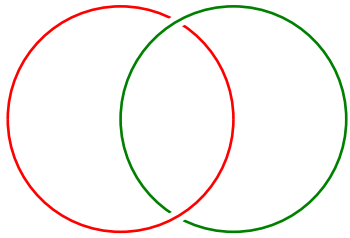

\tikzsetnextfilename{schhopflink}

\begin{figure}
\centering
\begin{tikzpicture}[every path/.style={thick knot,double=Red},every node/.style={knot crossing}]
\node (h1) at (0,0) {};
\node (h2) at (1,-.5) {};
\draw (h1.center) -- ++(0,.5) -- ++(1.5,0);
\draw (h1.center) |- (h2) -- ++(.5,0);
\draw (h2.center) -- ++(0,-.5) -- ++(-1.5,0);
\draw (h2.center) |- (h1) -- ++(-.5,0);
\end{tikzpicture}
\caption{Detail of the linking in the Hopf link}
\autolabel
\end{figure}
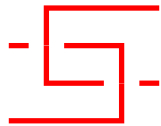

The two links that we are interested in which involve the Hopf linking are the Hopf rings and the Hopf chains.
Figure~\ref{fig:hopfchain} shows a Hopf chain of 8 components.
Figure~\ref{fig:hopfring} shows a Hopf ring also of 8 components.

\tikzsetnextfilename{hopfchain}

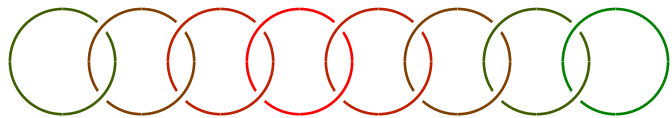
\begin{figure}
\centering
\begin{tikzpicture}
\pgfmathsetmacro{\brradius}{.7}
\pgfmathsetmacro{\brshift}{1.5*\brradius }
\foreach \brk in {1,...,8} {
\pgfmathsetmacro{\brtint}{(\brk > 4 ? 8 - \brk : \brk) * 25}
  \begin{scope}[xshift=\brk * \brshift cm]
  \draw[knot,double=Red!\brtint!Green] (0,0) arc (180:90:\brradius);
  \draw[knot,double=Red!\brtint!Green] (\brradius,-\brradius) arc (-90:0:\brradius);
  \begin{pgfonlayer}{back}
  \draw[knot,double=Red!\brtint!Green] (0,0) arc (-180:-90:\brradius);
  \draw[knot,double=Red!\brtint!Green] (\brradius,\brradius) arc (90:0:\brradius);
  \end{pgfonlayer}
  \end{scope}
}
\end{tikzpicture}
\caption{A Hopf chain of 8 components}
\autolabel
\end{figure}

\tikzsetnextfilename{hopfring}

\begin{figure}
\centering
\begin{tikzpicture}
\pgfmathsetmacro{\brradius}{1}
\pgfmathsetmacro{\brshift}{2*\brradius }
\foreach \k in {1,...,8} {
\pgfmathsetmacro{\brtint}{(\k > 4 ? 8 - \k : \k) * 25}
\begin{scope}[rotate=\k * 45,yshift=\brshift cm]
\draw[knot,double=Red!\brtint!Blue] (0,\brradius) arc (90:0:\brradius);
\draw[knot,double=Red!\brtint!Blue] (0,-\brradius) arc (-90:-180:\brradius);
\begin{pgfonlayer}{front}
\draw[knot,double=Red!\brtint!Blue] (0,\brradius) arc (90:180:\brradius);
\draw[knot,double=Red!\brtint!Blue] (0,-\brradius) arc (-90:0:\brradius);
\end{pgfonlayer}
\end{scope}
}
\end{tikzpicture}
\caption{Hopf ring with 8 components}
\autolabel
\end{figure}
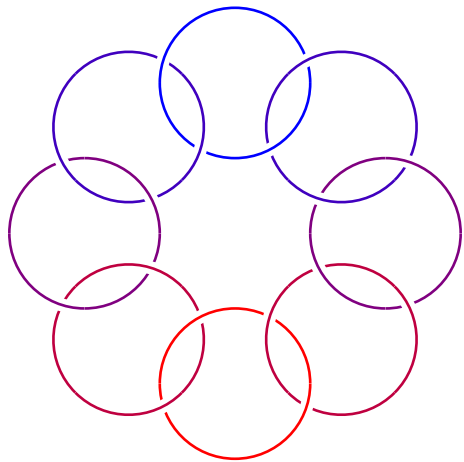

We are particularly interested in the ring, but computing the Kauffman bracket for the ring will involve computing it for the chain as well.
To compute this, we first resolve the Hopf link as in Figure~\ref{fig:hopfresolve}.

\tikzsetnextfilename{hopfresolve}

\begin{figure}
\centering
\begin{tikzpicture}[every path/.style={thick knot,double=Red},every node/.style={knot crossing}]
\node (h1) at (0,0) {};
\node (h2) at (1,-.5) {};
\draw (h1.center) -- ++(0,.5) -- ++(1.5,0);
\draw (h1.center) |- (h2) -- ++(.5,0);
\draw (h2.center) -- ++(0,-.5) -- ++(-1.5,0);
\draw (h2.center) |- (h1) -- ++(-.5,0);
\node[draw=none,text=black] at (2,0) {\(=\)};
\begin{scope}[xshift=4cm]
\node[draw=none,text=black] at (-1,0) {\(A\)};
\node (h1) at (0,0) {};
\node (h2) at (1,-.5) {};
\draw (h1) -- ++(0,.5) -- ++(1.5,0);
\draw (h1) |- (h2) -- ++(.5,0);
\draw (h1.south) -- (h1.west);
\draw (h1.north) -- (h1.east);
\draw (h2.center) -- ++(0,-.5) -- ++(-1.5,0);
\draw (h2.center) |- (h1) -- ++(-.5,0);
\node[draw=none,text=black] at (2,0) {\(+\)};
\begin{scope}[xshift=4cm]
\node[draw=none,text=black] at (-1,0) {\(A^{-1}\)};
\node (h1) at (0,0) {};
\node (h2) at (1,-.5) {};
\draw (h1) -- ++(0,.5) -- ++(1.5,0);
\draw (h1) |- (h2) -- ++(.5,0);
\draw (h1.south) -- (h1.east);
\draw (h1.north) -- (h1.west);
\draw (h2.center) -- ++(0,-.5) -- ++(-1.5,0);
\draw (h2.center) |- (h1) -- ++(-.5,0);
\end{scope}
\end{scope}
\begin{scope}[yshift=-2cm]
\node[draw=none,text=black] at (2,0) {\(=\)};
\begin{scope}[xshift=4cm]
\node[draw=none,text=black] at (-1,0) {\(A^2\)};
\node (h1) at (0,0) {};
\node (h2) at (1,-.5) {};
\draw (h1) -- ++(0,.5) -- ++(1.5,0);
\draw (h1) |- (h2) -- ++(.5,0);
\draw (h1.south) -- (h1.west);
\draw (h1.north) -- (h1.east);
\draw (h2) -- ++(0,-.5) -- ++(-1.5,0);
\draw (h2) |- (h1) -- ++(-.5,0);
\draw (h2.south) -- (h2.west);
\draw (h2.north) -- (h2.east);
\node[draw=none,text=black] at (2,0) {\(+\)};
\begin{scope}[xshift=4cm]
%\node[draw=none,text=black] at (-1,0) {\(A^{-1}\)};
\node (h1) at (0,0) {};
\node (h2) at (1,-.5) {};
\draw (h1) -- ++(0,.5) -- ++(1.5,0);
\draw (h1) |- (h2) -- ++(.5,0);
\draw (h1.south) -- (h1.west);
\draw (h1.north) -- (h1.east);
\draw (h2) -- ++(0,-.5) -- ++(-1.5,0);
\draw (h2) |- (h1) -- ++(-.5,0);
\draw (h2.south) -- (h2.east);
\draw (h2.north) -- (h2.west);
\end{scope}
\end{scope}
\begin{scope}[yshift=-2cm]
\node[draw=none,text=black] at (2.5,0) {\(+\)};
\begin{scope}[xshift=4cm]
%\node[draw=none,text=black] at (-1,0) {\(A^2\)};
\node (h1) at (0,0) {};
\node (h2) at (1,-.5) {};
\draw (h1) -- ++(0,.5) -- ++(1.5,0);
\draw (h1) |- (h2) -- ++(.5,0);
\draw (h1.south) -- (h1.east);
\draw (h1.north) -- (h1.west);
\draw (h2) -- ++(0,-.5) -- ++(-1.5,0);
\draw (h2) |- (h1) -- ++(-.5,0);
\draw (h2.south) -- (h2.west);
\draw (h2.north) -- (h2.east);
\node[draw=none,text=black] at (2,0) {\(+\)};
\begin{scope}[xshift=4cm]
\node[draw=none,text=black] at (-1,0) {\(A^{-2}\)};
\node (h1) at (0,0) {};
\node (h2) at (1,-.5) {};
\draw (h1) -- ++(0,.5) -- ++(1.5,0);
\draw (h1) |- (h2) -- ++(.5,0);
\draw (h1.south) -- (h1.east);
\draw (h1.north) -- (h1.west);
\draw (h2) -- ++(0,-.5) -- ++(-1.5,0);
\draw (h2) |- (h1) -- ++(-.5,0);
\draw (h2.south) -- (h2.east);
\draw (h2.north) -- (h2.west);
\end{scope}
\end{scope}
\begin{scope}[yshift=-2cm]
\node[draw=none,text=black] at (2,0) {\(=\)};
\begin{scope}[xshift=4cm]
\node[draw=none,text=black] at (-1,0) {\(A^2\)};
\node (h1) at (0,0) {};
\node (h2) at (1,-.5) {};
\draw (h1) -- ++(0,.5) -- ++(1.5,0);
\draw (h1.north) -- (h2.east) (h2) -- ++(.5,0);
\draw (h2) -- ++(0,-.5) -- ++(-1.5,0);
\draw (h2.south) -- (h1.west) (h1) -- ++(-.5,0);
\node[draw=none,text=black] at (2,0) {\(+\)};
\begin{scope}[xshift=4.5cm]
\node[draw=none,text=black,anchor=east] at (-.5,0) {\((1  - A^{-4})\)};
\node (h1) at (0,0) {};
\node (h2) at (1,-.5) {};
\draw (h1) -- ++(0,.5) -- ++(1.5,0);
\draw (h1.north) -- (h1.west) (h1) -- ++(-.5,0);
\draw (h2) -- ++(0,-.5) -- ++(-1.5,0);
\draw (h2.south) -- (h2.east) (h2) -- ++(.5,0);
\end{scope}
\end{scope}
\end{scope}
\end{scope}
\end{scope}
\end{tikzpicture}
\caption{Resolving the Hopf linking}
\autolabel
\end{figure}
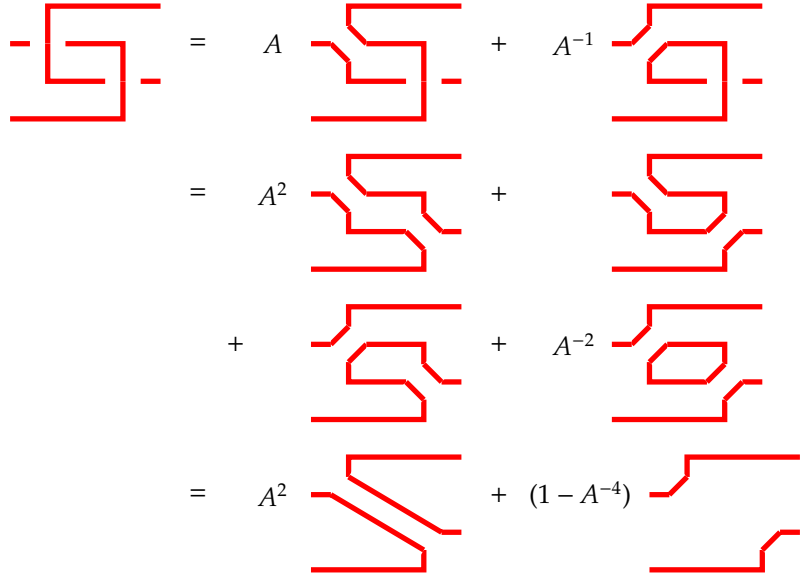

From Figure~\ref{fig:hopfresolve} we see that a Hopf linkage decomposes into two diagrams: one where the two circles involved in the linkage are disjoint and one where the two circles are replaced by just one circle.
Applying this to the Hopf chain, we get for the first a shorter Hopf chain together with a disjoint circle and for the second we get just the shorter Hopf chain.

Let us write \(\kb{hopfchain}{n}\) for the Kauffman bracket of a Hopf chain of \(n\) circles.
Then \(\kb{hopfchain}{1}\) is the Kauffman bracket of a single circle so \(\kb{hopfchain}{1} = 1\).
Resolving the first Hopf linkage, we get the recurrence relation
\[
  \kb{hopfchain}{n} = A^2(-A^2 - A^{-2})   \kb{hopfchain}{n-1} + (1 - A^{-4}) \kb{hopfchain}{n-1} = (-A^4 - A^{-4})  \kb{hopfchain}{n-1}.
\]
Hence
\[
 \kb{hopfchain}{n} = (-A^4 - A^{-4})^{n-1}.
\]

When we resolve a Hopf linkage in the Hopf ring, we obtain on the one hand a Hopf chain with the same number of components and on the other a Hopf ring with one less component.
Let us write \(\kb{hopfring}{n}\) for the Kauffman bracket of a Hopf ring of \(n\) circles.
Then we have the recurrence relation
\[
  \kb{hopfring}{n} = A^2 \kb{hopfchain}{n} + (1 - A^{-4}) \kb{hopfring}{n-1}.
\]
Substituting in for \(c_n(A)\) we obtain
\begin{equation}
\label{eq:rn}
  \kb{hopfring}{n} = A^2(-A^4 - A^{-4})^{n-1} + (1 - A^{-4}) \kb{hopfring}{n-1}.
\end{equation}

The initial term in this sequence is \(\kb{hopfring}{1}\) but this is not \(1\).  The polynomial \(\kb{hopfring}{1}\) is the Kauffman bracket of a single circle with a self Hopf linkage, as in Figure~\ref{fig:onehopfring}.
This is not trivial as a \emph{framed} link and we can resolve it one stage further using the resolution of the Hopf linkage.
The \(A^2\) term is of a single circle whilst the \((1 - A^{-4})\) term is of two circles.
Thus
\[
  \kb{hopfring}{1} = A^2 + (1 - A^{-4})(-A^2 - A^{-2}) = A^2 - A^2 - A^{-2} + A^{-2} + A^{-6} = A^{-6}.
\]
Alternatively, we could apply the recursion one stage further to get:
\[
\kb{hopfring}{1} = A^2 + (1 - A^{-4}) \kb{hopfring}{0}.
\]
By our numbering scheme, \(\kb{hopfring}{0}\) ought to be the Kauffman bracket of a Hopf ring with no components.
However, we could note that the numbering also refers to the number of Hopf linkages in the diagram.
Removing the linkage from Figure~\ref{fig:onehopfring} leaves a double circle.
Thus \(\kb{hopfring}{0} = -A^2 - A^{-2}\) and we compute:
\[
  \kb{hopfring}{1} = A^2 + (1 - A^{-4})(-A^2 - A^{-2}) = A^2 - A^2 - A^{-2} + A^{-6} = A^{-6}.
\]

\tikzsetnextfilename{onehopfring}

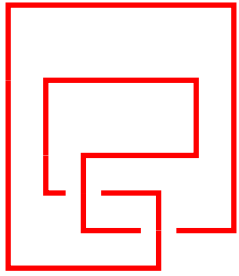
\begin{figure}
\centering
\begin{tikzpicture}[every path/.style={thick knot,double=Red},every node/.style={knot crossing}]
\node (h1) at (0,0) {};
\node (h2) at (1,-.5) {};
\draw (h1.center) -- ++(0,.5) -- ++(1.5,0) -- ++(0,1) -- ++(-2,0) -- ++(0,-1);
\draw (h1.center) |- (h2) -- ++(1,0) -- ++(0,3) -- ++(-3,0) -- ++(0,-1);
\draw (h2.center) -- ++(0,-.5) -- ++(-2,0) -- ++(0,2.5);
\draw (h2.center) |- (h1) -- ++(-.5,0) -- ++(0,.5);
\end{tikzpicture}
\caption{A Single Hopf Linkage}
\autolabel
\end{figure}

\begin{proposition}
\[
  \kb{hopfring}{n} = -A^2 \frac{(A^4 - A^{-8})(1 - A^{-4})^{n-1} + (-A^4 - A^{-4})^n}{1 + A^4}.
\]
\end{proposition}

\begin{proof}
Setting \(n = 1\), we simplify as follows.
\begin{align*}
\kb{hopfring}{1} &= -A^2 \frac{(A^4 - A^{-8}) + (-A^4 - A^{-4})}{1 + A^4} \\
&= -A^2 \frac{-A^{-4} - A^{-8}}{1 + A^4} \\
&= -A^{-6} \frac{-A^4 - 1}{1 + A^4} \\
&= A^{-6}
\end{align*}

For arbitrary \(n\), we substitute \(\kb{hopfring}{n-1}\) in to the right\hyp{}hand side of \eqref{eq:rn} and simplify.
We start with the \((1 - A^{-4})\kb{hopfring}{n-1}\) term.
\begin{align*}
(1 - A^{-4}) \kb{hopfring}{n-1} &= (1 - A^{-4}) \left( -A^2 \frac{(A^4 - A^{-8})(1 - A^{-4})^{n-2} + (-A^4 - A^{-4})^{n-1}}{1 + A^4}\right) \\
&= -A^2 \frac{(A^4 - A^{-8})(1 - A^{-4})^{n-1} + (1 - A^{-4})(-A^4 - A^{-4})^{n-1}}{1 + A^4}.
\end{align*}
The other term contributes \(A^2(-A^4 - A^{-4})^{n-1}\) which we rewrite as follows.
\[
A^2(-A^4 - A^{-4})^{n-1} = -A^2 \frac{-(1 + A^4)(-A^4 - A^{-4})^{n-1}}{1 + A^4}.
\]
If we add these two and concentrate on the pieces involving \((-A^4 - A^{-4})\) we simplify as follows.
\begin{align*}
(1 - A^{-4})(-A^4 - A^{-4})^{n-1} -(1 + A^4)(-A^4 - A^{-4})^{n-1} &= (-A^{-4} - A^4)(-A^4 - A^{-4})^{n-1} \\
&= (-A^4 - A^{-4})^n.
\end{align*}
And thus:
\begin{align*}
  A^2(-A^4 - A^{-4})^{n-1} &+ (1 - A^{-4}) r_{n-1}(A) \\
  &=  -A^2 \frac{(A^4 - A^{-8})(1 - A^{-4})^{n-1} + (-A^4 - A^{-4})^n}{1 + A^4} \\
&= \kb{hopfring}{n}
\end{align*}
as required.
\end{proof}

\textbf{Note to self}: In the ring, it is possible to have different orientations of the Hopf linkages.
In the chain, these can be flipped to all be the same, but not in the ring.

\section{Level Two Hopf Links}

Using Hopf rings, we can build level two links.
In these, the Hopf rings are the building blocks for higher linking.
Starting with the Hopf link, we replace each of the circles by a Hopf ring.
In order to draw this in a reasonable fashion, we first deform the circles in the Hopf ring, flattening them somewhat.
In Figure~\ref{fig:flathopfring}, we show the undeformed and deformed Hopf rings side by side.
This makes the level two Hopf link easier to see, as shown in Figure~\ref{fig:l2hopflink}.

\tikzsetnextfilename{flathopfring}

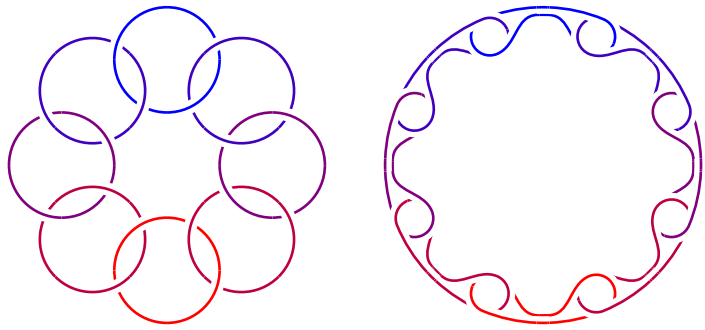
\begin{figure}
\centering
\begin{tikzpicture}
\pgfmathsetmacro{\brradius}{.7}
\pgfmathsetmacro{\brshift}{2*\brradius }
\foreach \k in {1,...,8} {
\pgfmathsetmacro{\brtint}{(\k > 4 ? 8 - \k : \k) * 25}
\begin{scope}[rotate=\k * 45,yshift=\brshift cm]
\draw[knot,double=Red!\brtint!Blue] (0,\brradius) arc (90:0:\brradius);
\draw[knot,double=Red!\brtint!Blue] (0,-\brradius) arc (-90:-180:\brradius);
\begin{pgfonlayer}{front}
\draw[knot,double=Red!\brtint!Blue] (0,\brradius) arc (90:180:\brradius);
\draw[knot,double=Red!\brtint!Blue] (0,-\brradius) arc (-90:0:\brradius);
\end{pgfonlayer}
\end{scope}
}
\begin{scope}[xshift=5cm,rotate=22.5]
\setbrstep{.1}
\foreach \k in {1,...,8} {
\pgfmathsetmacro{\brtint}{(\k > 4 ? 8 - \k : \k) * 25}
\xglobal\colorlet{tmpcol\k}{ring\k}
\xglobal\colorlet{ring\k}{Blue!\brtint!Red}
}
\hopfring{2}{8}
\foreach \k in {1,...,8} {
\xglobal\colorlet{ring\k}{tmpcol\k}
}
\end{scope}
\end{tikzpicture}
\caption{Flattened Hopf ring}
\autolabel
\end{figure}

\tikzsetnextfilename{l2hopflink}

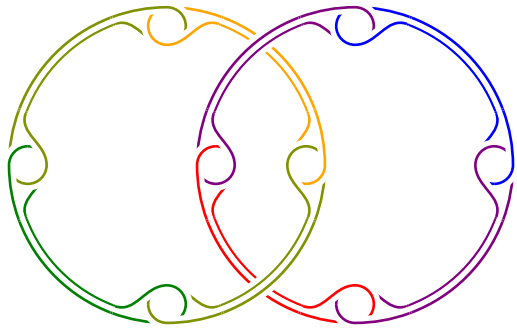
\begin{figure}
\centering
\begin{tikzpicture}
\foreach \k in {1,...,4} {
\pgfmathsetmacro{\brtint}{(\k > 2 ? 4 - \k : \k) * 50}
\xglobal\colorlet{tmpcol\k}{ring\k}
\xglobal\colorlet{ring\k}{Orange!\brtint!Green}
}
\setbrstep{.1}
\hopfring{2}{4}
\begin{scope}[xshift=2.5cm]
\foreach \k in {1,...,4} {
\pgfmathsetmacro{\brtint}{(\k > 2 ? 4 - \k : \k) * 50}
\xglobal\colorlet{tmpcol\k}{ring\k}
\xglobal\colorlet{ring\k}{Blue!\brtint!Red}
}
\hopfring{2}{4}
\foreach \k in {1,...,4} {
\xglobal\colorlet{ring\k}{tmpcol\k}
}
\end{scope}
\end{tikzpicture}
\caption{Level Two Hopf Link}
\autolabel
\end{figure}

We can now consider chains of these rings, or rings of rings.
As each circle can have an arbitrary number of components, we get a myriad of possibilities.
Let us write \(\kb{ltwohopfchain}{n_1,\dotsc,n_k}\) for the Kauffman bracket of a level two chain where the \(j\)th component is itself a Hopf ring with \(n_j\) components.
Let us write \(\kb{ltwohopfring}{n_1,\dotsc,n_k}\) for the corresponding ring.
Note that \(\kb{ltwohopfring}{n_1,\dotsc,n_k}\) is invariant under cyclic permutations of its indices.

Let us consider a level two Hopf chain, say of length \(k\) with the \(j\)th component having \(n_j\) components.
Choose a component, say the \(j\)th.
This is a Hopf ring so its circles are connected by Hopf linkages.
If we resolve one of those Hopf linkages we get a Hopf chain (of the same number of circles) and a Hopf ring (of one fewer circles).
The Hopf chain can be disconnected from the level two Hopf chain, leaving two segments of the level two chain, whilst the Hopf ring remains attached.
The Kauffman bracket of two disjoint links is \((-A^2 - A^{-2})\) times the product of the brackets of the two links in isolation.
Thus we obtain the recurrence relation:
\begin{dmath}
\label{eq:l2hopfchain}
  \kb{ltwohopfchain}{n_1,\dotsc,n_k} = A^2(-A^2 - A^{-2})^{1 - \delta_{1 j} + 1 - \delta_{k j}} \kb{ltwohopfchain}{n_1,\dotsc,n_{j-1}} \kb{hopfchain}{n_j} \kb{ltwohopfchain}{n_{j+1},\dotsc,n_k} + (1 - A^{-4})   \kb{ltwohopfchain}{n_1,\dotsc,n_j - 1,\dotsc,n_k}.
\end{dmath}

At one extreme, we obtain products of the form
\[
\kb{ltwohopfchain}{n_1} \kb{hopfchain}{n_2} \dotsb \kb{ltwohopfchain}{n_{k-1}} \kb{hopfchain}{n_k},
\]
depending on where we have broken our level two Hopf chains.
The factors \(\kb{ltwohopfchain}{n}\) in this are level two Hopf chains with a single component (at the second level).
This is just an ordinary Hopf ring.
That is, \(\kb{ltwohopfchain}{n} = \kb{hopfring}{n}\).

Taking the other path, we soon reach things like \(\kb{ltwohopfchain}{1,n_2,\dotsc,n_k}\).
A simple example of the corresponding link is in Figure~\ref{fig:l2hopfone}.
Resolving the final Hopf linkage in the first component, on the one hand we have a free circle and a shorter chain, whilst on the other, we have two circles still linked to the rest of the chain as in Figure~\ref{fig:l2hopfzero}.
If we refer to the Kauffman bracket of this as \(\kb{ltwohopfchain}{0,n_2,\dotsc,n_k}\) then \eqref{eq:l2hopfchain} still holds.
Thus we can proceed with our iteration until all the level one components have no links and are just double circles.
The linking between them is what could be called a \emph{double Hopf linkage}.
We draw it schematically in Figure~\ref{fig:schdblehopf}.

\tikzsetnextfilename{l2hopfone}

\begin{figure}
\centering
\begin{tikzpicture}
\setbrstep{.1}
\hopfring{2}{1}
\begin{scope}[xshift=2.5cm]
\hopfring{2}{4}
\end{scope}
\end{tikzpicture}
\caption{Level two Hopf link with one circle in the first component}
\autolabel
\end{figure}
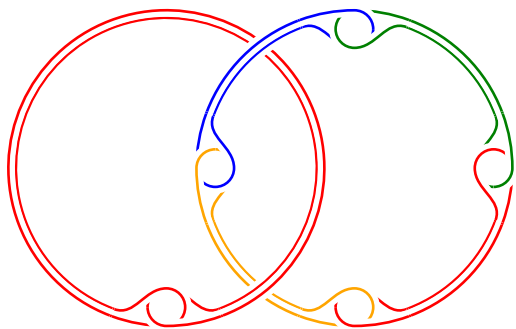

\tikzsetnextfilename{l2hopfzero}

\begin{figure}
\centering
\begin{tikzpicture}
\setbrstep{.1}
\draw[knot,double=ring1] (0,0) circle (2);
\draw[knot,double=ring1] (0,0) circle (2 + .5*\brstep);
\begin{scope}[xshift=2.5cm]
\hopfring{2}{4}
\end{scope}
\end{tikzpicture}
\caption{Level two Hopf link with no junctions in the first component}
\autolabel
\end{figure}
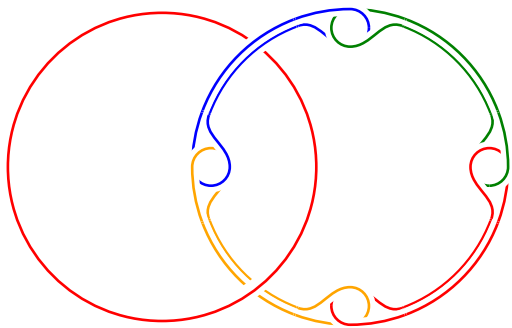

\tikzsetnextfilename{schdblehopf}

\begin{figure}
\centering
\begin{tikzpicture}[every path/.style={thick knot,double=Red}]
\pgfmathsetmacro{\brover}{.5}
\pgfmathsetmacro{\bradj}{.1}
\foreach \bri/\brdir/\brnud/\brv in {
  1/-1/2/-1,
  2/-1/1/-1,
  3/1/2/-1,
  4/1/1/-1,
  5/-1/1/1,
  6/-1/2/1,
  7/1/1/1,
  8/1/2/1
} {
\coordinate (\bri) at (\brover * \brdir * \brnud,-\bri/2 + \brv * \bradj);
\draw (\bri) -- ++(0,-\brv * \bradj) -- +(-\brdir * 2 - \brover * \brdir * \brnud,0);
}
\draw (1) -- (6);
\draw (2) -- (5);
\draw (3) -- (8);
\draw (4) -- (7);
\end{tikzpicture}
\caption{Detail of the linking in the double Hopf link}
\autolabel
\end{figure}
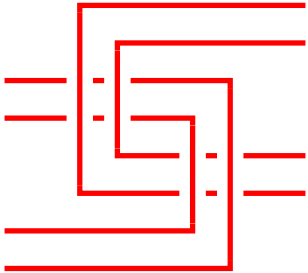

Before considering this diagram, let us consider a level two Hopf ring as in Figure~\ref{fig:l2hopfring}.
Again, for drawing it is convenient to flatten the lower level structure and render the link as in Figure~\ref{fig:l2hopfringflat}.

\tikzsetnextfilename{l2hopfring}

\begin{figure}
\centering
\begin{tikzpicture}
\setbrstep{.1}
\pgfmathsetmacro{\last}{10}
\pgfmathsetmacro{\angle}{360/10}
\foreach \l in {1,...,\last} {
\begin{scope}[rotate=\angle * \l,yshift=-4.5cm]
\begin{scope}[rotate=90]
\hopfring{2}{4}
\end{scope}
\end{scope}
}
\end{tikzpicture}
\caption{Level two Hopf ring}
\autolabel
\end{figure}
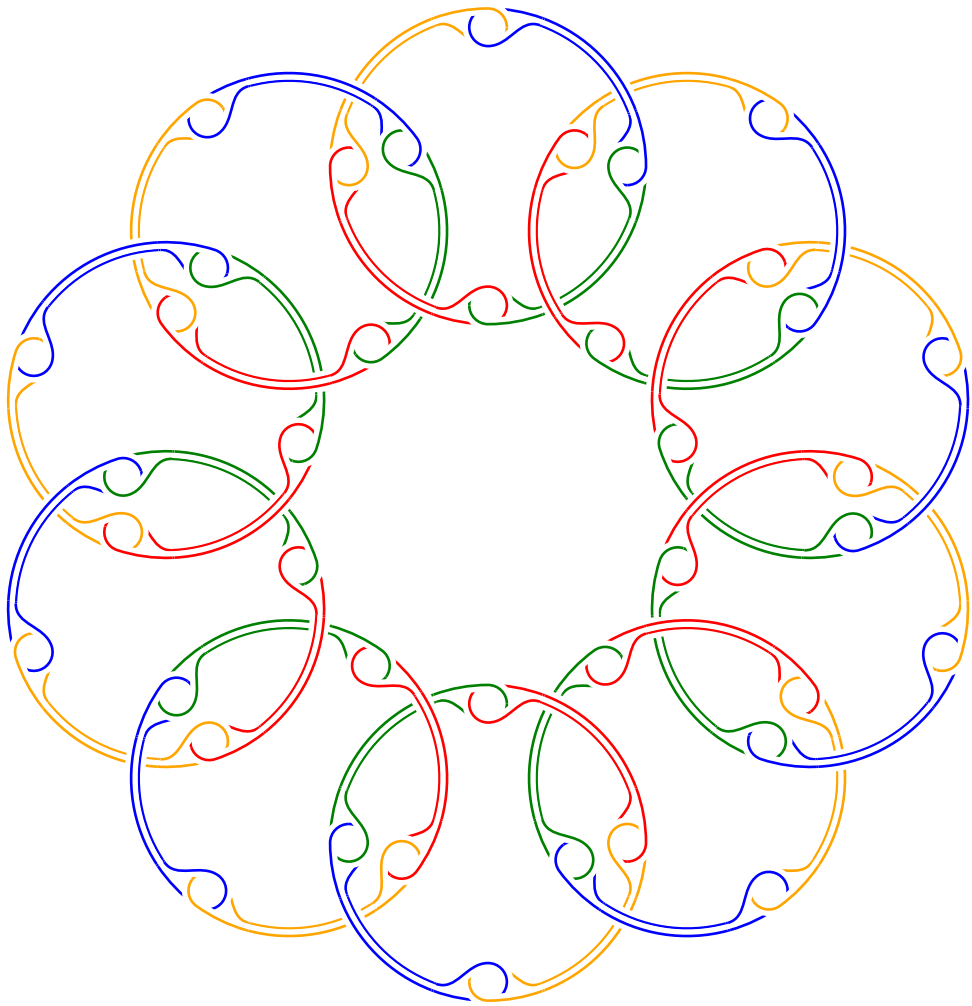

\tikzsetnextfilename{l2hopfringflat}

\begin{figure}
\centering
\begin{tikzpicture}
\setbrstep{.1}
\hopftwo{2}{3}{3}
\end{tikzpicture}
\caption{Flattened level two Hopf ring}
\autolabel
\end{figure}
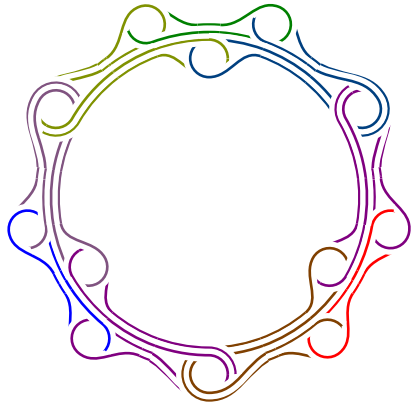

Recall that we write \(\kb{ltwohopfring}{n_1,\dotsc, n_k}\) for the Kauffman bracket of the level two Hopf ring formed from \(k\) Hopf rings of \(n_1, \dotsc, n_k\) components.
We choose a component, say \(j\), and resolve one of the Hopf linkages in this component.
As for the Hopf chain, this leaves us with two diagrams; one in which the Hopf ring under consideration becomes a Hopf chain and one in which it becomes a Hopf ring with one fewer components.
In the first case, the level two Hopf ring becomes a level two Hopf chain with a disjoint level one Hopf chain.
In the second, we still have a level two Hopf ring.
Thus our first recurrence relation is:
\begin{dmath*}
\kb{ltwohopfring}{n_1,\dotsc,n_k} = A^2(-A^2 - A^{-2})\kb{hopfchain}{n_j} \kb{ltwohopfchain}{n_{j+1},\dotsc,n_k,n_1,\dotsc,n_{j-1}} + (1 - A^{-4}) \kb{ltwohopfring}{n_1,\dotsc,n_j-1,\dotsc,n_k}.
\end{dmath*}
Arguing as with the level two Hopf chain, we see that we need to consider \(\kb{ltwohopfring}{0,\dotsc,0}\) which is the Kauffman bracket of a link formed by linking pairs of circles in a ring.

Thus resolving the linkages in the level two Hopf chain and ring lead us to consider the links shown in Figure~\ref{fig:dbleHopfchain} and Figure~\ref{fig:dbleHopfring}.

\tikzsetnextfilename{dbleHopfchain}

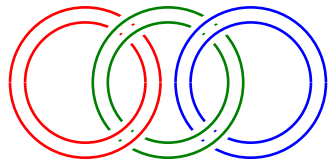
\begin{figure}
\centering
\begin{tikzpicture}[every path/.style={thick knot,double=Red},every node/.style={text=black}]
\draw[knot,double=ring1] (-.1,0) arc(180:0:1);
\draw[knot,double=ring1] (.1,0) arc(180:0:.8);
\draw[knot,double=ring3] (2.1,0) arc(-180:0:1);
\draw[knot,double=ring3] (2.3,0) arc(-180:0:.8);
\draw[knot,double=ring2] (2,0) circle (1);
\draw[knot,double=ring2] (2,0) circle (.8);
\draw[knot,double=ring1] (-.1,0) arc(-180:0:1);
\draw[knot,double=ring1] (.1,0) arc(-180:0:.8);
\draw[knot,double=ring3] (2.1,0) arc(180:0:1);
\draw[knot,double=ring3] (2.3,0) arc(180:0:.8);
\end{tikzpicture}
\caption{The double Hopf chain}
\autolabel
\end{figure}

\tikzsetnextfilename{dbleHopfring}

\begin{figure}
\centering
\begin{tikzpicture}
\setbrstep{.1}

\foreach \brk in {1,...,3} {
  \pgfmathparse{int(Mod(\brk,3) + 1)}
  \edef\brm{\pgfmathresult}
  \begin{scope}[rotate = \brk * 120]
   \colorlet{outer}{ring\brm}
   \colorlet{inner}{ring\brk}
  \hopfjunction{2}{2}{-1}{60}{60}
  \end{scope}
}
\end{tikzpicture}
\caption{The double Hopf ring}
\autolabel
\end{figure}
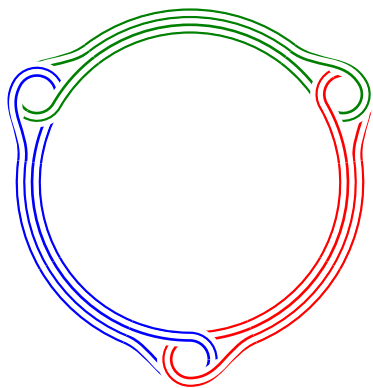

The double Hopf chain is the simpler to resolve.
We separate the circles at one end and deal with them one at a time.
This leads us to consider the ``one and a half'' Hopf linkage.
This resolves as in Figure~\ref{fig:onehalfHopfres}.

\tikzsetnextfilename{onehalfHopfres}

\begin{figure}
\centering
\begin{tikzpicture}[every path/.style={thick knot,double=Red},every node/.style={text=black}]
\pgfmathsetmacro{\brover}{.25}
\pgfmathsetmacro{\bradj}{.05}
\pgfmathsetmacro{\brsep}{.25}
\pgfmathsetmacro{\brlen}{1}
\foreach \bri/\brdir/\brnud/\brv in {
  1/-1/2/-1,
  2/-1/1/-1,
  3/1/1/-1,
  4/-1/1/1,
  5/-1/2/1,
  6/1/1/1
} {
\coordinate (\bri) at (\brover * \brdir * \brnud,-\bri * \brsep + \brv * \bradj);
\draw (\bri) -- ++(0,-\brv * \bradj) -- +(-\brdir * \brlen - \brover * \brdir * \brnud,0) coordinate (e\bri);
}
\draw
  (1) -- (5)
  (2) -- (4)
  (3) -- (6)
  (e3) -- ++(-.2,0) |- (e6)
;

\begin{scope}[xshift = 3.5cm]
\node[anchor=mid east] at (-1.3,-1) {\(= A^2\)};
\foreach \bri/\brdir/\brnud/\brv in {
  1/-1/2/-1,
  2/-1/1/-1,
  3/-1/1/1,
  4/1/1/-1,
  5/-1/2/1,
  6/1/1/1
} {
\coordinate (\bri) at (\brover * \brdir * \brnud,-\bri * \brsep + \brv * \bradj);
\draw (\bri) -- ++(0,-\brv * \bradj) -- +(-\brdir * \brlen - \brover * \brdir * \brnud,0) coordinate (e\bri);
}
\draw
  (1) -- (5)
  (2) -- (3)
  (4) -- (6)
  (e4) -- ++(-.2,0) |- (e6)
;
\begin{scope}[xshift=4.5cm]
\node[anchor=mid east] at (-1.3,-1) {\(+(1 - A^{-4})\)};
\foreach \bri/\brdir/\brnud/\brv in {
  1/-1/2/-1,
  2/-1/1/-1,
  3/1/-1/1,
  4/-1/-1/-1,
  5/-1/2/1,
  6/1/1/1
} {
\coordinate (\bri) at (\brover * \brdir * \brnud,-\bri * \brsep + \brv * \bradj);
\draw (\bri) -- ++(0,-\brv * \bradj) -- +(-\brdir * \brlen - \brover * \brdir * \brnud,0) coordinate (e\bri);
}
\draw
  (1) -- (5)
  (2) -- (3)
  (4) -- (6)
  (e3) -- ++(-.2,0) |- (e6)
;
\end{scope}
\begin{scope}[yshift=-3cm]
\node[anchor=mid east] at (-1.3,-1) {\(= A^4\)};
\foreach \bri/\brdir/\brnud/\brv in {
  1/-1/2/-1,
  2/-1/1/-1,
  3/-1/1/1,
  4/-1/2/1,
  5/1/1/-1,
  6/1/1/1
} {
\coordinate (\bri) at (\brover * \brdir * \brnud,-\bri * \brsep + \brv * \bradj);
\draw (\bri) -- ++(0,-\brv * \bradj) -- +(-\brdir * \brlen - \brover * \brdir * \brnud,0) coordinate (e\bri);
}
\draw
  (1) -- (4)
  (2) -- (3)
  (5) -- (6)
  (e5) -- ++(-.2,0) |- (e6)
;
\begin{scope}[xshift=5cm]
\node[anchor=mid east] at (-1.3,-1) {\(+ A^2(1 - A^{-4})\)};
\foreach \bri/\brdir/\brnud/\brv in {
  1/-1/2/-1,
  2/-1/1/-1,
  3/-1/1/1,
  4/1/-2/1,
  5/-1/-1/-1,
  6/1/1/1
} {
\coordinate (\bri) at (\brover * \brdir * \brnud,-\bri * \brsep + \brv * \bradj);
\draw (\bri) -- ++(0,-\brv * \bradj) -- +(-\brdir * \brlen - \brover * \brdir * \brnud,0) coordinate (e\bri);
}
\draw
  (1) -- (4)
  (2) -- (3)
  (5) -- (6)
  (e4) -- ++(-.2,0) |- (e6)
;
\end{scope}
\begin{scope}[yshift=-3cm]
\node[anchor=mid east] at (-1.3,-1) {\(+(1 - A^{-4})A^{-2}\)};
\foreach \bri/\brdir/\brnud/\brv in {
  1/-1/2/-1,
  2/-1/1/-1,
  3/1/-2/1,
  4/-1/1/1,
  5/-1/-1/-1,
  6/1/1/1
} {
\coordinate (\bri) at (\brover * \brdir * \brnud,-\bri * \brsep + \brv * \bradj);
\draw (\bri) -- ++(0,-\brv * \bradj) -- +(-\brdir * \brlen - \brover * \brdir * \brnud,0) coordinate (e\bri);
}
\draw
  (1) -- (3)
  (2) -- (4)
  (5) -- (6)
  (e3) -- ++(-.2,0) |- (e6)
;
\begin{scope}[xshift=5cm]
\node[anchor=mid east] at (-.3,-1) {\(+(1 - A^{-4})(1 - A^4)\)};
\foreach \bri/\brdir/\brnud/\brv in {
  1/-1/2/-1,
  2/-1/2/1,
  4/-1/-1/-1,
  5/-1/-1/1
} {
\coordinate (\bri) at (\brover * \brdir * \brnud,-\bri * \brsep + \brv * \bradj);
\draw (\bri) -- ++(0,-\brv * \bradj) -- +(-\brdir * \brlen - \brover * \brdir * \brnud,0) coordinate (e\bri);
}
\draw
  (1) -- (2)
  (4) -- (5)
;
\end{scope}
\begin{scope}[yshift=-3cm]
\node[anchor=mid east] at (-.7,-.7) {\(= (- A^{6} - A^{-6})\)};
\foreach \bri/\brdir/\brnud/\brv in {
  1/-1/2/-1,
  2/-1/1/-1,
  3/-1/1/1,
  4/-1/2/1
} {
\coordinate (\bri) at (\brover * \brdir * \brnud,-\bri * \brsep + \brv * \bradj);
\draw (\bri) -- ++(0,-\brv * \bradj) -- +(-\brdir * \brlen - \brover * \brdir * \brnud,0) coordinate (e\bri);
}
\draw
  (1) -- (4)
  (2) -- (3)
;
\begin{scope}[xshift=4.5cm]
\node[anchor=mid east] at (-.7,-.7) {\(+ (2 - A^{4} - A^{-4})\)};
\foreach \bri/\brdir/\brnud/\brv in {
  1/-1/1/-1,
  2/-1/1/1,
  3/-1/1/-1,
  4/-1/1/1
} {
\coordinate (\bri) at (\brover * \brdir * \brnud,-\bri * \brsep + \brv * \bradj);
\draw (\bri) -- ++(0,-\brv * \bradj) -- +(-\brdir * \brlen - \brover * \brdir * \brnud,0) coordinate (e\bri);
}
\draw
  (1) -- (2)
  (4) -- (3)
;
\end{scope}
\end{scope}
\end{scope}
\end{scope}
\end{scope}
\end{tikzpicture}
\caption{Resolving the one\hyp{}and\hyp{}a\hyp{}half Hopf linkage}
\autolabel
\end{figure}

In the double Hopf chain, the \((-A^6 - A^{-6})\)\enhyp{}factor still has another linked circle to be removed.
In the \((2 - A^4 - A^{-4})\)\enhyp{}factor, this extra circle now disconnects.
The two capped\hyp{}off ends from Figure~\ref{fig:onehalfHopfres} are actually part of the same component, which also disconnects from the main link.
Writing \(\kb{dblhopfchain}{2n}\) for the Kauffman bracket of the double Hopf chain with \(n\) doubled components (and thus \(2 n\) actual components), we therefore find that we have the relation:
\begin{dmath}
\label{eq:dblhchain}
  \kb{dblhopfchain}{2n} = (A^6 + A^{-6})^2  \kb{dblhopfchain}{2n-2} + (-A^{12} + A^4 + A^{-4} - A^{-12}) \kb{dblhopfchain}{2n-4}
\end{dmath}
with starting point \(\kb{dblhopfchain}{2} = (-A^2 - A^{-2})\) and
\begin{align*}
  \kb{dblhopfchain}{4} &= (A^6 + A^{-6})^2 (-A^2 - A^{-2}) + (-A^6 - A^2 - A^{-2} - A^{-6})(2 - A^4 - A^{-4}) \\
&= - A^{14} - A^{6} - 2 A^{2} - 2 A^{-2} - A^{-6} - A^{-14}.
\end{align*}
(Or we could declare \(\kb{dblhopfchain}{0} = (-A^2 - A^{-2})^{-1}\).)

When working with the double Hopf ring, things are not so simple.
The double Hopf chain simplified nicely because we could start at one end and work towards the other.
With the double Hopf ring, this obviously isn't possible.
We therefore need to resolve the double Hopf junction without assuming anything about how it continues on either side.
This is where the simplification technique of Section~\ref{sec:simple} comes to the fore.
When resolving the double Hopf link, there eight crossings so there will be \(256\) terms.
However, there are at most fourteen distinct diagrams that can result.
Therefore we can collapse the \(256\) down to just \(14\) and consider only those.
Each of the \(14\) diagrams will come with a polynomial factor which counts how it contributes to the Kauffman bracket of the original link.

To work out the contributions, we use a computer program.
The source code of the program is listed in Appendix~\ref{ap:kauffman-code}.
Let us briefly explain the program.
We start with Figure~\ref{fig:schdblehopf}.
In this diagram, we label the strands between the crossings.
This includes the strands that leave the linkage (for simplicity, we label these first).
This produces Figure~\ref{fig:schdblehopflbl}.
Each crossing is therefore associated with four labels.
We list those labels in clockwise order, starting with one of the strands corresponding to the over part of the crossing (it doesn't matter which is chosen).
The program then iterates over all resolutions of the crossings.
For each crossing, it links two of the strands.
It then starts at an entry point, follows the strands, and finds where it leaves.
The list of pairings of entry and exit points determines which of the possible diagrams is produced.
There is a further complication in that some isolated circles may also be produced.
The program checks for these by looking for unused strands.

\tikzsetnextfilename{schdblehopflbl}

\begin{figure}
\centering
\begin{tikzpicture}[every path/.style={thick knot,double=Red},every node/.style={text=black}]
\pgfmathsetmacro{\brover}{.5}
\pgfmathsetmacro{\bradj}{.1}
\foreach \bri/\brdir/\brnud/\brv in {
  1/-1/2/-1,
  2/-1/1/-1,
  3/1/2/-1,
  4/1/1/-1,
  5/-1/1/1,
  6/-1/2/1,
  7/1/1/1,
  8/1/2/1
} {
\coordinate (\bri) at (\brover * \brdir * \brnud,-\bri/2 + \brv * \bradj);
\draw (\bri) -- ++(0,-\brv * \bradj) -- +(-\brdir * 2 - \brover * \brdir * \brnud,0) coordinate (e\bri);
}
\draw (1) -- (6);
\draw (2) -- (5);
\draw (3) -- (8);
\draw (4) -- (7);
\node at (e1) {5};
\node at (e2) {6};
\node at (e3) {1};
\node at (e4) {2};
\node at (e5) {7};
\node at (e6) {8};
\node at (e7) {3};
\node at (e8) {4};
\node at ($(e3)!.45!(3)$) {9};
\node at ($(e4)!.5!(4)$) {10};
\node at ($(1)!.5!(6)$) {11};
\node at ($(2)!.5!(5)$) {12};
\node at (3) {13};
\node at (4) {14};
\node at (5) {15};
\node at (6) {16};
\node at ($(e5)!.5!(5)$) {17};
\node at ($(e6)!.45!(6)$) {18};
\node at ($(3)!.5!(8)$) {19};
\node at ($(4)!.5!(7)$) {20};
\end{tikzpicture}
\caption{Labelled double Hopf link}
\autolabel
\end{figure}
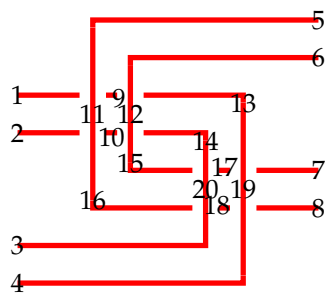

Having found the contributions of the fourteen possible diagrams, we  then put back the information that we have temporarily forgotten: namely that these form part of a double Hopf ring.
Remembering that, we can simplify the resulting fourteen diagrams considerably.
We assume that we start with a Hopf ring of \(2n\) components, thus \(n\) doubled components.
We are considering a segment of this ring where we show three linkages, as in Figure~\ref{fig:dblhopf3link}.
What is important to remember from this diagram is that there are four complete components shown and thus removing these four leaves a link with \(2n - 4\) components.

\tikzsetnextfilename{dblhopf3link}

\begin{figure}
\centering
\begin{tikzpicture}[every path/.style={thick knot,double=Red},every node/.style={text=black},baseline=0cm]
\pgfmathsetmacro{\brover}{.25}
\pgfmathsetmacro{\bradj}{.05}
\pgfmathsetmacro{\brsep}{.25}
\pgfmathsetmacro{\brlen}{1}
\foreach \brj in {0,...,2} {
\begin{scope}[xshift=3 * \brj cm]
\foreach \bri/\brdir/\brnud/\brv in {
  1/-1/2/-1,
  2/-1/1/-1,
  3/1/2/-1,
  4/1/1/-1,
  5/-1/1/1,
  6/-1/2/1,
  7/1/1/1,
  8/1/2/1
} {
\coordinate (\brj\bri) at (\brover * \brdir * \brnud,-\bri * \brsep + \brv * \bradj);
\draw (\brj\bri) -- ++(0,-\brv * \bradj) -- +(-\brdir * \brlen - \brover * \brdir * \brnud,0) coordinate (e\brj\bri);
}
\draw (\brj1) -- (\brj6);
\draw (\brj2) -- (\brj5);
\draw (\brj3) -- (\brj8);
\draw (\brj4) -- (\brj7);
\end{scope}
}
\draw
 (e01) -- (e13)
 (e02) -- (e14)
 (e05) -- (e17)
 (e06) -- (e18)
 (e11) -- (e23)
 (e12) -- (e24)
 (e15) -- (e27)
 (e16) -- (e28)
;
\end{tikzpicture}
\caption{Template for the replacement diagrams in the double Hopf ring}
\autolabel
\end{figure}
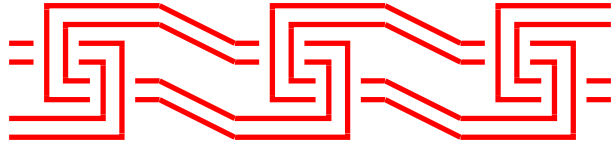

\begin{enumerate}
\pgfmathsetmacro{\brover}{.25}
\pgfmathsetmacro{\bradj}{.05}
\pgfmathsetmacro{\brsep}{.25}
\pgfmathsetmacro{\brlen}{1}

\item
\tikzsetnextfilename{dblhopfdia1}
\autolabel
\begin{tikzpicture}[every path/.style={thick knot,double=Red},every node/.style={text=black},baseline=0cm]
\dblhopftwoside
\draw[double=none,Red,line width=2pt] (1) -- (2) (3) -- (5) (4) -- (6) (7) -- (8) ;
\end{tikzpicture}
\[
A^{-2}-A^{2}
\]

We remove the central component (at a cost of \((-A^2 - A^{-2})\)), leaving a double Hopf chain with \(2n-4\) components.

\item
\tikzsetnextfilename{dblhopfdia2}
\autolabel
\begin{tikzpicture}[every path/.style={thick knot,double=Red},every node/.style={text=black},baseline=0cm]
\dblhopftwoside
\draw[double=none,Red,line width=2pt] (1) -- (5) (2) -- (3) (4) -- (8) (6) -- (7) ;
\end{tikzpicture}
\[
A^{-6}+A^{-2}-A^{2}-A^{6}
\]

This link admits no immediate simplification and so is analysed later.

\item
\tikzsetnextfilename{dblhopfdia3}
\autolabel
\begin{tikzpicture}[every path/.style={thick knot,double=Red},every node/.style={text=black},baseline=0cm]
\dblhopftwoside
\draw[double=none,Red,line width=2pt] (1) -- (5) (2) -- (8) (3) -- (4) (6) -- (7) ;
\end{tikzpicture}
\[
A^{-4}-A^{4}
\]

Unthreading the left\hyp{}hand hook results in a double Hopf chain of \(2n-2\) components.

\item
\tikzsetnextfilename{dblhopfdia4}
\autolabel
\begin{tikzpicture}[every path/.style={thick knot,double=Red},every node/.style={text=black},baseline=0cm]
\dblhopftwoside
\draw[double=none,Red,line width=2pt] (1) -- (5) (2) -- (6) (3) -- (4) (7) -- (8) ;
\end{tikzpicture}
\[
A^{-2}-2A^{6}+A^{10}
\]

We remove the central component (at a cost of \((-A^2 - A^{-2})\)), leaving a double Hopf chain with \(2n-4\) components.

\item
\tikzsetnextfilename{dblhopfdia5}
\autolabel
\begin{tikzpicture}[every path/.style={thick knot,double=Red},every node/.style={text=black},baseline=0cm]
\dblhopftwoside
\draw[double=none,Red,line width=2pt] (1) -- (2) (3) -- (4) (5) -- (6) (7) -- (8) ;
\end{tikzpicture}
\[
A^{-4}+1-3A^{4}+A^{8}
\]

We remove the two components (at a cost of \((-A^2 - A^{-2})^2\)), leaving a double Hopf chain with \(2n-4\) components.

\item
\tikzsetnextfilename{dblhopfdia6}
\autolabel
\begin{tikzpicture}[every path/.style={thick knot,double=Red},every node/.style={text=black},baseline=0cm]
\dblhopftwoside
\draw[double=none,Red,line width=2pt] (1) -- (2) (3) -- (7) (4) -- (8) (5) -- (6) ;
\end{tikzpicture}
\[
A^{-2}-2A^{6}+A^{10}
\]

We remove the central component (at a cost of \((-A^2 - A^{-2})\)), leaving a double Hopf chain with \(2n-4\) components.

\item
\tikzsetnextfilename{dblhopfdia7}
\autolabel
\begin{tikzpicture}[every path/.style={thick knot,double=Red},every node/.style={text=black},baseline=0cm]
\dblhopftwoside
\draw[double=none,Red,line width=2pt] (1) -- (2) (3) -- (5) (4) -- (8) (6) -- (7) ;
\end{tikzpicture}
\[
A^{-4}-A^{4}
\]

Unthreading the left\hyp{}hand hook results in a double Hopf chain of \(2n-2\) components.

\item
\tikzsetnextfilename{dblhopfdia8}
\autolabel
\begin{tikzpicture}[every path/.style={thick knot,double=Red},every node/.style={text=black},baseline=0cm]
\dblhopftwoside
\draw[double=none,Red,line width=2pt] (1) -- (5) (2) -- (3) (4) -- (6) (7) -- (8) ;
\end{tikzpicture}
\[
A^{-4}-A^{4}
\]

Unthreading the right\hyp{}hand hook results in a double Hopf chain of \(2n-2\) components.

\item
\tikzsetnextfilename{dblhopfdia9}
\autolabel
\begin{tikzpicture}[every path/.style={thick knot,double=Red},every node/.style={text=black},baseline=0cm]
\dblhopftwoside
\draw[double=none,Red,line width=2pt] (1) -- (4) (2) -- (3) (5) -- (8) (6) -- (7) ;
\end{tikzpicture}
\[
A^{-8}
\]

This is a double Hopf chain of \(2n\) components.

\item
\tikzsetnextfilename{dblhopfdia10}
\autolabel
\begin{tikzpicture}[every path/.style={thick knot,double=Red},every node/.style={text=black},baseline=0cm]
\dblhopftwoside
\draw[double=none,Red,line width=2pt] (1) -- (7) (2) -- (8) (3) -- (4) (5) -- (6) ;
\end{tikzpicture}
\[
A^{-2}-A^{2}
\]

We remove the central component (at a cost of \((-A^2 - A^{-2})\)), leaving a double Hopf chain with \(2n-4\) components.

\item
\tikzsetnextfilename{dblhopfdia11}
\autolabel
\begin{tikzpicture}[every path/.style={thick knot,double=Red},every node/.style={text=black},baseline=0cm]
\dblhopftwoside
\draw[double=none,Red,line width=2pt] (1) -- (7) (2) -- (3) (4) -- (8) (5) -- (6) ;
\end{tikzpicture}
\[
A^{-4}-A^{4}
\]

Unthreading the right\hyp{}hand hook results in a double Hopf chain of \(2n-2\) components.

\item
\tikzsetnextfilename{dblhopfdia12}
\autolabel
\begin{tikzpicture}[every path/.style={thick knot,double=Red},every node/.style={text=black},baseline=0cm]
\dblhopftwoside
\draw[double=none,Red,line width=2pt] (1) -- (4) (2) -- (3) (5) -- (6) (7) -- (8) ;
\end{tikzpicture}
\[
A^{-6}-A^{2}
\]

We remove the right\hyp{}hand component (at a cost of \((-A^2 - A^{-2})\)), leaving a double Hopf chain with \(2n-2\) components.

\item
\tikzsetnextfilename{dblhopfdia13}
\autolabel
\begin{tikzpicture}[every path/.style={thick knot,double=Red},every node/.style={text=black},baseline=0cm]
\dblhopftwoside
\draw[double=none,Red,line width=2pt] (1) -- (2) (3) -- (4) (5) -- (8) (6) -- (7) ;
\end{tikzpicture}
\[
A^{-6}-A^{2}
\]

We remove the left\hyp{}hand component (at a cost of \((-A^2 - A^{-2})\)), leaving a double Hopf chain with \(2n-2\) components.

\item
\tikzsetnextfilename{dblhopfdia14}
\autolabel
\begin{tikzpicture}[every path/.style={thick knot,double=Red},every node/.style={text=black},baseline=0cm]
\dblhopftwoside
\draw[double=none,Red,line width=2pt] (1) -- (5) (2) -- (6) (3) -- (7) (4) -- (8) ;
\end{tikzpicture}
\[
1-A^{4}-A^{8}+A^{12}
\]

This is a double Hopf ring of \(2n-2\) components.
\end{enumerate}

Thus, with the exception of Link~\ref{fig:dblhopfdia2} and Link~\ref{fig:dblhopfdia14}, we get a mixture of double Hopf chains of varying lengths.

\begin{enumerate}
\item Double Hopf chain with \(2n-4\) components: contributions from Link~\ref{fig:dblhopfdia1}, Link~\ref{fig:dblhopfdia4}, Link~\ref{fig:dblhopfdia5}, Link~\ref{fig:dblhopfdia6}, Link~\ref{fig:dblhopfdia10}.

\[
A^{-8}-A^{-4}-2+2A^{4}+A^{8}-A^{12}.
\]

\item Double Hopf chain with \(2n-2\) components: contributions from Link~\ref{fig:dblhopfdia3}, Link~\ref{fig:dblhopfdia7}, Link~\ref{fig:dblhopfdia8}, Link~\ref{fig:dblhopfdia11}, Link~\ref{fig:dblhopfdia12}, and Link~\ref{fig:dblhopfdia13}.

\[
-2A^{-8}+2A^{-4}+2-2A^{4}.
\]

\item Double Hopf chain with \(2n\) components: contribution from Link~\ref{fig:dblhopfdia9}.
\[
 A^{-8}.
\]

Using the recursion in \eqref{eq:dblhchain} we could rewrite this in terms of the double Hopf chains with \(2n-2\) and \(2n-4\) components.
\end{enumerate}

This leaves Link~\ref{fig:dblhopfdia2} to resolve.
It is halfway between the double Hopf chain and double Hopf ring: imagine bringing the two ends of the Hopf chain together and merging them as in Figure~\ref{fig:dblhopfmerge}, therefore we shall call it the double Hopf half\hyp{}ring.
Let us write \(\kb{dblhopfhalfring}{n}\) for the Kauffman bracket of this link with \(n\) components (note that \(n\) will be odd, and if the original double Hopf ring has \(2n\) components then in Link~\ref{fig:dblhopfdia2} we have \(2n-1\) components in this half\hyp{}ring).
To resolve this link, we look at the linkage between the ``different'' component and the double component to the right.
We shall assume that we are starting with \(2n+1\) components.
Our template is in Figure~\ref{fig:dblhopfhalftemplate}.
Here we have \(5\) complete components and thus the remainder consists of \(2n - 4\) components.

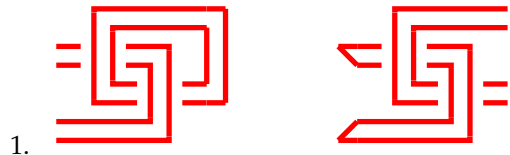
\tikzsetnextfilename{dblhopfhalftemplate}

\begin{figure}
\centering
\begin{tikzpicture}[every path/.style={thick knot,double=Red},every node/.style={text=black},baseline=0cm]
\pgfmathsetmacro{\brover}{.25}
\pgfmathsetmacro{\bradj}{.05}
\pgfmathsetmacro{\brsep}{.25}
\pgfmathsetmacro{\brlen}{1}
\foreach \brj in {0,...,2} {
\begin{scope}[xshift=3 * \brj cm]
\foreach \bri/\brdir/\brnud/\brv in {
  1/-1/2/-1,
  2/-1/1/-1,
  3/1/2/-1,
  4/1/1/-1,
  5/-1/1/1,
  6/-1/2/1,
  7/1/1/1,
  8/1/2/1
} {
\coordinate (\brj\bri) at (\brover * \brdir * \brnud,-\bri * \brsep + \brv * \bradj);
\draw (\brj\bri) -- ++(0,-\brv * \bradj) -- +(-\brdir * \brlen - \brover * \brdir * \brnud,0) coordinate (e\brj\bri);
}
\draw (\brj1) -- (\brj6);
\draw (\brj2) -- (\brj5);
\draw (\brj3) -- (\brj8);
\draw (\brj4) -- (\brj7);
\end{scope}
}
\draw[double=none,Red,line width=2pt]
 (e01) -- (e13)
 (e02) -- (e05)
 (e14) -- (e17)
 (e06) -- (e18)
 (e11) -- (e23)
 (e12) -- (e24)
 (e15) -- (e27)
 (e16) -- (e28)
;
\end{tikzpicture}
\caption{Template for the replacement diagrams in the double Hopf half\hyp{}ring}
\autolabel
\end{figure}

\tikzsetnextfilename{dblhopfmerge}

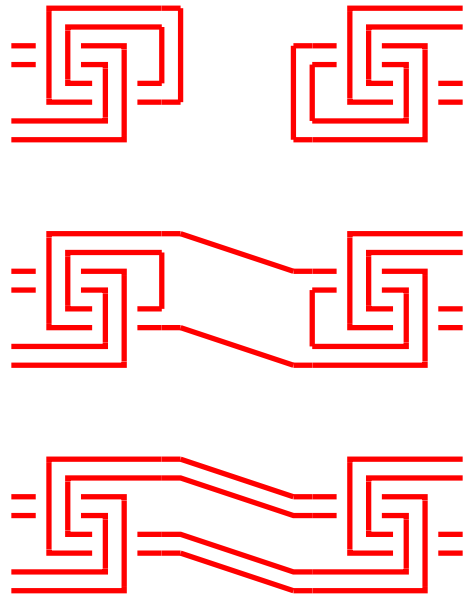
\begin{figure}
\centering
\begin{tikzpicture}[every path/.style={thick knot,double=Red},every node/.style={text=black},baseline=0cm]
\pgfmathsetmacro{\brover}{.25}
\pgfmathsetmacro{\bradj}{.05}
\pgfmathsetmacro{\brsep}{.25}
\pgfmathsetmacro{\brlen}{1}
\dblhopftwoside
\draw[double=none,Red,line width=2pt]
 (1) -- (4)
 (2) -- (3)
 (5) -- (8)
 (6) -- (7)
;
\begin{scope}[yshift=-3cm]
\dblhopftwoside
\draw[double=none,Red,line width=2pt]
 (1) -- (5)
 (2) -- (3)
 (4) -- (8)
 (6) -- (7)
;
\begin{scope}[yshift=-3cm]
\dblhopftwoside
\draw[double=none,Red,line width=2pt]
 (1) -- (5)
 (2) -- ++(\brsep,0) -- ($(6)-(\brsep,0)$) -- (6)
 (3) -- ++(\brsep,0) -- ($(7)-(\brsep,0)$) -- (7)
 (4) -- (8)
;
\end{scope}
\end{scope}
\end{tikzpicture}
\caption{Merging a double Hopf chain into a double Hopf ring.}
\autolabel
\end{figure}

\begin{enumerate}
\pgfmathsetmacro{\brover}{.25}
\pgfmathsetmacro{\bradj}{.05}
\pgfmathsetmacro{\brsep}{.25}
\pgfmathsetmacro{\brlen}{1}

\item
\tikzsetnextfilename{dblhopfrcdia1}
\autolabel
\begin{tikzpicture}[every path/.style={thick knot,double=Red},every node/.style={text=black},baseline=0cm]
\dblhopfrctwoside
\draw[double=none,Red,line width=2pt] (1) -- (4) (5) -- (6) (7) -- (8) ;
\end{tikzpicture}
\[
-A^{-8}+2-2A^{4}+A^{8}
\]

After removing the circle, this is the double Hopf chain with \(2n-2\) components.

\item
\tikzsetnextfilename{dblhopfrcdia2}
\autolabel
\begin{tikzpicture}[every path/.style={thick knot,double=Red},every node/.style={text=black},baseline=0cm]
\dblhopfrctwoside
\draw[double=none,Red,line width=2pt] (1) -- (4) (5) -- (8) (6) -- (7) ;
\end{tikzpicture}
\[
-A^{-10}-A^{2}
\]

This is the double Hopf chain with \(2n\) components.

\item
\tikzsetnextfilename{dblhopfrcdia3}
\autolabel
\begin{tikzpicture}[every path/.style={thick knot,double=Red},every node/.style={text=black},baseline=0cm]
\dblhopfrctwoside
\draw[double=none,Red,line width=2pt] (1) -- (5) (4) -- (6) (7) -- (8) ;
\end{tikzpicture}
\[
-A^{-6}+A^{-2}-A^{6}+A^{10}
\]

This is the double Hopf chain with \(2n-2\) components.

\item
\tikzsetnextfilename{dblhopfrcdia4}
\autolabel
\begin{tikzpicture}[every path/.style={thick knot,double=Red},every node/.style={text=black},baseline=0cm]
\dblhopfrctwoside
\draw[double=none,Red,line width=2pt] (1) -- (7) (4) -- (8) (5) -- (6) ;
\end{tikzpicture}
\[
-A^{-6}+A^{-2}-A^{6}+A^{10}
\]

This is the double Hopf chain with \(2n-2\) components.

\item
\tikzsetnextfilename{dblhopfrcdia5}
\autolabel
\begin{tikzpicture}[every path/.style={thick knot,double=Red},every node/.style={text=black},baseline=0cm]
\dblhopfrctwoside
\draw[double=none,Red,line width=2pt] (1) -- (5) (4) -- (8) (6) -- (7) ;
\end{tikzpicture}
\[
-A^{-8}+1-A^{4}+A^{12}
\]

This is the same as we started with, but with two fewer components (thus \(2n -1\)).

\end{enumerate}

Gathering together terms, we find that:
\begin{dmath*}
  \kb{dblhopfhalfring}{2n+1} = (A^{-10}-A^{-6}-A^{6}+A^{10})\kb{dblhopfchain}{2n} - (A^{-10} + A^2) \kb{dblhopfchain}{2n-2} + (-A^{-8} + 1 - A^4 + A^{12}) \kb{dblhopfhalfring}{2n-1}.
\end{dmath*}

We need to establish a starting point.
In the template for our analysis, we had the special linkage flanked on either side by ordinary double Hopf linkages.
We can safely assume that those two double Hopf linkages are actually the same Hopf linkage.
This means that our recursive formula holds for \(2n+1 = 5\) and so we need to compute \(\kb{dblhopfhalfring}{3}\).
A simple calculation shows that:
\[
  \kb{dblhopfhalfring}{3} = A^{16} + A^8 + 2.
\]

Returning to the double Hopf ring, we have the recursive formula:
\begin{dmath*}
\kb{dblhopfring}{2n} = A^{-8} \kb{dblhopfchain}{2n} + (-2 A^{-8} + 2 A^{-4} + 2 - 2 A^4) \kb{dblhopfchain}{2n-2} + (A^{-8} - A^{-4} -2 + 2 A^4 + A^8 - A^{12}) \kb{dblhopfchain}{2n-4} + (A^{-6} + A^{-2} - A^2 - A^6) \kb{dblhopfhalfring}{2n-1} + (1 - A^4 - A^8 + A^{12}) \kb{dblhopfring}{2n-2}.
\end{dmath*}

Again, we need to establish a starting point and again, we can assume that in our analysis the flanking linkages were the same.
With the convention that \(\kb{dblhopfchain}{0} = (-A^2 - A^{-2})^{-1}\), our recursive formula still holds for \(\kb{dblhopfring}{4}\).
That is to say,
\begin{dmath*}
\kb{dblhopfring}{4} = A^{-8} \kb{dblhopfchain}{4} + (-2 A^{-8} + 2 A^{-4} + 2 - 2 A^4) \kb{dblhopfchain}{2} + (A^{-8} - A^{-4} -2 + 2 A^4 + A^8 - A^{12}) \kb{dblhopfchain}{0} + (A^{-6} + A^{-2} - A^2 - A^6) \kb{dblhopfhalfring}{3} + (1 - A^4 - A^8 + A^{12}) \kb{dblhopfring}{2}.
\end{dmath*}

A simple calculation shows that:
\[
\kb{dblhopfring}{2} = - A^{18} - A^{10} + A^{6} - A^{2}.
\]
Using this as our starting point would mean that in the recursive formula for the double hopf ring, we would have to use \(\kb{dblhopfchain}{0}\).
We would prefer not to have this, and so we compute the first iteration as well to see that:
\[
  \kb{dblhopfring}{4} = - A^{30} + A^{26} - 2 A^{18} - A^{14} + A^{10} + A^{6} - 2 A^{2} - 3 A^{-2} - A^{-14} - A^{-22}.
\]

Finally, we return to the level two Hopf chain and ring and observe that the double Hopf chain and ring relate to the level two structures in the following way:
\begin{align*}
  \kb{ltwohopfchain}{0,\dotsc,0} &= \kb{dblhopfchain}{2 n} \\
  \kb{ltwohopfring}{0,\dotsc,0} &= \kb{dblhopfring}{2 n}
\end{align*}
where in each case there are \(n\) zeros in the index.

\section{The Brunnian Link}
\label{sec:brlink}

Now we turn to Brunnian linkages.
A Brunnian linkage between two circles is shown in Figure~\ref{fig:brunnianlink}.
As it stands, this is unlinked.
To make it linked, one has to consider this as part of a larger diagram.
With only pure Brunnian linkages, one of the simplest such diagrams is the ring\hyp{}like link in Figure~\ref{fig:brunnianring}.

\tikzsetnextfilename{brunnianlink}

\begin{figure}
\centering
\begin{tikzpicture}
\pgfmathsetmacro{\brscale}{1}
\foreach \brk in {1,...,2} {
\pgfmathsetmacro{\brtint}{(\brk > 4 ? 8 - \brk : \brk) * 25}
\pgfmathsetmacro{\brl}{\brscale * \brk}
\begin{scope}[xshift=\brl cm]
\colorlet{chain}{Red!\brtint!Green}
\flatbrunnianlink{\brscale}
\end{scope}
}
\end{tikzpicture}
\caption{The Brunnian linkage}
\autolabel
\end{figure}
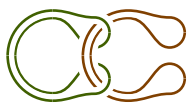

\tikzsetnextfilename{brunnianring}
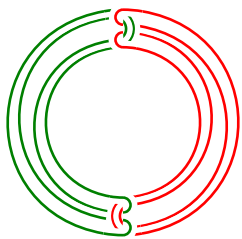
\begin{figure}
\centering
\begin{tikzpicture}
\setbrstep{.1}
\foreach \brk in {1,2} {
\begin{scope}[rotate=\brk * 180]
\colorlet{chain}{ring\brk}
\brunnianlink{1}{180}
\end{scope}
}
\end{tikzpicture}
\caption{The Brunnian ring with two components}
\autolabel
\end{figure}

It is also possible to have a ``half\hyp{}Brunnian'' linkage, as in Figure~\ref{fig:halfbrunnian}.
Again, this is unlinked unless there is a larger diagram.
The simplest case now is with a ``half\hyp{}Brunnian'' linkage at either end, as in Figure~\ref{fig:brunnianchain}.
This is, incidentally, isotopic to the Borromean rings.
From this, one can add more links in the middle to form a \emph{Brunnian chain}, as in Figure~\ref{fig:longbrunnianchain}.

\tikzsetnextfilename{halfbrunnian}
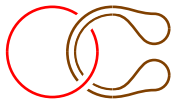
\begin{figure}
\centering
\begin{tikzpicture}
\colorlet{chain}{Red!50!Green}
\begin{pgfonlayer}{back}
\draw[knot,double=Red] (-.2,0) arc(0:180:.6);
\end{pgfonlayer}
\flatbrunnianlink{1}
\draw[knot,double=Red] (-.2,0) arc(0:-180:.6);
\end{tikzpicture}
\caption{Half\hyp{}Brunnian linkage}
\autolabel
\end{figure}

\tikzsetnextfilename{brunnianchain}
\begin{figure}
\centering
\begin{tikzpicture}
\colorlet{chain}{Red!50!Green}
\begin{pgfonlayer}{back}
\draw[knot,double=Red] (-.2,0) arc(0:180:.6);
\end{pgfonlayer}
\flatbrunnianlink{1}
\draw[knot,double=Red] (-.2,0) arc(0:-180:.6);
\draw[knot,double=Green] (1,0) circle (.6);
\end{tikzpicture}
\caption{Brunnian chain}
\autolabel
\end{figure}

\tikzsetnextfilename{longbrunnianchain}
\begin{figure}
\centering
\begin{tikzpicture}
\colorlet{chain}{Red!50!Green}
\begin{pgfonlayer}{back}
\draw[knot,double=Red] (.8,0) arc(0:180:.6);
\end{pgfonlayer}
\pgfmathsetmacro{\brscale}{1}
\foreach \brk in {1,...,8} {
\pgfmathsetmacro{\brtint}{\brk * 11}
\pgfmathsetmacro{\brl}{\brscale * \brk}
\begin{scope}[xshift=\brl cm]
\colorlet{chain}{Green!\brtint!Red}
\flatbrunnianlink{\brscale}
\end{scope}
}
\draw[knot,double=Red] (.8,0) arc(0:-180:.6);
\draw[knot,double=Green] (9,0) circle (.6);
\end{tikzpicture}
\caption{Brunnian chain with 8 middle components}
\autolabel
\end{figure}
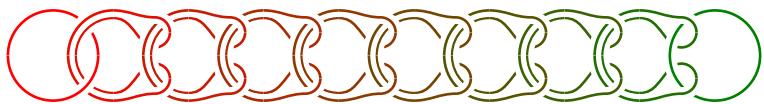

We already know how to start with the Brunnian chain since the linkage between the first two components is the ``one\hyp{}and\hyp{}a\hyp{}half Hopf'' linkage from Figure~\ref{fig:onehalfHopfres}.
This therefore resolves into two diagrams where the first circle is removed.
In the first diagram, which has multiplier \(-A^6 - A^{-6}\), the diagram is the same as the original except for the removal of the first circle.
This diagram is now unlinked and so is \(n-1\) disjoint circles, where \(n\) is the number of components in the original link.
Removing each circle (save for the last) contributes a factor of \((-A^2 - A^{-2})\), resulting in \((-A^6 - A^{-6})(-A^2 - A^{-2})^{n-2}\).
The second diagram, which has multiplier \(2 - A^{4} - A^{-4}\), is formed by capping off the double strand.
This produces a link as in Figure~\ref{fig:brchain2end}.
Note that although we lose the extreme component, we also split the next component in to two and thus have the same number of components as at the start.

\tikzsetnextfilename{brchain2end}
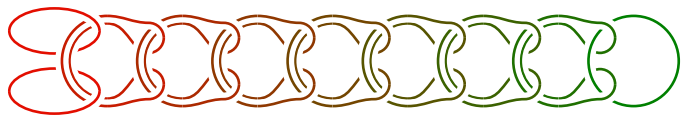
\begin{figure}
\centering
\begin{tikzpicture}
\begin{pgfonlayer}{back}
\draw[knot,double=Green!11!Red] (1.9,.4) arc(0:-180:.6 and .3);
\draw[knot,double=Green!11!Red] (1.9,-.4) arc(0:180:.6 and .3);
\end{pgfonlayer}
\pgfmathsetmacro{\brscale}{1}
\foreach \brk in {2,...,8} {
\pgfmathsetmacro{\brtint}{\brk * 11}
\pgfmathsetmacro{\brl}{\brscale * \brk}
\begin{scope}[xshift=\brl cm]
\colorlet{chain}{Green!\brtint!Red}
\flatbrunnianlink{\brscale}
\end{scope}
}
\draw[knot,double=Green!11!Red] (1.9,.4) arc(0:180:.6 and .3);
\draw[knot,double=Green!11!Red] (1.9,-.4) arc(0:-180:.6 and .3);
\draw[knot,double=Green] (9,0) circle (.6);
\end{tikzpicture}
\caption{Brunnian chain with 7 middle components and a double end loop}
\autolabel
\end{figure}

Let us write \(\kb{brunnianchain}{n}\) for the Kauffman bracket of a Brunnian chain with \(n\) components.
Let us write \(\kb{brunnianchaindblleft}{n}\) for the Kauffman bracket of a Brunnian chain with \(n\) components where there are \(2\) rings at the end as in Figure~\ref{fig:brchain2end}.
Then we have shown that we have the following partial recurrence relation:
\[
  \kb{brunnianchain}{n} = (-A^6 - A^{-6})(-A^2 - A^{-2})^{n-2} + (2 - A^{4} - A^{-4}) \kb{brunnianchaindblleft}{n}.
\]

Looking at the Brunnian chain with the double circle, we can resolve this using the computer program.
We feed in the entire left\hyp{}hand end and obtain two diagrams at the conclusion, as shown in Figure~\ref{fig:brunniancappedend}.
The link corresponding to the first diagram unlinks leaving \(n-2\) circles.
The link corresponding to the second diagram is the Brunnian chain with the double circle with one fewer component.
Hence:
\[
  \kb{brunnianchaindblleft}{n} = (A^{12} + 2 + A^{-12})(-A^2 - A^{-2})^{n-3} + (A^{10} - A^{6} - A^{-6} + A^{-10})   \kb{brunnianchaindblleft}{n-1}.
\]

\tikzsetnextfilename{brunniancappedend}

\begin{figure}
\centering
\begin{tikzpicture}[every path/.style={thick knot,double=Red},every node/.style={text=black}]
\pgfmathsetmacro{\brover}{.5}
\pgfmathsetmacro{\bradj}{.1}
\pgfmathsetmacro{\brlen}{1}
\foreach \bri/\brdir/\brnud/\brv in {
  1/1/1/-1,
  2/-1/1/-1,
  3/-1/0/-1,
  4/1/1/1,
  5/1/1/-1,
  6/-1/0/1,
  7/-1/1/1,
  8/1/1/1
} {
\coordinate (\bri) at (\brover * \brdir * \brnud,-\bri/2 + \brv * \bradj);
\draw (\bri) -- ++(0,-\brv * \bradj) -- +(-\brdir * \brlen - \brover * \brdir * \brnud,0) coordinate (e\bri);
}
\draw (1) -- (4);
\draw (2) -- (7);
\draw (3) -- (6);
\draw (5) -- (8);
\draw (e1) -- ++(-.2,0) |- (e4);
\draw (e5) -- ++(-.2,0) |- (e8);
\draw[double=none,line width=2pt,black,->] (1.5,-2.25) -- ++(1,0);
\begin{scope}[xshift=3cm]
\foreach \bri/\brdir/\brnud/\brv in {
  2/-1/0/-1,
  3/-1/-1/-1,
  6/-1/-1/1,
  7/-1/0/1
} {
\coordinate (\bri) at (\brover * \brdir * \brnud,-\bri/2 + \brv * \bradj);
\draw (\bri) -- ++(0,-\brv * \bradj) -- +(-\brdir * \brlen - \brover * \brdir * \brnud,0) coordinate (e\bri);
}
\draw (2) -- (7);
\draw (3) -- (6);
\node at (2,-2.25) {and};
\begin{scope}[xshift=3cm]
\foreach \bri/\brdir/\brnud/\brv in {
  2/-1/0/-1,
  3/-1/0/1,
  6/-1/0/-1,
  7/-1/0/1
} {
\coordinate (\bri) at (\brover * \brdir * \brnud,-\bri/2 + \brv * \bradj);
\draw (\bri) -- ++(0,-\brv * \bradj) -- +(-\brdir * \brlen - \brover * \brdir * \brnud,0) coordinate (e\bri);
}
\draw (2) -- (3);
\draw (7) -- (6);
\end{scope}
\end{scope}
\end{tikzpicture}
\caption{The end of the capped Brunnian chain}
\autolabel
\end{figure}
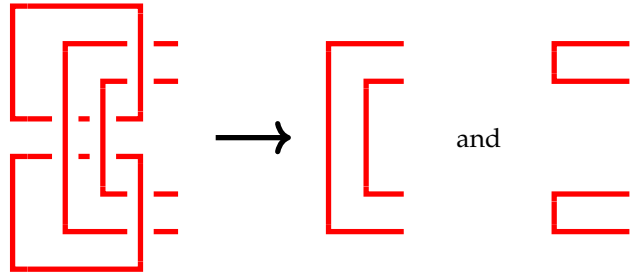

For our starting points, we note that \(\kb{brunnianchain}{2} = \kb{hopfchain}{2} = -A^4 - A^{-4}\) and \(\kb{brunnianchaindblleft}{3} = \kb{hopfchain}{3} = (-A^4 - A^{-4})^2\).

Now let us consider the Brunnian rings wherein the linkages are all Brunnian linkages.
Figure~\ref{fig:brunnianring} contains the Brunnian ring with \(2\) components.
Figure~\ref{fig:brunnianring345} contains the Brunnian rings with  \(3\), \(4\), and \(5\) components.
Let us write \(\kb{brunnianring}{n}\) for the Kauffman bracket of the Brunnian ring with \(n\) components.

\tikzsetnextfilename{brunnianring345}
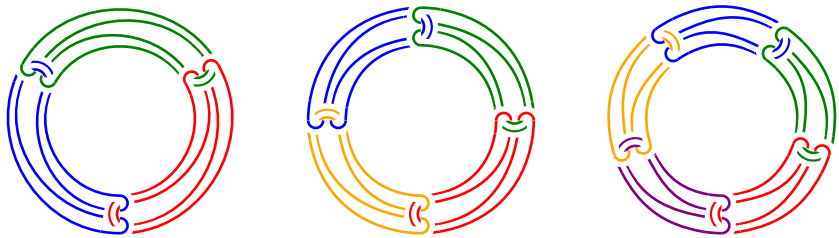
\begin{figure}
\centering
\begin{tikzpicture}
\setbrstep{.1}
\foreach \brl in {3,4,5} {
\begin{scope}[xshift=\brl * 4cm]
\pgfmathsetmacro{\brang}{360/\brl}
\foreach \brk in {1,...,\brl} {
\begin{scope}[rotate=\brk * \brang]
\colorlet{chain}{ring\brk}
\brunnianlink{1}{\brang}
\end{scope}
}
\end{scope}
}
\end{tikzpicture}
\caption{The Brunnian rings with three, four, and five components}
\autolabel
\end{figure}

The Brunnian linkage can be drawn schematically as in Figure~\ref{fig:schbrunnian}.
As there are eight entry\hyp{}exit strands, there will be (at most) \(14\) diagrams after resolving this linkage.
To find the contributions of each diagram, we use the computer program to scan through the possibilities.
The labelled diagram is Figure~\ref{fig:schbrlabelled}.

Note that one possible diagram does not actually occur.
The diagram wherein \(1\) is connected to \(4\), \(2\) to \(3\), \(5\) to \(8\), and \(6\) to \(7\) cannot be obtained from the Brunnian linkage.
The reason is that only two strands can cross the ``half\hyp{}way'' point but that diagram requires four.

The results from the other thirteen diagrams follow.
As before, we show the diagram flanked by Brunnian linkages on either side.
Starting with a Brunnian ring with \(n\) components, the part of the diagram that we can see originally had two full components and four parts of two more components.
The partial strands do not alter under this replacement process, and thus the number of components in one of the following links is the number of full components visible plus \(n-2\).

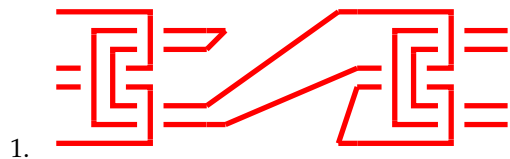
\tikzsetnextfilename{schbrunnian}

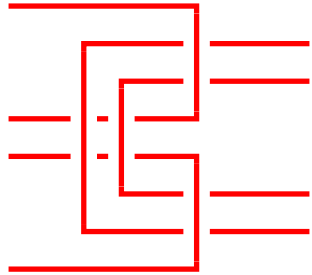
\begin{figure}
\centering
\begin{tikzpicture}[every path/.style={thick knot,double=Red}]
\pgfmathsetmacro{\brover}{.5}
\pgfmathsetmacro{\bradj}{.1}
\foreach \bri/\brdir/\brnud/\brv in {1/1/1/-1,2/-1/2/-1,3/-1/1/-1,4/1/1/1,5/1/1/-1,6/-1/1/1,7/-1/2/1,8/1/1/1} {
\coordinate (\bri) at (\brover * \brdir * \brnud,-\bri/2 + \brv * \bradj);
\draw (\bri) -- ++(0,-\brv * \bradj) -- +(-\brdir * 2 - \brover * \brdir * \brnud,0);
}
\draw (1) -- (4);
\draw (2) -- (7);
\draw (3) -- (6);
\draw (5) -- (8);
\end{tikzpicture}
\caption{Schematic version of the Brunnian linkage}
\autolabel
\end{figure}

\tikzsetnextfilename{schbrlabelled}

\begin{figure}
\centering
\begin{tikzpicture}[every path/.style={thick knot,double=Red},every node/.style={text=black}]
\pgfmathsetmacro{\brover}{.5}
\pgfmathsetmacro{\bradj}{.1}
\foreach \bri/\brdir/\brnud/\brv in {
  1/1/1/-1,
  2/-1/2.5/-1,
  3/-1/.5/-1,
  4/1/1/1,
  5/1/1/-1,
  6/-1/.5/1,
  7/-1/2.5/1,
  8/1/1/1
} {
\coordinate (\bri) at (\brover * \brdir * \brnud,-\bri/2 + \brv * \bradj);
\draw (\bri) -- ++(0,-\brv * \bradj) -- +(-\brdir * 2 - \brover * \brdir * \brnud,0) coordinate (e\bri);
}
\draw (1) -- (4);
\draw (2) -- (7);
\draw (3) -- (6);
\draw (5) -- (8);
\path
node at (e1) {1}
node at (e4) {2}
node at (e5) {3}
node at (e8) {4}
node at (e2) {5}
node at (e3) {6}
node at (e6) {7}
node at (e7) {8}
node at (2) {9}
node at (3) {10}
node at (4) {11}
node at (5) {12}
node at (6) {13}
node at (7) {14}
node at ($(1)!.5!(4)$) {15}
node at ($(2)!.5!(7)$) {16}
node at ($(3)!.5!(6)$) {17}
node at ($(5)!.5!(8)$) {18}
node at ($(4)!.5!(e4)$) {19}
node at ($(5)!.5!(e5)$) {20}
;
\end{tikzpicture}
\caption{The labelled Brunnian linkage}
\autolabel
\end{figure}
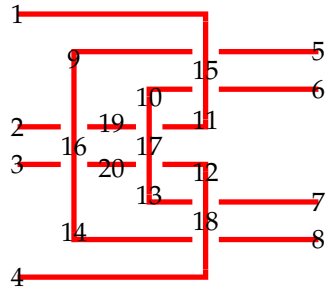

\begin{enumerate}
\pgfmathsetmacro{\brover}{.25}
\pgfmathsetmacro{\bradj}{.05}
\pgfmathsetmacro{\brsep}{.25}
\pgfmathsetmacro{\brlen}{1}

\item
\tikzsetnextfilename{brdia1}
\autolabel
\begin{tikzpicture}[every path/.style={thick knot,double=Red},every node/.style={text=black},baseline=0cm]
\brunniantwoside
\draw[double=none,Red,line width=2pt] (1) -- (2) (3) -- (5) (4) -- (6) (7) -- (8) ;
\end{tikzpicture}
\begin{dmath*}
A^4 - A^8
\end{dmath*}

This diagram unlinks.
The left\hyp{}hand link clearly unlinks by threading the upper spike back through the loops.
As the part of the diagram that is not shown consists entirely of Brunnian linkages, once one is unlinked the rest fall apart also.
This leads round to the linkage on the right, which also comes apart.
The result is \(n\) unlinked circles. 

\item
\tikzsetnextfilename{brdia2}
\autolabel
\begin{tikzpicture}[every path/.style={thick knot,double=Red},every node/.style={text=black},baseline=0cm]
\brunniantwoside
\draw[double=none,Red,line width=2pt] (1) -- (5) (2) -- (3) (4) -- (8) (6) -- (7) ;
\end{tikzpicture}
\begin{dmath*}
-A^{-8}+2-A^{8}
\end{dmath*}

This diagram unlinks, starting on the right.
At the end, there are \(n\) unlinked circles.

\item
\tikzsetnextfilename{brdia3}
\autolabel
\begin{tikzpicture}[every path/.style={thick knot,double=Red},every node/.style={text=black},baseline=0cm]
\brunniantwoside
\draw[double=none,Red,line width=2pt] (1) -- (5) (2) -- (8) (3) -- (4) (6) -- (7) ;
\end{tikzpicture}
\begin{dmath*}
-A^{-10}+A^{-2}
\end{dmath*}

This diagram unlinks, starting on both sides.
At the end, there are \((n-1)\) unlinked circles.

\item
\tikzsetnextfilename{brdia4}
\autolabel
\begin{tikzpicture}[every path/.style={thick knot,double=Red},every node/.style={text=black},baseline=0cm]
\brunniantwoside
\draw[double=none,Red,line width=2pt] (1) -- (5) (2) -- (6) (3) -- (4) (7) -- (8) ;
\end{tikzpicture}
\begin{dmath*}
-A^{-8}+2-A^{4}
\end{dmath*}

This diagram unlinks, starting on the left.
At the end, there are \(n\) unlinked circles.

\item
\tikzsetnextfilename{brdia5}
\autolabel
\begin{tikzpicture}[every path/.style={thick knot,double=Red},every node/.style={text=black},baseline=0cm]
\brunniantwoside
\draw[double=none,Red,line width=2pt] (1) -- (2) (3) -- (4) (5) -- (6) (7) -- (8) ;
\end{tikzpicture}
\begin{dmath*}
-A^{-6}+A^{-2}+A^{2}-A^{6}
\end{dmath*}

This diagram unlinks, starting on the left.
At the end, there are \((n+1)\) unlinked circles.

\item
\tikzsetnextfilename{brdia6}
\autolabel
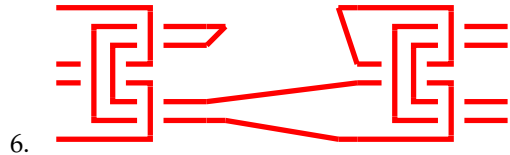
\begin{tikzpicture}[every path/.style={thick knot,double=Red},every node/.style={text=black},baseline=0cm]
\brunniantwoside
\draw[double=none,Red,line width=2pt] (1) -- (2) (3) -- (7) (4) -- (8) (5) -- (6) ;
\end{tikzpicture}
\begin{dmath*}
-A^{-4}+2-A^{8}
\end{dmath*}

This diagram unlinks, starting on the left.
At the end, there are \(n\) unlinked circles.

\item
\tikzsetnextfilename{brdia7}
\autolabel
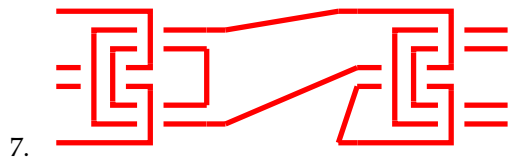
\begin{tikzpicture}[every path/.style={thick knot,double=Red},every node/.style={text=black},baseline=0cm]
\brunniantwoside
\draw[double=none,Red,line width=2pt] (1) -- (5) (2) -- (3) (4) -- (6) (7) -- (8) ;
\end{tikzpicture}
\begin{dmath*}
-A^{-6}+A^{-2}+A^{2}-A^{6}
\end{dmath*}

This link does not unlink and is not one that we have seen before, so we postpone its analysis for the moment.
Note that we actually have \(n+1\) components in this link.

\item
\tikzsetnextfilename{brdia8}
\autolabel
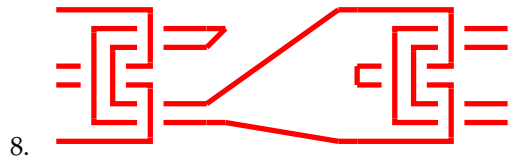
\begin{tikzpicture}[every path/.style={thick knot,double=Red},every node/.style={text=black},baseline=0cm]
\brunniantwoside
\draw[double=none,Red,line width=2pt] (1) -- (2) (3) -- (5) (4) -- (8) (6) -- (7) ;
\end{tikzpicture}
\begin{dmath*}
A^{2}-A^{10}
\end{dmath*}

This diagram unlinks, starting on the left.
At the end, there are \((n-1)\) unlinked circles.

\item
\tikzsetnextfilename{brdia9}
\autolabel
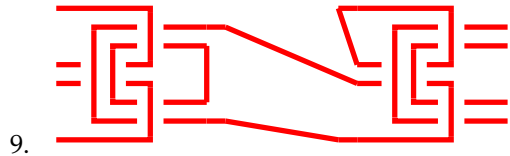
\begin{tikzpicture}[every path/.style={thick knot,double=Red},every node/.style={text=black},baseline=0cm]
\brunniantwoside
\draw[double=none,Red,line width=2pt] (1) -- (7) (2) -- (3) (4) -- (8) (5) -- (6) ;
\end{tikzpicture}
\begin{dmath*}
-A^{-6}+A^{-2}+A^{2}-A^{6}
\end{dmath*}

This link does not unlink and is not one that we have seen before, so we postpone its analysis for the moment.
Note that we actually have \(n+1\) components in this link.

\item
\tikzsetnextfilename{brdia10}
\autolabel
\begin{tikzpicture}[every path/.style={thick knot,double=Red},every node/.style={text=black},baseline=0cm]
\brunniantwoside
\draw[double=none,Red,line width=2pt] (1) -- (7) (2) -- (8) (3) -- (4) (5) -- (6) ;
\end{tikzpicture}
\begin{dmath*}
-A^{-8}+A^{-4}
\end{dmath*}

This diagram unlinks, starting on the left.
At the end, there are \(n\) unlinked circles.

\item
\tikzsetnextfilename{brdia11}
\autolabel
\begin{tikzpicture}[every path/.style={thick knot,double=Red},every node/.style={text=black},baseline=0cm]
\brunniantwoside
\draw[double=none,Red,line width=2pt] (1) -- (4) (2) -- (3) (5) -- (6) (7) -- (8) ;
\end{tikzpicture}
\begin{dmath*}
-A^{-4}+2-A^{4}
\end{dmath*}

This is like the modified Brunnian chain, except that there are double loops at \emph{both} ends, not just one.

\item
\tikzsetnextfilename{brdia12}
\autolabel
\begin{tikzpicture}[every path/.style={thick knot,double=Red},every node/.style={text=black},baseline=0cm]
\brunniantwoside
\draw[double=none,Red,line width=2pt] (1) -- (2) (3) -- (4) (5) -- (8) (6) -- (7) ;
\end{tikzpicture}
\begin{dmath*}
 1 
\end{dmath*}

This unlinks, starting on the left.
At the end, there are \(n\) unlinked circles.

\item
\tikzsetnextfilename{brdia13}
\autolabel
\begin{tikzpicture}[every path/.style={thick knot,double=Red},every node/.style={text=black},baseline=0cm]
\brunniantwoside
\draw[double=none,Red,line width=2pt] (1) -- (5) (2) -- (6) (3) -- (7) (4) -- (8) ;
\end{tikzpicture}
\begin{dmath*}
-A^{-6}+A^{-2}+A^{2}-A^{6}
\end{dmath*}

This is a Brunnian ring of one fewer components.

\end{enumerate}

Many of the diagrams unlink to either \((n-1)\), \(n\), or \((n+1)\) unlinked circles.
Taking \((n-1)\) as our ``base'', we can gather together all of these terms.
This results in:

\begin{dmath*}
(A^{-10}+2A^{-6}-4A^{-2}-4A^{2}+2A^{6}+A^{10})(-A^2-A^{-2})^{n-2} = \big(-(-A^{-2} - A^2)^4 + 3(-A^{-2} - A^2)^2 + 5\big)(-A^2-A^{-2})^{n-1}
\end{dmath*}

This leaves the links in \ref{fig:brdia7}, \ref{fig:brdia9}, \ref{fig:brdia11}, and \ref{fig:brdia13}.
The last of these is the Brunnian ring with one fewer components; the penultimate one is the Brunnian chain with double rings at both ends.
Thus we are left with the similar\hyp{}looking \ref{fig:brdia7} and \ref{fig:brdia9}.
The left\hyp{}hand linkage is the same in both diagrams, so we feed that back in to our computer program.
Let us, for the sake of definiteness, consider \ref{fig:brdia7}.
We assume that we start with \(n\) components.
Figure~\ref{fig:brunnianhalfseg} is the segment of the link diagram that we are considering.
On the left\hyp{}hand side we have an ordinary Brunnian linkage.
On the right\hyp{}hand side we have the right\hyp{}hand linkage from link \ref{fig:brdia7} above.
The remainder of the diagram consists of Brunnian linkages.
Starting with \(n\) components, we can see \(4\) complete components.
Thus the rest of the diagram consists of \(n - 4\) components.

\tikzsetnextfilename{brunnianhalfseg}

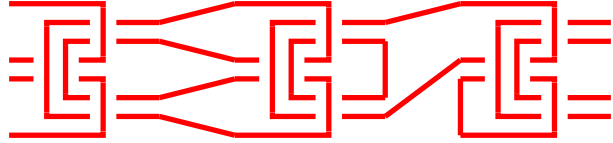
\begin{figure}
\centering
\begin{tikzpicture}[every path/.style={thick knot,double=Red},every node/.style={text=black},baseline=0cm]
\pgfmathsetmacro{\brover}{.25}
\pgfmathsetmacro{\bradj}{.05}
\pgfmathsetmacro{\brsep}{.25}
\pgfmathsetmacro{\brlen}{1}
\foreach \bri/\brdir/\brnud/\brv in {
  1/1/1/-1,
  2/-1/2/-1,
  3/-1/1/-1,
  4/1/1/1,
  5/1/1/-1,
  6/-1/1/1,
  7/-1/2/1,
  8/1/1/1
} {
  \coordinate (l\bri) at (\brover * \brdir * \brnud,-\bri * \brsep + \brv * \bradj);
  \draw (l\bri) -- ++(0,-\brv * \bradj) -- +(-\brdir * \brlen - \brover * \brdir * \brnud,0) coordinate (le\bri);
}
\draw (l1) -- (l4);
\draw (l2) -- (l7);
\draw (l3) -- (l6);
\draw (l5) -- (l8);
\begin{scope}[xshift=3cm]
\foreach \bri/\brdir/\brnud/\brv in {
  1/1/1/-1,
  2/-1/2/-1,
  3/-1/1/-1,
  4/1/1/1,
  5/1/1/-1,
  6/-1/1/1,
  7/-1/2/1,
  8/1/1/1
} {
  \coordinate (m\bri) at (\brover * \brdir * \brnud,-\bri * \brsep + \brv * \bradj);
  \draw (m\bri) -- ++(0,-\brv * \bradj) -- +(-\brdir * \brlen - \brover * \brdir * \brnud,0) coordinate (me\bri);
}
\draw (m1) -- (m4);
\draw (m2) -- (m7);
\draw (m3) -- (m6);
\draw (m5) -- (m8);
\begin{scope}[xshift=3cm]
\foreach \bri/\brdir/\brnud/\brv in {
  1/1/1/-1,
  2/-1/2/-1,
  3/-1/1/-1,
  4/1/1/1,
  5/1/1/-1,
  6/-1/1/1,
  7/-1/2/1,
  8/1/1/1
} {
  \coordinate (r\bri) at (\brover * \brdir * \brnud,-\bri * \brsep + \brv * \bradj);
  \draw (r\bri) -- ++(0,-\brv * \bradj) -- +(-\brdir * \brlen - \brover * \brdir * \brnud,0) coordinate (re\bri);
}
\draw (r1) -- (r4);
\draw (r2) -- (r7);
\draw (r3) -- (r6);
\draw (r5) -- (r8);
\end{scope}
\end{scope}
\draw[double=none,Red,line width=2pt]
 (le2) -- (me1)
 (le3) -- (me4)
 (le6) -- (me5)
 (le7) -- (me8)
 (me2) -- (re1)
 (me3) -- (me6)
 (me7) -- (re4)
 (re5) -- (re8)
;
\end{tikzpicture}
\caption{The template for the half\hyp{}Brunnian link}
\autolabel
\end{figure}

\begin{enumerate}

\pgfmathsetmacro{\brover}{.25}
\pgfmathsetmacro{\bradj}{.05}
\pgfmathsetmacro{\brsep}{.25}
\pgfmathsetmacro{\brlen}{1}

\item
\tikzsetnextfilename{brunnianhalfdia1}
\autolabel
\begin{tikzpicture}[every path/.style={thick knot,double=Red},every node/.style={text=black},baseline=0cm]
\brunnianhalftwoside
\draw[double=none,Red,line width=2pt] (1) -- (2) (3) -- (4) (5) -- (8) ;
\end{tikzpicture}

\[
- A^{6} - A^{-6}
\]

Once the left\hyp{}hand component is removed, the rest of the diagram unlinks leaving \(n-1\) components.

\item
\tikzsetnextfilename{brunnianhalfdia2}
\autolabel
\begin{tikzpicture}[every path/.style={thick knot,double=Red},every node/.style={text=black},baseline=0cm]
\brunnianhalftwoside
\draw[double=none,Red,line width=2pt] (1) -- (5) (2) -- (8) (3) -- (4) ;
\end{tikzpicture}

\[
- A^{4} + 1 - A^{-8} + A^{-12}
\]

The left\hyp{}hand hook ``unhooks'', and the diagram unlinks leaving \(n-2\) components.

\item
\tikzsetnextfilename{brunnianhalfdia3}
\autolabel
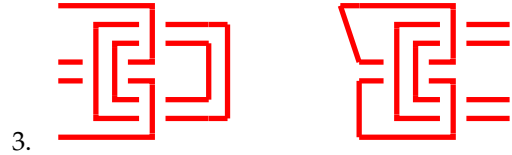
\begin{tikzpicture}[every path/.style={thick knot,double=Red},every node/.style={text=black},baseline=0cm]
\brunnianhalftwoside
\draw[double=none,Red,line width=2pt] (1) -- (4) (2) -- (3) (5) -- (8) ;
\end{tikzpicture}

\[
- A^{4} + 2 - A^{-4}
\]

This is similar to the Brunnian chain with the double end, except that both ends are doubled. 

\item
\tikzsetnextfilename{brunnianhalfdia4}
\autolabel
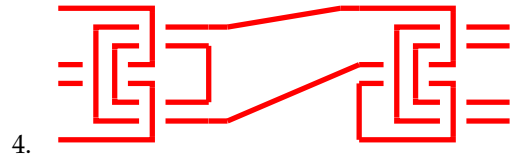
\begin{tikzpicture}[every path/.style={thick knot,double=Red},every node/.style={text=black},baseline=0cm]
\brunnianhalftwoside
\draw[double=none,Red,line width=2pt] (1) -- (5) (2) -- (3) (4) -- (8) ;
\end{tikzpicture}

\[
A^{10} - 2 A^{6} + A^{2} + A^{-2} - 2 A^{-6} + A^{-10}
\]

Here we have the same diagram as we started with, except with one fewer component.

\item
\tikzsetnextfilename{brunnianhalfdia5}
\autolabel
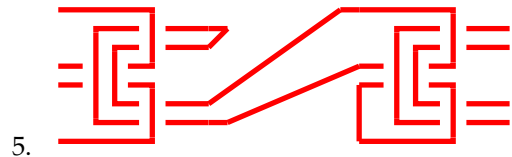
\begin{tikzpicture}[every path/.style={thick knot,double=Red},every node/.style={text=black},baseline=0cm]
\brunnianhalftwoside
\draw[double=none,Red,line width=2pt] (1) -- (2) (3) -- (5) (4) -- (8) ;
\end{tikzpicture}

\[
A^{12} - A^{8} + 1 - A^{-4}
\]

This unlinks, leaving \(n-2\) components.

\end{enumerate}

Taking \(n-2\) as our base, gathering the unlinked components together we obtain:
\[
(A^{12} + 2 + A^{-12})(-A^2 - A^{-2})^{n-3} = (A^{4} - 1 + A^{-4})^2(-A^2 - A^{-2})^{n-1}.
\]

The Brunnian chain with both ends doubled obeys the same recursion formula as the chain with only the left\hyp{}hand end doubled (except that there is one more component) because that analysis depended on resolving the links starting from the left\hyp{}hand end so the doubled right\hyp{}hand end does not come in to play until the end (except for counting components).
Writing \(\kb{brunnianchaindblends}{n}\) for this Brunnian chain with both ends doubled, we see that \(\kb{brunnianchaindblends}{4} = \kb{dblhopfchain}{4}\).

Putting this back, if we write \(\kb{brunnianhalfringplus}{n}\) for the Kauffman bracket of the Brunnian half\hyp{}ring as in Link~\ref{fig:brdia7} (and \(\kb{brunnianhalfringminus}{n}\) for \ref{fig:brdia9}) then:
\begin{dmath}
\label{eq:brhalfring}
\kb{brunnianhalfringplusminus}{n} =  (A^{4} - 1 + A^{-4})^2(-A^2 - A^{-2})^{n-1} + (-A^2 - A^{-2})^2 \kb{brunnianchaindblends}{n} + (A^{10} - 2 A^{6} + A^{2} + A^{-2} - 2 A^{-6} + A^{-10}) \kb{brunnianhalfringplusminus}{n-1}.
\end{dmath}
The starting points are
\begin{align*}
\kb{brunnianhalfringplus}{3} &= 2 + A^{-8} + A^{-16}, \\
\kb{brunnianhalfringminus}{3} &= A^{16} + A^{8} + 2.
\end{align*}

Returning to the Brunnian ring, writing \(\kb{brunnianring}{n}\) for the Kauffman bracket of the Brunnian ring with \(n\) components, we find that:
\begin{dmath}
\label{eq:brring}
\kb{brunnianring}{n} = \Big(-(-A^{-2} - A^2)^4 + 3(-A^{-2} - A^2)^2 + 5\Big)(-A^2-A^{-2})^{n-1} - (-A^{-2} - A^2)^2 \kb{brunnianchaindblends}{n+2} + (- A^{6} + A^{2} + A^{-2} - A^{-6}) (\kb{brunnianhalfringplus}{n+1} +\kb{brunnianhalfringminus}{n+1})  + (- A^{6} + A^{2} + A^{-2} - A^{-6}) \kb{brunnianring}{n-1}.
\end{dmath}

To find the starting point, as per usual we can assume that in the resolution diagrams then the two flanking linkages are in fact the same one.
This yields our starting point as:
\[
\kb{brunnianring}{2} = A^{4} + 2 + A^{-4} = (-A^2 - A^{-2})^2.
\]

\section{Level Two Brunnian Links}

As we did with the Hopf links, we now move to level two Brunnian links.
Viewing the Brunnian rings as the basic ring\hyp{}shape, we can weave these into Brunnian chains and rings.

We apply the same analysis to these as we did to the level two Hopf structures.
We choose a particular component.
This is a Brunnian ring, so we choose a particular Brunnian linkage and resolve it.
This leads to five diagrams (with appropriate weights), corresponding to the terms in \eqref{eq:brring}.
In the diagram corresponding to the first term, the Brunnian ring under consideration has unlinked completely.
The resulting level two structure therefore disassembles, leaving only a disconnected family of Brunnian rings.
In the diagram corresponding to the second term, the Brunnian ring under consideration becomes a Brunnian chain.
This also causes the level two structure to disassemble.
In the diagrams corresponding to the other terms, the Brunnian ring retains its ring\hyp{}like structure and so the level two structure remains.
However, the Brunnian ring under consideration is replaced, either by a half\hyp{}ring or by a shorter Brunnian ring.
We also need to consider the case where some of the Brunnian rings are actually Brunnian half\hyp{}rings.
A similar analysis based on \eqref{eq:brhalfring} holds.

Let us write \(\kb{ltwobrunnianchain}{n_1,\dotsc,n_k}\) and  \(\kb{ltwobrunnianring}{n_1,\dotsc,n_k}\) for the Kauffman bracket of a level two Brunnian chain or ring composed of \(k\) Brunnian rings where the \(j\)th has \(n_j\) components.
If one of the components is adorned, such as \(n^+\) or \(n^-\), then this indicates that this component is actually a Brunnian half\hyp{}ring of the corresponding polarity.
We have four formulae, depending on whether we are considering a ring or a chain, and whether the component that we are resolving is a Brunnian ring or a Brunnian half\hyp{}ring.
In these, \(n^\epsilon\) means one of \(n\), \(n^+\), or \(n^-\).
For convenience, we let \(\kb{brunnianring}{n^\pm}\) denote \(\kb{brunnianhalfringplusminus}{n}\).
In the recursive formulae, there is no difference between the level two Brunnian chain and the level two Brunnian ring, so we shall just state the formulae for the rings.

\begin{dgroup}
\begin{dmath*}
\kb{ltwobrunnianring}{n_1^{\epsilon_1},\dotsc,n_i,\dotsc,n_k^{\epsilon_k}} = 
\Big(\big(-(-A^{-2} - A^2)^4 + 3(-A^{-2} - A^2)^2 + 5\big)(-A^2-A^{-2})^{n_j-1} - (-A^{-2} - A^2)^2 \kb{brunnianchaindblends}{n_j+2}\Big) (-A^{-2} - A^2)^{k-1} \prod_{\substack{i=1 \\ i \ne j}}^k\kb{brunnianring}{n_i^{\epsilon_i}}
 + (- A^{6} + A^{2} + A^{-2} - A^{-6}) (\kb{ltwobrunnianring}{n_1^{\epsilon_1},\dotsc,(n_j + 1)^+,\dotsc,n_k^{\epsilon_k}} + \kb{ltwobrunnianring}{n_1^{\epsilon_1},\dotsc,(n_j + 1)^-,\dotsc,n_k^{\epsilon_k}})
+ (- A^{6} + A^{2} + A^{-2} - A^{-6}) \kb{ltwobrunnianring}{n_1^{\epsilon_1},\dotsc,n_j-1,\dotsc,n_k^{\epsilon_k}}.
\end{dmath*}
\begin{dmath*}
\kb{ltwobrunnianring}{n_1^{\epsilon_1},\dotsc,n_i^{\pm},\dotsc,n_k^{\epsilon_k}} =\Big(   (A^{4} - 1 + A^{-4})^2(-A^2 - A^{-2})^{n_i-1} + (-A^2 - A^{-2})^2 \kb{brunnianchaindblends}{n_i}\Big) (-A^{-2} - A^2)^{k-1} \prod_{\substack{i=1 \\ i \ne j}}^k\kb{brunnianring}{n_i^{\epsilon_i}}
 + (A^{10} - 2 A^{6} + A^{2} + A^{-2} - 2 A^{-6} + A^{-10}) \kb{ltwobrunnianring}{n_1^{\epsilon_1},\dotsc,(n_i - 1)^{\pm},\dotsc,n_k^{\epsilon_k}}.
\end{dmath*}
\end{dgroup}

We turn our attention to the starting points.
The recursion for the Brunnian half\hyp{}ring continues until we reach \(3\) components.
At this point we have the link in Figure~\ref{fig:brhalfring3}.
The other components of the level two structure pass through the upper region (marked in grey).
We want to resolve this diagram but keeping track of the grey region.
To do this, we resolve the linkage outlined in green.
There are two possibilities with either two circles enclosing the grey region or one circle disjoint from it.
In the latter case, the level two structure disassembles.
In the former, we have now replaced this component by two circles.
Putting in the coefficients, we therefore have:
\begin{dgroup}
\begin{dmath*}
\kb{ltwobrunnianring}{n_1^{\epsilon_1},\dotsc,3^+,\dotsc,n_k^{\epsilon_k}} = (A^{10} - A^{6} + A^{2} - A^{-6} + A^{-10} - A^{-14}) \kb{ltwobrunnianring}{n_1^{\epsilon_1},\dotsc,0^+,\dotsc,n_k^{\epsilon_k}} + (A^{12} + 3 - A^{-4} + A^{-8}) (-A^2 - A^{-2}) \prod_{\substack{i=1 \\ i \ne j}}^k\kb{brunnianring}{n_i^{\epsilon_i}}
\end{dmath*}
\begin{dmath*}
\kb{ltwobrunnianring}{n_1^{\epsilon_1},\dotsc,3^-,\dotsc,n_k^{\epsilon_k}} = (- A^{14} + A^{10} - A^{6} + A^{-2} - A^{-6} + A^{-10}) \kb{ltwobrunnianring}{n_1^{\epsilon_1},\dotsc,0^-,\dotsc,n_k^{\epsilon_k}} + (A^{8} - A^{4} + 3 + A^{-12}) (-A^2 - A^{-2}) \prod_{\substack{i=1 \\ i \ne j}}^k\kb{brunnianring}{n_i^{\epsilon_i}}
\end{dmath*}
\end{dgroup}
The \(0^\pm\) is to be interpreted as meaning that this component consists of two (unlinked) circles.

\tikzsetnextfilename{brhalfring3}

\begin{figure}
\centering
\begin{tikzpicture}[every path/.style={thick knot,double=Red},every node/.style={text=black},baseline=0cm]
\pgfmathsetmacro{\brover}{.25}
\pgfmathsetmacro{\bradj}{.05}
\pgfmathsetmacro{\brsep}{.25}
\pgfmathsetmacro{\brlen}{1}
\foreach \bri/\brdir/\brnud/\brv in {
  1/1/1/-1,
  2/-1/2/-1,
  3/-1/1/-1,
  4/1/1/1,
  5/1/1/-1,
  6/-1/1/1,
  7/-1/2/1,
  8/1/1/1
} {
  \coordinate (\bri) at (\brover * \brdir * \brnud,-\bri * \brsep + \brv * \bradj);
  \draw (\bri) -- ++(0,-\brv * \bradj) -- +(-\brdir * \brlen - \brover * \brdir * \brnud,0) coordinate (e\bri);
}
\draw (1) -- (4);
\draw (2) -- (7);
\draw (3) -- (6);
\draw (5) -- (8);
\draw[double=none,Red,line width=2pt]
  (e1) -- ++(0,2*\brsep) -| (e2)
  (e3) -- (e6)
  (e4) -- ++(-\brsep,0) -- ++(0,6*\brsep) -| ($(e7)+(\brsep,0)$) -- (e7)
  (e5) -- (e8)
;
\fill[gray] ($(e1)+(.5*\brsep,.5*\brsep)$) rectangle ($(e2)+(-.5*\brsep,2.5*\brsep)$);
\draw[double=none,Green,opacity=.5,rounded corners] ($(e1)+(-2*\brsep,.5*\brsep)$) rectangle ($(e7)+(2*\brsep,-2*\brsep)$);
\end{tikzpicture}
\caption{Brunnian half\hyp{}ring with \(3\) components}
\autolabel
\end{figure}

The recursion for the Brunnian ring continues until we reach the Brunnian ring with one component, as in Figure~\ref{fig:brring1}.
Again, the other components of the level two structure pass through the grey region so we resolve the linkage keeping track of this area.

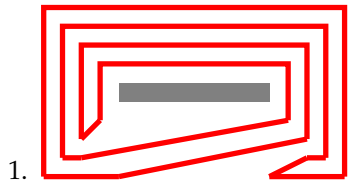
\tikzsetnextfilename{brring1}

\begin{figure}
\centering
\begin{tikzpicture}[every path/.style={thick knot,double=Red},every node/.style={text=black},baseline=0cm]
\pgfmathsetmacro{\brover}{.25}
\pgfmathsetmacro{\bradj}{.05}
\pgfmathsetmacro{\brsep}{.25}
\pgfmathsetmacro{\brlen}{1}
\foreach \bri/\brdir/\brnud/\brv in {
  1/1/1/-1,
  2/-1/2/-1,
  3/-1/1/-1,
  4/1/1/1,
  5/1/1/-1,
  6/-1/1/1,
  7/-1/2/1,
  8/1/1/1
} {
  \coordinate (\bri) at (\brover * \brdir * \brnud,-\bri * \brsep + \brv * \bradj);
  \draw (\bri) -- ++(0,-\brv * \bradj) -- +(-\brdir * \brlen - \brover * \brdir * \brnud,0) coordinate (e\bri);
}
\draw (1) -- (4);
\draw (2) -- (7);
\draw (3) -- (6);
\draw (5) -- (8);
\draw[double=none,Red,line width=2pt]
  (e1) -- ++(0,2*\brsep) -| (e2)
  (e4) -- ++(-\brsep,0) -- ++(0,6*\brsep) -| ($(e3)+(\brsep,0)$) -- (e3)
  (e5) -- ++(-2*\brsep,0) -- ++(0,8*\brsep) -| ($(e6)+(2*\brsep,0)$) -- (e6)
  (e8) -- ++(-3*\brsep,0) -- ++(0,12*\brsep) -| ($(e7)+(3*\brsep,0)$) -- (e7)
;
\fill[gray] ($(e1)+(.5*\brsep,.5*\brsep)$) rectangle ($(e2)+(-.5*\brsep,2.5*\brsep)$);
\draw[double=none,Green,opacity=.5,rounded corners] ($(e1)+(-4*\brsep,.5*\brsep)$) rectangle ($(e7)+(4*\brsep,-2*\brsep)$);
\end{tikzpicture}
\caption{Brunnian ring with \(1\) component}
\autolabel
\end{figure}
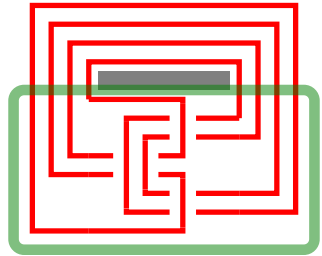

\begin{enumerate}

\pgfmathsetmacro{\brover}{.25}
\pgfmathsetmacro{\bradj}{.05}
\pgfmathsetmacro{\brsep}{.25}
\pgfmathsetmacro{\brlen}{1}

\item
\tikzsetnextfilename{brunnianringonedia1}
\autolabel
\begin{tikzpicture}[every path/.style={thick knot,double=Red},every node/.style={text=black},baseline=0cm]
\brunnianringonetwoside
\draw[double=none,Red,line width=2pt] (1) -- (2) (3) -- (5) (4) -- (6) (7) -- (8) ;
\end{tikzpicture}

\[
- A^{8} + A^{4}
\]

This unravels, leaving one circle.

\item
\tikzsetnextfilename{brunnianringonedia2}
\autolabel
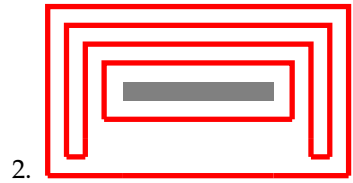
\begin{tikzpicture}[every path/.style={thick knot,double=Red},every node/.style={text=black},baseline=0cm]
\brunnianringonetwoside
\draw[double=none,Red,line width=2pt] (1) -- (5) (2) -- (3) (4) -- (8) (6) -- (7) ;
\end{tikzpicture}

\[
- A^{8} + 2 - A^{-8}
\]

The central circle comes out, leaving two circles in place.

\item
\tikzsetnextfilename{brunnianringonedia3}
\autolabel
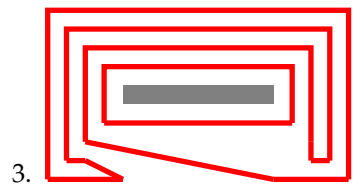
\begin{tikzpicture}[every path/.style={thick knot,double=Red},every node/.style={text=black},baseline=0cm]
\brunnianringonetwoside
\draw[double=none,Red,line width=2pt] (1) -- (5) (2) -- (8) (3) -- (4) (6) -- (7) ;
\end{tikzpicture}

\[
A^{-2} - A^{-10}
\]

The central part unravels slightly, leaving two circles surrounding the grey area.

\item
\tikzsetnextfilename{brunnianringonedia4}
\autolabel
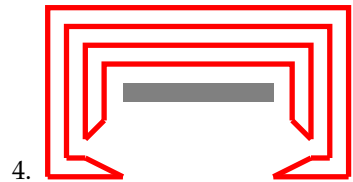
\begin{tikzpicture}[every path/.style={thick knot,double=Red},every node/.style={text=black},baseline=0cm]
\brunnianringonetwoside
\draw[double=none,Red,line width=2pt] (1) -- (2) (3) -- (4) (5) -- (6) (7) -- (8) ;
\end{tikzpicture}

\[
- A^{6} + A^{2} + A^{-2} - A^{-6}
\]

This unravels leaving two circles.

\item
\tikzsetnextfilename{brunnianringonedia5}
\autolabel
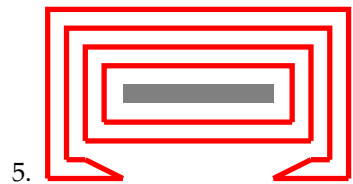
\begin{tikzpicture}[every path/.style={thick knot,double=Red},every node/.style={text=black},baseline=0cm]
\brunnianringonetwoside
\draw[double=none,Red,line width=2pt] (1) -- (5) (2) -- (6) (3) -- (4) (7) -- (8) ;
\end{tikzpicture}

\[
- A^{4} + 2 - A^{-8}
\]

The outermost circle can be removed, leaving two encircling the grey area.

\item
\tikzsetnextfilename{brunnianringonedia6}
\autolabel
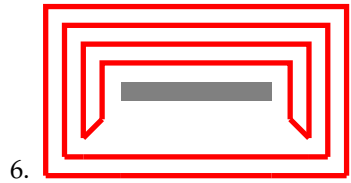
\begin{tikzpicture}[every path/.style={thick knot,double=Red},every node/.style={text=black},baseline=0cm]
\brunnianringonetwoside
\draw[double=none,Red,line width=2pt] (1) -- (2) (3) -- (7) (4) -- (8) (5) -- (6) ;
\end{tikzpicture}

\[
- A^{8} + 2 - A^{-4}
\]

The innermost circle can be removed, leaving two encircling the grey area.

\item
\tikzsetnextfilename{brunnianringonedia7}
\autolabel
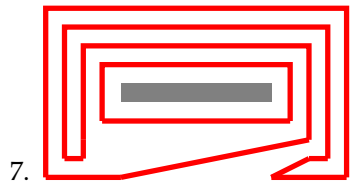
\begin{tikzpicture}[every path/.style={thick knot,double=Red},every node/.style={text=black},baseline=0cm]
\brunnianringonetwoside
\draw[double=none,Red,line width=2pt] (1) -- (5) (2) -- (3) (4) -- (6) (7) -- (8) ;
\end{tikzpicture}

\[
- A^{6} + A^{2} + A^{-2} - A^{-6}
\]

The outermost part unwinds, leaving two circles surrounding the grey area.

\item
\tikzsetnextfilename{brunnianringonedia8}
\autolabel
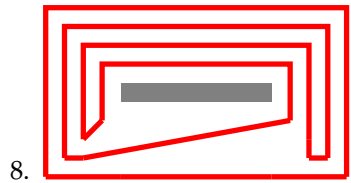
\begin{tikzpicture}[every path/.style={thick knot,double=Red},every node/.style={text=black},baseline=0cm]
\brunnianringonetwoside
\draw[double=none,Red,line width=2pt] (1) -- (2) (3) -- (5) (4) -- (8) (6) -- (7) ;
\end{tikzpicture}

\[
- A^{10} + A^{2}
\]

This simplifies to two circles surrounding the grey area.

\item
\tikzsetnextfilename{brunnianringonedia9}
\autolabel
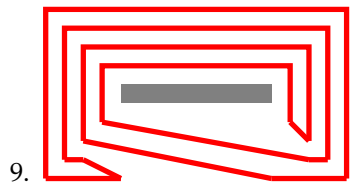
\begin{tikzpicture}[every path/.style={thick knot,double=Red},every node/.style={text=black},baseline=0cm]
\brunnianringonetwoside
\draw[double=none,Red,line width=2pt] (1) -- (7) (2) -- (8) (3) -- (4) (5) -- (6) ;
\end{tikzpicture}

\[
A^{-4} - A^{-8}
\]

This unravels, leaving one circle.

\item
\tikzsetnextfilename{brunnianringonedia10}
\autolabel
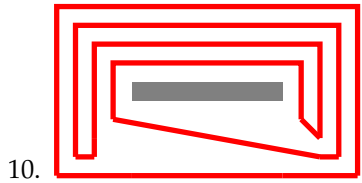
\begin{tikzpicture}[every path/.style={thick knot,double=Red},every node/.style={text=black},baseline=0cm]
\brunnianringonetwoside
\draw[double=none,Red,line width=2pt] (1) -- (7) (2) -- (3) (4) -- (8) (5) -- (6) ;
\end{tikzpicture}

\[
- A^{6} + A^{2} + A^{-2} - A^{-6}
\]

This simplifies to two circles surrounding the grey area.

\item
\tikzsetnextfilename{brunnianringonedia11}
\autolabel
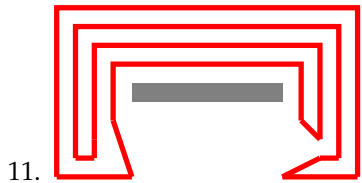
\begin{tikzpicture}[every path/.style={thick knot,double=Red},every node/.style={text=black},baseline=0cm]
\brunnianringonetwoside
\draw[double=none,Red,line width=2pt] (1) -- (4) (2) -- (3) (5) -- (6) (7) -- (8) ;
\end{tikzpicture}

\[
- A^{4} + 2 - A^{-4}
\]

This unravels leaving one circle.

\item
\tikzsetnextfilename{brunnianringonedia12}
\autolabel
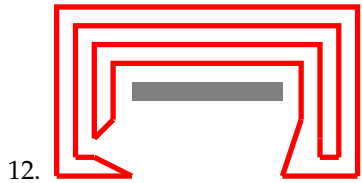
\begin{tikzpicture}[every path/.style={thick knot,double=Red},every node/.style={text=black},baseline=0cm]
\brunnianringonetwoside
\draw[double=none,Red,line width=2pt] (1) -- (2) (3) -- (4) (5) -- (8) (6) -- (7) ;
\end{tikzpicture}

\[
1
\]

This unravels leaving one circle.

\item
\tikzsetnextfilename{brunnianringonedia13}
\autolabel
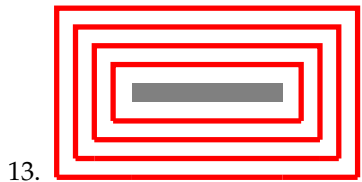
\begin{tikzpicture}[every path/.style={thick knot,double=Red},every node/.style={text=black},baseline=0cm]
\brunnianringonetwoside
\draw[double=none,Red,line width=2pt] (1) -- (5) (2) -- (6) (3) -- (7) (4) -- (8) ;
\end{tikzpicture}

\[
- A^{6} + A^{2} + A^{-2} - A^{-6}
\]

Here, there are four circles surrounding the grey area.

\end{enumerate}

Looking at the resulting diagrams, we see that there are three possibilities.
Either the component unlinks, resulting in the level two structure dissassembling, or it becomes a doubled component, or it becomes a quadrupled component.
We already have the notation \(0^\pm\) for the doubled component so we use a \(0\) to denote the quadruple component.
Putting in the coefficients, we therefore have:
\begin{dmath*}
\kb{ltwobrunnianring}{n_1^{\epsilon_1},\dotsc,1,\dotsc,n_k^{\epsilon_k}} =
(- A^{6} + A^{2} + A^{-2} - A^{-6}) \kb{ltwobrunnianring}{n_1^{\epsilon_1},\dotsc,0,\dotsc,n_k^{\epsilon_k}} + (A^{10} + A^{6} - 2 A^{2} - 2 A^{-2} + A^{-6} + A^{-10}) \kb{ltwobrunnianring}{n_1^{\epsilon_1},\dotsc,0^+,\dotsc,n_k^{\epsilon_k}} + (-A^2 - A^{-2}) \prod_{\substack{i=1 \\ i \ne j}}^k\kb{brunnianring}{n_i^{\epsilon_i}}
\end{dmath*}

From the above, our new starting points are links based on the Brunnian chain and Brunnian ring but where the components are doubled or quadrupled.
The doubling and quadrupling can be mixed, meaning that we cannot simply define the ``doubled Brunnian chain''.
Rather we have to define the doubled Brunnian chain (and ring) with a specification of which components are doubled and which quadrupled.
Thus we write, for example, \(\kb{dblbrunnianchain}{2,4,4,2,2}\) and \(\kb{dblbrunnianring}{2,4,4,2,2}\).
Thus we need to consider Brunnian linkages where the strands are doubled or quadrupled.
There are four possibilities to consider, as we need to consider doubling or quadrupling both sides independently.
Merely doubling both leads to a linkage with \(32\) crossings and \(16\) entry\enhyp{}exit points (thus \(8\) strands).
There are \(1430\) resulting diagrams.
That is rather a lot.

\appendix
\section{The kauffman.pl Program}
\label{ap:kauffman-code}

\lstset{language=Perl,breaklines=true}
\begin{lstlisting}
#!/usr/bin/perl

use strict;
use feature 'switch';
use Math::Polynomial::Laurent;

my $prefix = 'link';
my $split = ' ';

my $tikz = 0;
my $tex = 1;
my $raw = 0;
my $unlinked = 0;
my $debug = 0;

my $tex_config = {
    "power_prefix" => q({),
    "power_suffix" => q(}),
    "prefix" => q(
\[
),
    "suffix" => q(
\]
),
    "fold_sign" => 1,
    "variable" => "A"
};

my $p =  Math::Polynomial::Laurent->new();

my $deloop = $p->new(-2,[-1,0,0,0,-1]);

while (@ARGV) {
    my $arg = shift(@ARGV);

    given($arg) {
	when ('-prefix') {
	    $prefix = shift(@ARGV);
	}
	when ('-split') {
	    $split = shift(@ARGV);
	}
	when ('-tikz') {
	    $tikz = 1;
	}
	when ('-raw') {
	    $raw = 1;
	}
	when ('-notex') {
	    $tikz = 0;
	    $tex = 0;
	}
	when ('-debug') {
	    $debug = 1;
	}
    }

}

$tex && $p->string_config($tex_config);

my $crossings;
my $gatherings = [];

while (<>) {
    chomp;
    given ($_) {
	when (/^#/) {
	    # ignore comment lines
	}
	when (/^ *$/) {
	    # ignore blank lines
	}
	when (/^[0-9 ]+$/) {
	    # crossing specification
	    my @tmp = split(/$split/,$_);
	    push @$crossings, \@tmp;
	}
	when (/"(.*)"\s*=>\s*(\d+)/) {
	    # gathering specification
	    $$gatherings[0]->{$1} = $2;
	}
	when (/^\s*{/) {
	    # start of next gathering
	    unshift @$gatherings, {};
	}
	when (/^\s*prefix\s*=\s*(\w+)/) {
	    # prefix declaration
	    $prefix = $1;
	}
    }
}

my $numcross = @$crossings;
my $numdiag = 2**$numcross;
my %kauffman;
my %ends;
my $numgather = @$gatherings;

debug("$numcross crossings");
debug("$numdiag diagrams");
debug("$numgather gatherings");

foreach my $cross (@$crossings) {
    for (my $i = 0; $i < 4; $i++) {
	if (defined($ends{$cross->[$i]})) {
	    delete($ends{$cross->[$i]});
	} else {
	    $ends{$cross->[$i]} = 1;
	}
    }
}

debug("Ends: " . join(" ", keys(%ends)));

my $numstrans = keys(%ends);

for (my $i = 0; $i < $numdiag; $i++) {
    my $m = $i;
    my $n = 0;
    my $numcpts = 0;
    my $diag = [];
    for (my $j = 0; $j < $numcross; $j++) {
	my $k = $m % 2;
	$m = ($m - $k)/2;
	$n += $k;
	push @$diag, [$$crossings[$j][0], $$crossings[$j][2*$k + 1]];
	push @$diag, [$$crossings[$j][2], $$crossings[$j][3 - 2*$k]];
    }

    my %links = ();

    foreach my $link (@$diag) {
	if (exists $links{$$link[0]}) {
	    $links{$$link[0]} = [$links{$$link[0]},$$link[1]];
	} else {
	    $links{$$link[0]} = $$link[1];
	}
	if (exists $links{$$link[1]}) {
	    $links{$$link[1]} = [$links{$$link[1]},$$link[0]];
	} else {
	    $links{$$link[1]} = $$link[0];
	}
    }

    my $diagram = "";

    # Trace the strands that connect the ends

    foreach my $j (sort {$a <=> $b} (keys %ends))
    {
	next unless defined($links{$j});
	$diagram .= "($j)";
	my $g = $j;
	my $h = $links{$j};
	delete $links{$j};
	while (!exists($ends{$h})) {
	    if ($links{$h}[0] == $g) {
		$g = $h;
		$h = $links{$h}[1];
	    } else {
		$g = $h;
		$h = $links{$h}[0];
	    }
	    delete $links{$g};
	}
	delete $links{$h};
	$diagram .= " -- ($h) ";
    }

    # Find out how many loops we have left

    my @left = keys(%links);
    my $loops = 0;
    foreach my $j (@left)
    {
	next unless defined($links{$j});
	$loops++;
	my $g = $j;
	my $h = $links{$j}[0];
	delete $links{$j};
	while ($h != $j) {
	    if ($links{$h}[0] == $g) {
		$g = $h;
		$h = $links{$h}[1];
	    } else {
		$g = $h;
		$h = $links{$h}[0];
	    }
	    delete $links{$g};
	}
	delete $links{$h};
    }

    debug("Diagram: $diagram\nLoops: $loops");

    # Monomial according to the index of the crossing resolutions

    my $poly = $p->monomial(2 * $n - $numcross,1);

    debug("Monomial: " . $poly);

    # Multiply by the looping factor

    for (my $j = 0; $j < $loops; $j++ )
    {
	$poly = $poly * $deloop;
    }

    debug("With loops removed: " . $poly);

    # Initialise this component if we haven't seen it before 

    if (!exists $kauffman{$diagram}) {
	$kauffman{$diagram} = $p->new();
    }

    debug("Bracket before: " . $kauffman{$diagram});

    $kauffman{$diagram} = $kauffman{$diagram} + $poly;

    debug("Bracket after: " . $kauffman{$diagram});

}

my $item = 0;
foreach my $diagram (keys %kauffman) {
    $item++;
    print ($tikz ? '\item
\tikzsetnextfilename{' . $prefix . 'dia' . $item . '}
\autolabel
\begin{tikzpicture}[every path/.style={thick knot,double=Red},every node/.style={text=black},baseline=0cm]
\\' . $prefix . 'twoside
\draw[double=none,Red,line width=2pt] ' . $diagram . ';
\end{tikzpicture}
' : '%' . $diagram . "\n");

    print ($raw ? $kauffman{$diagram}->as_raw . "\n" : $kauffman{$diagram});
}

if ($numgather) {
    print "% Gatherings\n";
    for (my $i = 0; $i < $numgather; $i++) {
	debug("Gathering " . ($i + 1) . ":");
	my $poly = $p->new();
	foreach my $diag (keys %{$$gatherings[$i]}) {
	    debug($diag);
	    my $tpol = $kauffman{$diag} * $deloop ** $$gatherings[$i]->{$diag};
	    debug('(' . $tpol . ')(' . $deloop . ')^{' .  $$gatherings[$i]->{$diag} . '}');

	    $poly = $poly + $tpol;
	}
	print ($raw ? $poly->as_raw . "\n" : $poly);

    }

}

exit 0;

sub debug {
    my ($msg) = @_;
    if ($debug) {
	print STDERR $msg . "\n";
    }
}
\end{lstlisting}

\section{The kauffman-calc.pl Program}
\label{ap:kauffman-calc-code}

\lstset{language=Perl}
\begin{lstlisting}
#! /usr/bin/perl -w

use strict;
use feature 'switch';
use Math::Polynomial::Laurent;
use FreezeThaw;
use File::Basename;

my $link;
my $length;
my $linksub;
my $p = Math::Polynomial::Laurent->new();

my $tex_config = {
    "power_prefix" => q({),
    "power_suffix" => q(}),
    "prefix" => q(\[
),
    "suffix" => q(
\]
),
    "fold_sign" => 1,
    "variable" => "A"
};

$p->string_config($tex_config);

my $deloop = $p->m(-2,-1) + $p->m(2,-1);
my $pzero = $p->new(0,[0]);
my $raw = 0;
my $noice = 0;
my $dirname = dirname(__FILE__);
my $frozen_file = $dirname . '/../share/kauffman-on-ice';

my ($savedpolys,$formulae) = init();

while (@ARGV) {
    my $arg = shift(@ARGV);

    given ($arg) {
	when (/-link/) {
	    $link = shift(@ARGV);
	}
	when (/-length/) {
	    my @length = split(/[[:punct:][:space:]]/, shift(@ARGV));
	    $length = (@length == 1 ? $length[0] : \@length);
	}
	when (/-raw/) {
	    $raw = 1;
	}
	when (/-noice/) {
	    $noice = 1;
	}
    }
}

$noice || defrost($frozen_file,$savedpolys);


if (defined($link) && defined($length)) {
    my $result = apply_formula($link,$length);
    print ($raw ? $result->as_raw : $result);
} elsif (defined($link)) {
    foreach my $key (sort keys %{$savedpolys->{$link}}) {
	print $key . ": ";
	my $result = $savedpolys->{$link}->{$key};
	print ($raw ? $result->as_raw : $result);
    }
} else {
    foreach my $link (keys %$savedpolys) {
	print "-----\n$link:\n-----\n";
	foreach my $key (sort keys %{$savedpolys->{$link}}) {
	    print $key . ": ";
	    my $result = $savedpolys->{$link}->{$key};
	    print ($raw ? $result->as_raw : $result);
	}
    }
}

onice($frozen_file,$savedpolys);

exit;

# This is a wrapper subroutine around the actual recursive routines, it
# defines bail-outs and saves stuff so that they don't have to
sub apply_formula {
    my ($form,$len) = @_;
    my $slen;
    if (ref $len eq 'ARRAY') {
	my $nlen = @$len;
	for (my $i = 0; $i < $nlen; $i++) {
	    if ($$len[$i] < 0) {
		# Bail out if length is negative
		print STDERR "Error: $$len[$i] is less than zero\n";
		return $pzero;
	    }
	}
	$slen = join(':',@$len);
    } else {
	if ($len < 0) {
	    # Bail out if length is negative
	    print STDERR "Error: $len is less than zero\n";
	    return $pzero;
	}
	$slen = $len;
    }
    if (exists $savedpolys->{$form}->{$slen}) {
	# See if we've saved it from before
	return  $savedpolys->{$form}->{$slen};
    }
    if (!exists $formulae->{$form}) {
	# Check that we have a formula for this one
	print STDERR "No formula known for $form (of length $slen)\n";
	return $pzero;
    }
    # Apply the formula
    my $ans = $formulae->{$form}($len);

    # Save the result
    $savedpolys->{$form}->{$slen} = $ans;
    return $ans;
}

sub defrost {
    my ($icefile,$savehash) = @_;
    my $ice;
    my $water;

    if (open(ICE,$frozen_file)) {
	while (!eof(ICE)) {
	    $ice .= <ICE>;
	}
	($water) = FreezeThaw::thaw($ice);
	foreach my $link (keys %$water) {
	    foreach my $key (keys %{$water->{$link}}) {
		$savehash->{$link}->{$key} ||= $water->{$link}->{$key}
	    }
	}
	close ICE;
    }
    return 1;
}

sub onice {
    my ($icefile,$water) = @_;
    if (-e $icefile) {
	my $bicefile = $icefile . ".bak";
	if (-e $bicefile) {
	    unlink $bicefile;
	}
	link $icefile, $bicefile;
	unlink $icefile;
    }
    my $ice = FreezeThaw::freeze($water);
    open (ICE, ">$icefile")
	or die "Couldn't open $icefile for freezing";
    print ICE $ice;
    close ICE;
    return;
}

sub init {
    my ($spol,$form);
    $spol->{"hopf_chain"} =
    {
	"0" => $p->m(0,1),
	"1" => $p->m(-4,-1) + $p->m(4,-1)
    };
    $spol->{"hopf_ring"} =
    {
	"0" => $deloop,
	"1" => $p->m(-6,1)
    };
    $spol->{"double_hopf_chain"} =
    {
	"2" => $deloop,
	"4" => $p->m(14,-1) + $p->m(6,-1) + $p->m(2,-2) + $p->m(-2,-2) + $p->m(-6,-1) + $p->m(-14,-1)
    };
    $spol->{"double_hopf_halfring"} =
    {
	"3" => $p->m(16,1) + $p->m(8,1) + $p->m(0,2)
    };
    $spol->{"double_hopf_ring"} =
    {
	"2" => $p->m(18,-1) + $p->m(10,-1) + $p->m(6,1) + $p->m(2,-1),
	"4" => $p->new(-22,[-1,0,0,0,0,0,0,0,-1,0,0,0,0,0,0,
			    0,0,0,0,0,-3,0,0,0,-2,0,0,0,1,0,
			    0,0,1,0,0,0,-1,0,0,0,-2,0,0,0,0,
			    0,0,0,1,0,0,0,-1])
    };
    $spol->{"brunnian_chain_dbl_left"} =
    {
	"3" => $p->m(8,1) + $p->m(0,2) + $p->m(-8,1)
    };
    $spol->{"brunnian_chain"} =
    {
	"2" => $p->m(4,-1) + $p->m(-4,-1)
    };
    $spol->{"brunnian_chain_dbl_ends"} =
    {
	"4" => $p->new(-14,[-1,0,0,0,0,0,0,0,-1,0,0,0,-2,0,
			    0,0,-2,0,0,0,-1,0,0,0,0,0,0,0,-1])
    };
    $spol->{"brunnian_halfring_plus"} = 
    {
	"3" => $p->m(0,2) + $p->m(-8,1) + $p->m(-16,1)
    };
    $spol->{"brunnian_halfring_minus"} = 
    {
	"3" => $p->m(0,2) + $p->m(8,1) + $p->m(16,1)
    };
    $spol->{"brunnian_ring"} =
    {
	"2" => $deloop ** 2
    };
    $form = {
	"hopf_chain" =>
	    sub {
		my ($n) = @_;
		my $c = $p->m(-4,-1) +  $p->m(4,-1);
		my $t = apply_formula("hopf_chain",$n - 1);
		my $ans = $c * $t;
		return $ans;
	},
	"hopf_ring" => 
	    sub {
		my ($n) = @_;
		my $c = $p->m(2,1);
		my $t = apply_formula("hopf_chain",$n);
		my $d = $p->m(-4,-1) + $p->m(0,1);
		my $s = apply_formula("hopf_ring",$n - 1);
		my $ans = $c * $t + $d * $s;
		return $ans;
	},
	"double_hopf_chain" =>
	    sub {
		my ($n) = @_;
		my $c = $p->m(12,1) + $p->m(0,2) + $p->m(-12,1);
		my $t = apply_formula("double_hopf_chain",$n - 2);
		my $d = $p->m(12,-1) + $p->m(4,1) + $p->m(-4,1) + $p->m(-12,-1);
		my $s = apply_formula("double_hopf_chain",$n - 4);
		my $ans = $c * $t + $d * $s;
		return $ans;
	},
	"double_hopf_halfring" => 
	    sub {
		my ($n) = @_;
		my $c = $p->m(-10,1) + $p->m(-6,-1) + $p->m(6,-1) + $p->m(10,1);
		my $t = apply_formula("double_hopf_chain",$n - 1);
		my $d = $p->m(-10,-1) + $p->m(2,-1);
		my $s = apply_formula("double_hopf_chain",$n - 3);
		my $e = $p->m(-8,-1) + $p->m(0,1) + $p->m(4,-1) + $p->m(12,1);
		my $r = apply_formula("double_hopf_halfring",$n - 2);
		my $ans = $c * $t + $d * $s + $e * $r;
	},
	"double_hopf_ring" => 
	    sub {
		my ($n) = @_;
		my $c = $p->m(-8,1),
		my $t = apply_formula("double_hopf_chain",$n);
		my $d = $p->m(-8,-2) + $p->m(-4,2) + $p->m(0,2) + $p->m(4,-2);
		my $s = apply_formula("double_hopf_chain",$n-2);
		my $e = $p->m(-8,1) + $p->m(-4,-1) + $p->m(0,-2) + $p->m(4,2) + $p->m(8,1) + $p->m(12,-1);
		my $r = apply_formula("double_hopf_chain",$n-4);
		my $f = $p->m(-6,1) + $p->m(-2,1) + $p->m(2,-1) + $p->m(6,-1);
		my $q = apply_formula("double_hopf_halfring",$n-1);
		my $g = $p->m(0,1) + $p->m(4,-1) + $p->m(8,-1) + $p->m(12,1);
		my $o = apply_formula("double_hopf_ring",$n - 2);
		my $ans = $c * $t + $d * $s + $e * $r + $f * $q + $g * $o;
		return $ans;
	},
	"level2_hopf_chain" =>
	    sub {
		my ($n) = @_;
		my $ans;
		if (ref $n ne "ARRAY") {
		    $n = [$n];
		}
		my $a = @$n;
		if ($a == 1) {
		    $ans = apply_formula("hopf_ring",$n->[0]);
		    return $ans;
		}
		my $k;
		for (my $i = 0; $i < $a; $i++) {
		    if ($n->[$i] != 0) {
			$k = $i;
			last;
		    }
		}
		if (!defined($k)) {
		    $ans = apply_formula("double_hopf_chain",2 * @$n);
		    return $ans;
		}
		my $exp = 2;
		if (($k == 0) || ($k == $a - 1)) {$exp = 1};
		my $b = $p->m(2,1);
		my $c = $deloop ** $exp;
		my $d;
		if ($k != 0) {
		    $d = apply_formula("double_hopf_chain",2 * ($k + 1));
		} else {
		    $d = $p->m(0,1);
		}
		my $e = apply_formula("hopf_chain", $n->[$k]);
		my $f;
		if ($k != $a - 1) {
		    my @m = @$n[$k+1..$a-1];
		    $f = apply_formula("level2_hopf_chain",\@m);
		} else {
		    $f = $p->m(0,1);
		}
		my $g = $p->m(0,1) + $p->m(-4,-1);
		$n->[$k]--;
		my $h = apply_formula("level2_hopf_chain",$n);
		$ans = $b * $c * $d * $e * $f + $g * $h;
		return $ans;
		
	},
	"level2_hopf_ring" =>
	    sub {
		my ($n) = @_;
		my $ans;
		if (ref $n ne "ARRAY") {
		    $n = [$n];
		}
		my $a = @$n;
		if ($a == 1) {
		    $ans = apply_formula("hopf_ring",$n->[0]);
		    return $ans;
		}
		my $k;
		for (my $i = 0; $i < $a; $i++) {
		    if ($n->[$i] != 0) {
			$k = $i;
			last;
		    }
		}
		if (!defined($k)) {
		    $ans = apply_formula("double_hopf_ring",2 * @$n);
		    return $ans;
		}
		my $b = $p->m(2,1);
		my $c = $deloop;
		my $e = apply_formula("hopf_chain", $n->[$k]);
		my @m = @$n[$k+1..$a-1];
		push @m, @$n[0..$k-1];
		my $f = apply_formula("level2_hopf_chain",\@m);
		my $g = $p->m(0,1) + $p->m(-4,-1);
		$n->[$k]--;
		my $h = apply_formula("level2_hopf_ring",$n);
		$ans = $b * $c * $e * $f + $g * $h;
		return $ans;
	},
	"brunnian_chain_dbl_left" =>
	    sub {
		my ($n) = @_;
		my $a = $p->m(12,1) + $p->m(0,2) + $p->m(-12,1);
		my $b = $deloop ** ($n - 3);
		my $c = $p->m(10,1) + $p->m(6,-1) + $p->m(-6,-1) + $p->m(-10,1);
		my $d = apply_formula("brunnian_chain_dbl_left",$n-1);
		my $ans = $a * $b + $c * $d;
		return $ans;
	},
	"brunnian_chain_dbl_ends" =>
	    sub {
		my ($n) = @_;
		my $a = $p->m(12,1) + $p->m(0,2) + $p->m(-12,1);
		my $b = $deloop ** ($n - 3);
		my $c = $p->m(10,1) + $p->m(6,-1) + $p->m(-6,-1) + $p->m(-10,1);
		my $d = apply_formula("brunnian_chain_dbl_ends",$n-1);
		my $ans = $a * $b + $c * $d;
		return $ans;
	},
	"brunnian_chain" =>
	    sub {
		my ($n) = @_;
		my $a = $p->m(6,-1) + $p->m(-6,-1);
		my $b = $deloop ** ($n - 2);
		my $c = $p->m(4,-1) + $p->m(0,2) + $p->m(-4,-1);
		my $d = apply_formula("brunnian_chain_dbl_left",$n);
		my $ans = $a * $b + $c * $d;
		return $ans;
	},
	"brunnian_halfring_plus" =>
	    sub {
		my ($n) = @_;
		my $a = $p->m(4,1) + $p->m(0,-1) + $p->m(-4,1);
		my $b = $a ** 2;
		my $c = $deloop ** ($n - 1);
		my $d = $deloop ** 2;
		my $e = apply_formula("brunnian_chain_dbl_ends",$n);
		my $f = $p->m(10,1) + $p->m(6,-2) + $p->m(2,1) + $p->m(-2,1) + $p->m(-6,-2) + $p->m(-10,1);
		my $g = apply_formula("brunnian_halfring_plus",$n-1);
		my $ans = $b * $c + $d * $e + $f * $g;
		return $ans;
	},
	"brunnian_halfring_minus" =>
	    sub {
		my ($n) = @_;
		my $a = $p->m(4,1) + $p->m(0,-1) + $p->m(-4,1);
		my $b = $a ** 2;
		my $c = $deloop ** ($n - 1);
		my $d = $deloop ** 2;
		my $e = apply_formula("brunnian_chain_dbl_ends",$n);
		my $f = $p->m(10,1) + $p->m(6,-2) + $p->m(2,1) + $p->m(-2,1) + $p->m(-6,-2) + $p->m(-10,1);
		my $g = apply_formula("brunnian_halfring_minus",$n-1);
		my $ans = $b * $c + $d * $e + $f * $g;
		return $ans;
	},
	"brunnian_ring" =>
	    sub {
		my ($n) = @_;
		my $a = $p->m(0,5) + $p->m(0,3) * $deloop ** 2 - $deloop ** 4;
		my $b = $deloop ** ($n - 1);
		my $c = $deloop ** 2;
		my $d = apply_formula("brunnian_chain_dbl_ends",$n + 2);
		my $e = $p->m(6,-1) + $p->m(2,1) + $p->m(-2,1) + $p->m(-6,-1);
		my $f = apply_formula("brunnian_halfring_plus",$n+1) + apply_formula("brunnian_halfring_minus",$n+1);
		my $g = apply_formula("brunnian_ring",$n-1);
		my $ans = $a * $b - $c * $d + $e * $f + $e * $g;
	}
    };
    return ($spol,$form);
}
\end{lstlisting}

\section{The Math::Polynomial::Laurent Module}
\label{ap:laurent-code}

\lstset{language=Perl}
\begin{lstlisting}
package Math::Polynomial::Laurent;

use 5.006;
use strict;
use warnings;
use Carp qw(croak);
use Math::Polynomial 1.000;

require overload;

overload->import(
    q{neg}      => 'neg',
    q{+}        => _binary('add'),
    q{-}        => _binary('sub_'),
    q{*}        => _binary('mul'),
    q{/}        => _binary('div'),
    q{%}        => _binary('mod'),
    q{**}       => _lefty('pow'),
    q{<<}       => _lefty('shift_up'),
    q{>>}       => _lefty('shift_down'),
    q{!}        => 'is_zero',
    q{bool}     => 'is_nonzero',
    q{==}       => _binary('is_equal'),
    q{!=}       => _binary('is_unequal'),
    q{""}       => 'as_string',
    q{fallback} => undef,       # auto-generate trivial substitutions
);

# ---- object definition ----

# Math::Polynomial::Laurent=ARRAY(...)
use constant F_OFFSET     => 0; # degree offset
use constant F_POLY       => 1; # reference to Math::Polynomial object
use constant F_POLY_COEFF => 0; # reference to coefficients in Math::Polynomial
use constant F_CONFIG     => 2; # reference to hash of options for converting to string
use constant NFIELDS      => 3;

# ---- static data ----

our $VERSION = '0.001';
our $max_degree = 10_000;



# default values for as_string options
my @string_defaults = (
    ascending     => 0,
    with_variable => 1,
    fold_sign     => 0,
    fold_zero     => 1,
    fold_one      => 1,
    fold_exp_zero => 1,
    fold_exp_one  => 1,
    convert_coeff => sub { "$_[0]" },
    sign_of_coeff => undef,
    plus          => q{ + },
    minus         => q{ - },
    leading_plus  => q{},
    leading_minus => q{- },
    times         => q{ },
    power         => q{^},
    power_prefix  => q(),
    power_suffix  => q(),
    variable      => q{x},
    prefix        => q{(},
    suffix        => q{)},
);

my $global_string_config = {};

# ---- private methods ----

# generic polynomial detection hook (see Math::Polynomial::Generic)
sub _is_generic {
    return 0;
}

# binary operator wrapper generator
# generates functions to be called via overload:
# - upgrading a non-polynomial operand to a compatible polynomial
# - casting a generic operand if appropriate
# - restoring the original operand order
sub _binary {
    my ($method) = @_;
    return sub {
        my ($this, $that, $reversed) = @_;
        if (!ref($that) || !eval { $that->isa('Math::Polynomial::Laurent') }) {
            if ($this->_is_generic) {
                $that = Math::Polynomial::Laurent->new($that);
            }
            else {
                $that = $this->new($that);
            }
        }
        if ($this->_is_generic) {
            if (!$that->_is_generic) {
                $this = $this->_cast($that);
            }
        }
        elsif ($that->_is_generic) {
            $that = $that->_cast($this);
        }
        if ($reversed) {
            ($this, $that) = ($that, $this);
        }
        return $this->$method($that);
    };
}

# asymmetrically prototyped binary operator wrapper generator
# generates functions to be called via overload:
# - disallowing reverse order of operands
sub _lefty {
    my ($method) = @_;
    return sub {
        my ($this, $that, $reversed) = @_;
        croak 'wrong operand type' if $reversed;
        return $this->$method($that);
    };
}

# integer argument checker (Math::Polynomial checks N, we check Z)

sub _check_int {
    foreach my $arg (@_) {
	eval {
	    use warnings FATAL => 'all';
	    $arg == int $arg
	} or croak 'integer argument expected';
    }
    return;
}

# positive integer argument checker

sub _check_pos_int {
    foreach my $arg (@_) {
	eval {
	    use warnings FATAL => 'all';
	    $arg == abs int $arg
	} or croak 'positive integer argument expected';
    }
    return;
}

# --- methods ---

sub new
{
    # Should check that $off is an integer ...
    my ($this, $off, $coeff) = @_;
    my $class = ref $this;
    my ($poly,$config);
    if ($class) {
	(undef,undef,$config) = @{$this};
    } else {
	$config = undef;
	$class = $this;
    }
    if (ref $coeff eq "Math::Polynomial") {
	$poly = $coeff->clone;
    } else {
	$poly = Math::Polynomial->new(@$coeff);
    }
    return bless [$off, $poly,$config], $class;
}

sub clone
{
    my ($this) = @_;
    return bless [@{$this}], ref $this;
}

sub monomial {
    my ($this, $degree, $coeff) = @_;
    _check_int($degree);
    return $this->new($degree, [$coeff]);
}

# "monomial" gets a bit much to type each time

sub m
{
    my ($this, $degree, $coeff) = @_;
    return $this->monomial($degree,$coeff);
}

sub string_config {
    my ($this, $config) = @_;
    my $have_arg = 2 <= @_;
    if (ref $this) {
        if ($have_arg) {
            $this->[F_CONFIG] = $config;
        }
        else {
            $config = $this->[F_CONFIG];
        }
    }
    else {
        if ($have_arg) {
            # note: do not leave ultimate fallback configuration undefined
            $global_string_config = $config || {};
        }
        else {
            $config = $global_string_config;
        }
    }
    return $config;
}

sub coeff {
    my ($this, $degree) = @_;
    _check_int($degree);
    if ($this->is_zero || $degree < $this->[F_OFFSET]) {
	return $this->[F_POLY]->coeff_zero;
    }
    return $this->[F_POLY]->coeff($degree - $this->[F_OFFSET]);
}

sub topdegree {
    my ($this) = @_;
    if ($this->is_zero) {
	return 0;
    }
    my @coeffs = $this->[F_POLY]->coefficients;
    while ($coeffs[-1] == 0) { pop @coeffs };
    return @coeffs - 1 + $this->[F_OFFSET];
}

sub botdegree {
    my ($this) = @_;
    if ($this->is_zero) {
	return 0;
    }
    my @coeffs = $this->[F_POLY]->coefficients;
    my $offset = $this->[F_OFFSET];
    while ($coeffs[0] == 0) { shift @coeffs; $offset++ };
    return $offset;
}

sub evaluate {
    my ($this,$x) = @_;
    my $result = $this->[F_POLY]->evaluate($x);
    $result = $result * $x ** $this->[F_OFFSET];
    return $result;
}

sub is_zero {
    my ($this) = @_;
    return $this->[F_POLY]->is_zero;
}

sub is_nonzero {
    my ($this) = @_;
    return $this->[F_POLY]->is_nonzero;
}

sub is_equal {
    my ($this,$that) = @_;
    my $thisis;
    my $thatis;
    if ($this->[F_OFFSET] == $that->[F_OFFSET]) {
	$thisis = $this->[F_POLY];
	$thatis = $that->[F_POLY];
    } elsif ($this->[F_OFFSET] > $that->[F_OFFSET]) {
	$thisis = $this->[F_POLY];
	$thatis = $that->[F_POLY]->shift_up($this->[F_OFFSET] - $that->[F_OFFSET]);
    } else {
	$thisis = $this->[F_POLY]->shift_up($that->[F_OFFSET] - $this->[F_OFFSET]);
	$thatis = $that->[F_POLY];
    }
    return $thisis->is_equal($thatis);
}

sub is_unequal {
    my ($this,$that) = @_;
    my $thisis;
    my $thatis;
    if ($this->[F_OFFSET] == $that->[F_OFFSET]) {
	$thisis = $this->[F_POLY];
	$thatis = $that->[F_POLY];
    } elsif ($this->[F_OFFSET] > $that->[F_OFFSET]) {
	$thisis = $this->[F_POLY];
	$thatis = $that->[F_POLY]->shift_up($this->[F_OFFSET] - $that->[F_OFFSET]);
    } else {
	$thisis = $this->[F_POLY]->shift_up($that->[F_OFFSET] - $this->[F_OFFSET]);
	$thatis = $that->[F_POLY];
    }
    return $thisis->is_unequal($thatis);
}

sub neg {
    my ($this) = @_;
    return $this->new($this->[F_OFFSET], $this->[F_POLY]->neg);
}

sub add {
    my ($this,$that) = @_;
    my $thisis;
    my $thatis;
    my $offset;
    if ($this->is_zero) {
	return $that;
    }
    if ($that->is_zero) {
	return $this;
    }
    if ($this->[F_OFFSET] == $that->[F_OFFSET]) {
	$thisis = $this->[F_POLY];
	$thatis = $that->[F_POLY];
	$offset = $this->[F_OFFSET];
    } elsif ($this->[F_OFFSET] > $that->[F_OFFSET]) {
	$thisis = $this->[F_POLY]->shift_up($this->[F_OFFSET] - $that->[F_OFFSET]);
	$thatis = $that->[F_POLY];
	$offset = $that->[F_OFFSET];
    } else {
	$thisis = $this->[F_POLY];
	$thatis = $that->[F_POLY]->shift_up($that->[F_OFFSET] - $this->[F_OFFSET]);
	$offset = $this->[F_OFFSET];
    }
    return $this->new($offset, $thisis->add($thatis));
}

sub sub_ {
    my ($this,$that) = @_;
    my $thisis;
    my $thatis;
    my $offset;
    if ($this->is_zero) {
	return $that;
    }
    if ($that->is_zero) {
	return $this;
    }
    if ($this->[F_OFFSET] == $that->[F_OFFSET]) {
	$thisis = $this->[F_POLY];
	$thatis = $that->[F_POLY];
	$offset = $this->[F_OFFSET];
    } elsif ($this->[F_OFFSET] > $that->[F_OFFSET]) {
	$thisis = $this->[F_POLY]->shift_up($this->[F_OFFSET] - $that->[F_OFFSET]);
	$thatis = $that->[F_POLY];
	$offset = $that->[F_OFFSET];
    } else {
	$thisis = $this->[F_POLY];
	$thatis = $that->[F_POLY]->shift_up($that->[F_OFFSET] - $this->[F_OFFSET]);
	$offset = $this->[F_OFFSET];
    }
    return $this->new($offset, $thisis->sub_($thatis));
}

sub mul {
    my ($this,$that) = @_;
    return $this->new($this->[F_OFFSET] + $that->[F_OFFSET], $this->[F_POLY]->mul($that->[F_POLY]));
}

# Convert to a canonical form: polynomial has non-zero constant term

sub canonical {
    my ($this) = @_;
    my $return = clone($this);
    if ($this->is_zero) {
	$return->[F_OFFSET] = 0;
    } else {
	while ($return->[F_POLY]->coeff(0) == 0) {
	    $return->[F_POLY] = $return->[F_POLY]->shift_down(1);
	    $return->[F_OFFSET]++;
	}
    }
    return $return;
}

sub divmod {
    my ($this,$that) = @_;
    croak 'array context required' if !wantarray;
    my $cthis = $this->canonical;
    my $cthat = $that->canonical;
    my ($cquot, $crem) = $cthis->[F_POLY]->divmod($cthat->[F_POLY]);
    my $quot = $this->new($cthis->[F_OFFSET] - $cthat->[F_OFFSET], $cquot);
    my $rem = $this->new(- $cthis->[F_OFFSET], $crem);
    return ($quot,$rem);
}

sub div {
    my ($this,$that) = @_;
    my ($quot,$rem) = $this->divmod($that);
    return $quot;
}

sub mod {
    my ($this,$that) = @_;
    my ($quot,$rem) = $this->divmod($that);
    return $rem;
}

# Adapted from Math::Polynomial to allow for degree shift and braces around exponent

sub _make_ltz {
    my ($config, $zero) = @_;
    return 0 if !$config->{'fold_sign'};
    my $sgn = $config->{'sign_of_coeff'};
    return
        defined($sgn)?
            sub { $sgn->($_[0]) < 0     }:
            sub {        $_[0]  < $zero };
}

sub as_string {
    my ($this, $params) = @_;
    my %config = (
        @string_defaults,
        %{$params || $this->string_config || (ref $this)->string_config},
    );
    my $max_exp = $this->topdegree;
    my $min_exp = $this->botdegree;
    if ($max_exp < $min_exp) {
        $max_exp = $min_exp;
    }
    my $result = q{};
    my $zero = $this->[F_POLY]->coeff_zero;
    my $ltz  = _make_ltz(\%config, $zero);
    my $one  = $this->[F_POLY]->coeff_one;
    my $with_variable = $config{'with_variable'};
    foreach my $exp ($config{'ascending'}? $min_exp..$max_exp: reverse $min_exp..$max_exp) {
        my $coeff = $this->coeff($exp);

        # skip term?
        if (
            $with_variable &&
            $exp < $max_exp &&
            $config{'fold_zero'} &&
            $coeff == $zero
        ) {
            next;
        }

        # plus/minus
        if ($ltz && $ltz->($coeff)) {
            $coeff = -$coeff;
            $result .= $config{q[] eq $result? 'leading_minus': 'minus'};
        }
        else{
            $result .= $config{q[] eq $result? 'leading_plus': 'plus'};
        }

        # coefficient
        if (
            !$with_variable ||
            !$config{'fold_one'} ||
            0 == $exp && $config{'fold_exp_zero'} ||
            $one != $coeff
        ) {
            $result .= $config{'convert_coeff'}->($coeff);
            next if !$with_variable;
            if (0 != $exp || !$config{'fold_exp_zero'}) {
                $result .= $config{'times'};
            }
        }

        # variable and exponent
        if (0 != $exp || !$config{'fold_exp_zero'}) {
            $result .= $config{'variable'};
            if (1 != $exp || !$config{'fold_exp_one'}) {
                $result .= $config{'power'} . $config{'power_prefix'} . $exp . $config{'power_suffix'};
            }
        }
    }
    return join q{}, $config{'prefix'}, $result, $config{'suffix'};
}

sub shift {
    my ($this, $exp) = @_;
    _check_int($exp);
    croak 'exponent too large' if
        defined($max_degree) && $this->topdegree + $exp > $max_degree;
    croak 'exponent too small' if
        defined($max_degree) && $this->botdegree + $exp < -$max_degree;
    return $this if !$exp;
    return $this->new($this->[F_OFFSET] + $exp, $this->[F_POLY]);
}

sub shift_up {
    my ($this,$exp) = @_;
    return $this->shift($exp);
}

sub shift_down {
    my ($this,$exp) = @_;
    return $this->shift(-$exp);
}

sub pow {
    my ($this,$exp) = @_;
    _check_pos_int($exp);
    my $poly = $this->[F_POLY]->pow($exp);
    my $off = $this->[F_OFFSET] * $exp;
    return $this->new($off,$poly);
}

sub as_raw {
    my ($this,$can) = @_;
    $can ||= 0;
    my $poly;
    my $ret;
    if ($can) {
	$poly = $this->canonical;
    } else {
	$poly = $this;
    }
    $ret = $poly->[F_OFFSET] . ",[" . join(",",@{$poly->[F_POLY]->[F_POLY_COEFF]}) . "]";
    return $ret;
}

\end{lstlisting}

\chapter{Homfly Computations}

\mymaketitle

\section{Introduction}

In this article we compute the HOMFLY\hyp{}PT polynomial of some Brunnian rings.
To do the computations, we use the \texttt{homfly} program.
As our links have many strands and crossings, we use a program to produce the input suitable for the \texttt{homfly} program.
The code for this program is in Appendix~\ref{ap:homfly-code}.
From the HOMFLY\hyp{}PT polynomial it is possible to compute the Jones polynomial and the Alexander\emhyp{}Conway polynomial.
The code in Appendix~\ref{ap:homfly-code} can apply the necessary substitutions to compute these invariants.

\section{Brunnian Links}

The links that we are considering are all constructed from Brunnian linkages.
The basic component, a Brunnian linkage between two circles, is shown in Figure~\ref{fig:brunnianlink}.
As it stands, this is unlinked.
To make it linked, one has to consider this as part of a larger diagram.
With only pure Brunnian linkages, one of the simplest such diagrams is the ring\hyp{}like link in Figure~\ref{fig:brunnianring}.

\tikzsetnextfilename{brunnianlink}

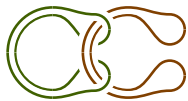
\begin{figure}
\centering
\begin{tikzpicture}
\pgfmathsetmacro{\brscale}{1}
\foreach \brk in {1,...,2} {
\pgfmathsetmacro{\brtint}{(\brk > 4 ? 8 - \brk : \brk) * 25}
\pgfmathsetmacro{\brl}{\brscale * \brk}
\begin{scope}[xshift=\brl cm]
\colorlet{chain}{Red!\brtint!Green}
\flatbrunnianlink{\brscale}
\end{scope}
}
\end{tikzpicture}
\caption{The Brunnian linkage}
\autolabel
\end{figure}

\tikzsetnextfilename{brunnianring}
\begin{figure}
\centering
\begin{tikzpicture}
\setbrstep{.1}
\foreach \brk in {1,2} {
\begin{scope}[rotate=\brk * 180]
\colorlet{chain}{ring\brk}
\brunnianlink{1}{180}
\end{scope}
}
\end{tikzpicture}
\caption{The Brunnian ring with two components}
\autolabel
\end{figure}
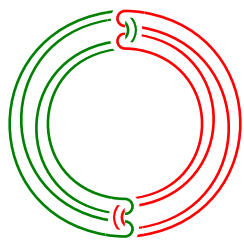

It is also possible to have a ``half\hyp{}Brunnian'' linkage, as in Figure~\ref{fig:halfbrunnian}.
Again, this is unlinked unless there is a larger diagram.
The simplest case now is with a ``half\hyp{}Brunnian'' linkage at either end, as in Figure~\ref{fig:brunnianchain}.
This is, incidentally, isotopic to the Borromean rings, shown in Figure~\ref{fig:borromean}.
From this, one can add more links in the middle to form a \emph{Brunnian chain}, as in Figure~\ref{fig:longbrunnianchain}.

\tikzsetnextfilename{halfbrunnian}
\begin{figure}
\centering
\begin{tikzpicture}
\colorlet{chain}{Red!50!Green}
\begin{pgfonlayer}{back}
\draw[knot,double=Red] (-.2,0) arc(0:180:.6);
\end{pgfonlayer}
\flatbrunnianlink{1}
\draw[knot,double=Red] (-.2,0) arc(0:-180:.6);
\end{tikzpicture}
\caption{Half\hyp{}Brunnian linkage}
\autolabel
\end{figure}
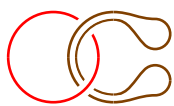

\tikzsetnextfilename{brunnianchain}
\begin{figure}
\centering
\begin{tikzpicture}
\colorlet{chain}{Red!50!Green}
\begin{pgfonlayer}{back}
\draw[knot,double=Red] (-.2,0) arc(0:180:.6);
\end{pgfonlayer}
\flatbrunnianlink{1}
\draw[knot,double=Red] (-.2,0) arc(0:-180:.6);
\draw[knot,double=Green] (1,0) circle (.6);
\end{tikzpicture}
\caption{Brunnian chain}
\autolabel
\end{figure}
\tikzsetnextfilename{borromean}

\begin{figure}
\centering
\begin{tikzpicture}
\pgfmathsetmacro{\brradius}{.6}
\pgfmathsetmacro{\csixty}{cos(60)}
\pgfmathsetmacro{\ssixty}{sin(60)}
\draw[knot,double=Red] (0,0) circle (\brradius);
\draw[knot,double=Green] (\brradius,0) circle (\brradius);
\begin{pgfonlayer}{back}
\draw[knot,double=Red!50!Green] (\brradius*\csixty,\brradius*\ssixty) ++(0,\brradius) arc(90:220:\brradius);
\draw[knot,double=Red!50!Green] (\brradius*\csixty,\brradius*\ssixty) ++(0,-\brradius) arc(-90:-10:\brradius);
\end{pgfonlayer}
\draw[knot,double=Red!50!Green] (\brradius*\csixty,\brradius*\ssixty) ++(0,\brradius) arc(90:-10:\brradius);
\draw[knot,double=Red!50!Green] (\brradius*\csixty,\brradius*\ssixty) ++(0,-\brradius) arc(-90:-140:\brradius);
\end{tikzpicture}
\caption{Borromean Rings}
\autolabel
\end{figure}
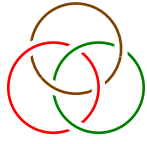

\tikzsetnextfilename{longbrunnianchain}
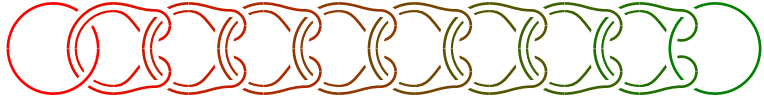
\begin{figure}
\centering
\begin{tikzpicture}
\colorlet{chain}{Red!50!Green}
\begin{pgfonlayer}{back}
\draw[knot,double=Red] (.8,0) arc(0:180:.6);
\end{pgfonlayer}
\pgfmathsetmacro{\brscale}{1}
\foreach \brk in {1,...,8} {
\pgfmathsetmacro{\brtint}{\brk * 11}
\pgfmathsetmacro{\brl}{\brscale * \brk}
\begin{scope}[xshift=\brl cm]
\colorlet{chain}{Green!\brtint!Red}
\flatbrunnianlink{\brscale}
\end{scope}
}
\draw[knot,double=Red] (.8,0) arc(0:-180:.6);
\draw[knot,double=Green] (9,0) circle (.6);
\end{tikzpicture}
\caption{Brunnian chain with 10 components}
\autolabel
\end{figure}

\subsection{Brunnain Chains}

The first links that we compute are the Brunnian chains.
The first non\hyp{}trivial Brunnian chain has three components and is isotopic to the Borromean rings (Figures~\ref{fig:brunnianchain} and~\ref{fig:borromean}).

\begin{enumerate}
\item Brunnian chain with three components (aka the Borromean rings).

\begin{itemize}
\item HOMFLY\hyp{}PT:
\begin{dmath*}
M^{-2}L^{-2} + 2M^{-2} + M^{-2}L^{2} - M^{2}L^{-2} - 2M^{2} - M^{2}L^{2} + M^{4}
\end{dmath*}

\item Jones polynomial:
\begin{dmath*}
- q^{3} + 3 q^{2} - 2 q + 4 - 2 q^{-1} + 3 q^{-2} - q^{-3}
\end{dmath*}

\item Alexander polynomial:
\begin{dmath*}
t^{2} - 4 t + 6 - 4 t^{-1} + t^{-2}
\end{dmath*}
\end{itemize}

\item Brunnian chain with four components.

\begin{itemize}
\item HOMFLY\hyp{}PT:
\begin{dmath*}
- M^{-3}L^{-3} - 3M^{-3}L^{-1} - 3M^{-3}L - M^{-3}L^{3} + 2M^{3}L^{-3} + 6M^{3}L^{-1} + 6M^{3}L + 2M^{3}L^{3} - M^{5}L^{-3} - 5M^{5}L^{-1} - 5M^{5}L - M^{5}L^{3} + M^{7}L^{-1} + M^{7}L
\end{dmath*}

\item Jones polynomial:
\begin{dmath*}
q^{-1/2} (- q^{6} + 4 q^{5} - 6 q^{4} + 5 q^{3} - 5 q^{2} - q - 1 - 5 q^{-1} + 5 q^{-2} - 6 q^{-3} + 4 q^{-4} - q^{-5})
\end{dmath*}

\item Alexander polynomial:
\begin{dmath*}
0
\end{dmath*}
\end{itemize}

\item Brunnian chain with five components.

\begin{itemize}
\item HOMFLY\hyp{}PT:
\begin{dmath*}
M^{-4}L^{-4} + 4M^{-4}L^{-2} + 6M^{-4} + 4M^{-4}L^{2} + M^{-4}L^{4} - 4M^{4}L^{-4} - 16M^{4}L^{-2} - 24M^{4} - 16M^{4}L^{2} - 4M^{4}L^{4} + 4M^{6}L^{-4} + 20M^{6}L^{-2} + 32M^{6} + 20M^{6}L^{2} + 4M^{6}L^{4} - M^{8}L^{-4} - 8M^{8}L^{-2} - 14M^{8} - 8M^{8}L^{2} - M^{8}L^{4} + M^{10}L^{-2} + 2M^{10} + M^{10}L^{2}
\end{dmath*}

\item Jones polynomial:
\begin{dmath*}
- q^{8} + 5 q^{7} - 10 q^{6} + 11 q^{5} - 8 q^{4} + q^{3} + 11 q^{2} - 13 q + 24 - 13 q^{-1} + 11 q^{-2} + q^{-3} - 8 q^{-4} + 11 q^{-5} - 10 q^{-6} + 5 q^{-7} - q^{-8}
\end{dmath*}

\item Alexander polynomial:
\begin{dmath*}
0
\end{dmath*}
\end{itemize}
\end{enumerate}

\subsection{Brunnian Rings}

The next links that we compute are the Brunnian rings.
The first non\hyp{}trivial Brunnian ring has two components and can be seen in Figure~\ref{fig:brunnianring}.
It can be deformed to a slightly simpler diagram which can be seen in Figure~\ref{fig:brunnianring2simple}.
The Brunnian ring with three components is shown in Figure~\ref{fig:brunnianring3}.

\begin{enumerate}
\item Brunnian ring with two components.

\begin{itemize}
\item HOMFLY\hyp{}PT:
\begin{dmath*}
- M^{-1}L^{-1} - M^{-1}L + ML^{-5} + ML^{-3} - 2ML^{-1} - 2ML + ML^{3} + ML^{5} - M^{3}L^{-3} + M^{3}L^{-1} + M^{3}L - M^{3}L^{3}
\end{dmath*}

\item Jones polynomial:
\begin{dmath*}
q^{-1/2} (- q^{6} + 2 q^{5} - 2 q^{4} + 3 q^{3} - 2 q^{2} - q - 1 - 2 q^{-1} + 3 q^{-2} - 2 q^{-3} + 2 q^{-4} - q^{-5})
\end{dmath*}

\item Alexander polynomial:
\begin{dmath*}
0
\end{dmath*}

\end{itemize}
\item Brunnian ring with three components.

\begin{itemize}
\item HOMFLY\hyp{}PT:
\begin{dmath*}
M^{-2}L^{-2} + 2M^{-2} + M^{-2}L^{2} + 2M^{2}L^{-4} + 8M^{2}L^{-2} + 12M^{2} + 8M^{2}L^{2} + 2M^{2}L^{4} + M^{4}L^{-6} + M^{4}L^{-4} - 7M^{4}L^{-2} - 14M^{4} - 7M^{4}L^{2} + M^{4}L^{4} + M^{4}L^{6} - M^{6}L^{-4} + M^{6}L^{-2} + 4M^{6} + M^{6}L^{2} - M^{6}L^{4}
\end{dmath*}

\item Jones polynomial:
\begin{dmath*}
- q^{8} + 5 q^{7} - 11 q^{6} + 14 q^{5} - 10 q^{4} + 11 q^{2} - 18 q + 24 - 18 q^{-1} + 11 q^{-2} - 10 q^{-4} + 14 q^{-5} - 11 q^{-6} + 5 q^{-7} - q^{-8}
\end{dmath*}

\item Alexander polynomial:
\begin{dmath*}
0
\end{dmath*}

\end{itemize}
\item Brunnian ring with four components.

\begin{itemize}
\item HOMFLY\hyp{}PT:
\begin{dmath*}
- M^{-3}L^{-3} - 3M^{-3}L^{-1} - 3M^{-3}L - M^{-3}L^{3} - 2M^{3}L^{-5} - 10M^{3}L^{-3} - 20M^{3}L^{-1} - 20M^{3}L - 10M^{3}L^{3} - 2M^{3}L^{5} - M^{5}L^{-7} + 5M^{5}L^{-5} + 17M^{5}L^{-3} + 19M^{5}L^{-1} + 19M^{5}L + 17M^{5}L^{3} + 5M^{5}L^{5} - M^{5}L^{7} + M^{7}L^{-7} + M^{7}L^{-5} - 7M^{7}L^{-3} + 5M^{7}L^{-1} + 5M^{7}L - 7M^{7}L^{3} + M^{7}L^{5} + M^{7}L^{7} - M^{9}L^{-5} - 9M^{9}L^{-1} - 9M^{9}L - M^{9}L^{5} + 2M^{11}L^{-1} + 2M^{11}L
\end{dmath*}

\item Jones polynomial:
\begin{dmath*}
q^{-1/2} (- q^{11} + 7 q^{10} - 24 q^{9} + 49 q^{8} - 56 q^{7} + 18 q^{6} + 51 q^{5} - 111 q^{4} + 131 q^{3} - 100 q^{2} + 32 q + 32 - 100 q^{-1} + 131 q^{-2} - 111 q^{-3} + 51 q^{-4} + 18 q^{-5} - 56 q^{-6} + 49 q^{-7} - 24 q^{-8} + 7 q^{-9} - q^{-10})
\end{dmath*}

\item Alexander polynomial:
\begin{dmath*}
0
\end{dmath*}

\end{itemize}
\end{enumerate}

Its HOMFLY\hyp{}PT polynomial is:
\begin{dmath*}
  M^{-2}L^{-2} + 2M^{-2} + M^{-2}L^{2} + 2M^{2}L^{-4} + 8M^{2}L^{-2} + 12M^{2} + 8M^{2}L^{2} + 2M^{2}L^{4} + M^{4}L^{-6} + M^{4}L^{-4} - 7M^{4}L^{-2} - 14M^{4} - 7M^{4}L^{2} + M^{4}L^{4} + M^{4}L^{6} - M^{6}L^{-4} + M^{6}L^{-2} + 4M^{6} + M^{6}L^{2} - M^{6}L^{4}.
\end{dmath*}

\tikzsetnextfilename{brunnianring2simple}

\begin{figure}
\centering
\begin{tikzpicture}[every path/.style={knot}]
\draw[double=ring1] (2.3,0) arc (0:45:2.3) (2.3,0) arc (0:-45:2.3);
\draw[double=ring1] (-2.3,0) arc (180:135:2.3) (-2.3,0) arc (-180:-135:2.3);
\draw[double=ring2] (-2,2) to[out=-135,in=-45] (2,2);
\draw[double=ring2] (2,-2) to[out=45,in=135] (-2,-2);
\draw[double=ring2] (2,2) to[out=135,in=-135] (2,-2);
\draw[double=ring2] (-2,-2) to[out=-45,in=45] (-2,2);
\draw[double=ring1] (0,2.3) arc (90:45:2.3) (0,2.3) arc (90:135:2.3);
\draw[double=ring1] (0,-2.3) arc (-90:-45:2.3) (0,-2.3) arc (-90:-135:2.3);
\end{tikzpicture}
\caption{The simplified form of the Brunnian ring with 2 components.}
\autolabel
\end{figure}
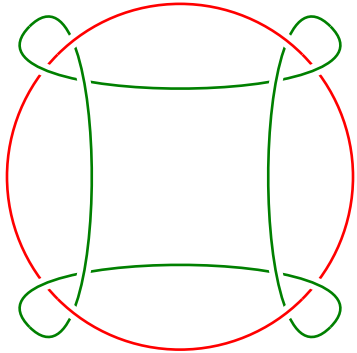

\tikzsetnextfilename{brunnianring3}

\begin{figure}
\centering
\begin{tikzpicture}
\foreach \brk in {1,2,3} {
\begin{scope}[rotate=\brk * 120 - 120]
\colorlet{chain}{ring\brk}
\brunnianlink{2.5}{120}
\end{scope}
}
\end{tikzpicture}
\caption{Brunnian ring with 3 components}
\autolabel
\end{figure}
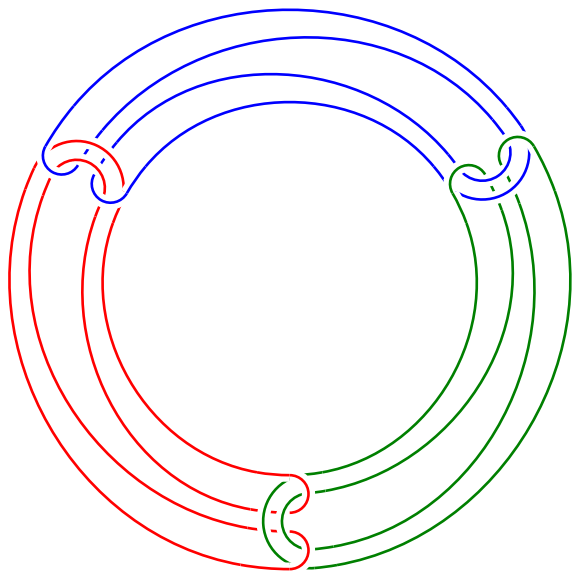

\appendix
\section{The homfly.pl Program}
\label{ap:homfly-code}

\lstset{language=Perl,breaklines=true}
\lstinputlisting{bin/homfly.pl}

\long\def\documentclass#1\begin#2{
\begingroup
}

\chapter{Brunnian Surfaces}

\mymaketitle

\section{Introduction}

The main purpose of this paper is to suggest new geometric forms with interesting properties for chemical synthesis.
The discussion is mathematical but we hope it will inspire chemists to work in this direction.

In \cite{nbnsas} we considered a method of building surfaces out of links.
The basic idea is similar to that of knitting: to ensure that the holes are sufficiently small that the appearance is that a genuine surface is created.
The distinction to knitting is that instead of using a single thread we use a multitude of components.

The inspiration for this was a generalisation of the family of so\hyp{}called \emph{rubberband links} or \emph{Brunnian rings} which are themselves a generalisation of the Borromean rings, see \cref{fig:borromean}.
It is straightforward to take a plethora of circles and link them together to make a surface \emhyp{} indeed, this is exactly what chain mail is \emhyp{} the challenge was to do so in such a way that it retained the key property of the Brunnian rings: that the removal of a single component caused the entire structure to disconnect.

\tikzsetnextfilename{borromean}
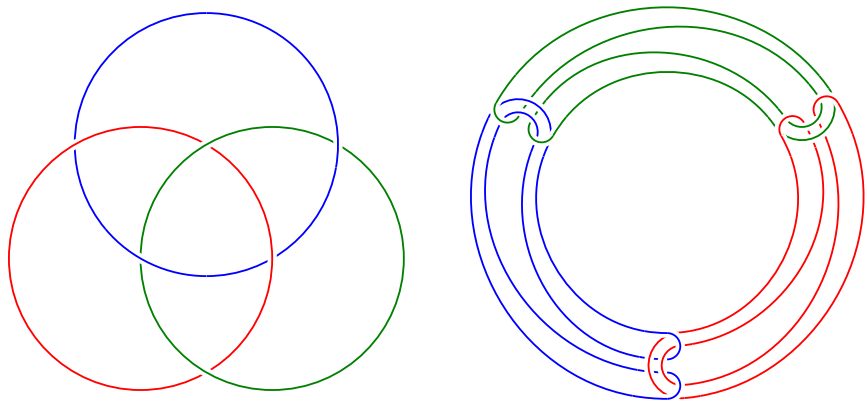
\begin{figure}
\centering
\scalebox{.7}{%
\begin{tikzpicture}[every path/.style={knot}]
\pgfmathsetmacro{\csixty}{cos(60)}
\pgfmathsetmacro{\ssixty}{sin(60)}
\pgfmathsetmacro{\brscale}{2.5}
\draw[double=ring1] (0,0) circle (\brscale);
\draw[double=ring2] (\brscale,0) circle (\brscale);
\begin{pgfonlayer}{back}
\draw[double=ring3] (\brscale *\csixty,\brscale *\ssixty) ++(0,\brscale) arc(90:220:2.5);
\draw[double=ring3] (\brscale *\csixty,\brscale *\ssixty) ++(0,-\brscale) arc(-90:-10:2.5);
\end{pgfonlayer}
\draw[double=ring3] (\brscale *\csixty,\brscale *\ssixty) ++(0,\brscale) arc(90:-10:\brscale);
\draw[double=ring3] (\brscale *\csixty,\brscale *\ssixty) ++(0,-\brscale) arc(-90:-140:\brscale);
\pgfmathsetmacro{\yshift}{\brscale * \ssixty/2}
\begin{scope}[xshift=10cm,yshift=\yshift cm]
\foreach \brk in {1,2,3} {
\begin{scope}[rotate=\brk * 120]
\colorlet{chain}{ring\brk}
\brunnianlink{\brscale}{120}
\end{scope}
}
\end{scope}
\end{tikzpicture}}
\caption{The Borromean Rings and the Brunnian Ring of Length \(3\)}
\autolabel
\end{figure}

The key to realising this is to observe that when a component is removed from the Brunnian ring then the way in which the rest of the components disconnect has some redundancy.
Not only does the disconnection proceed in both directions around the ring but also it is not necessary to fully disconnect one component in order to start on the next.
One can see this by the following.
In \cref{fig:borromean}, instead of removing the red component altogether, simply remove one of its loops from the green component.
This does not immediately separate the red and green components, but is enough to allow the blue component to slide off the red.
Once that has happened, the green can be disentangled from the blue and finally the last loop of the red can be removed from the green.

This leads one to consider building links from small components as in \cref{fig:chain} (possibly with small deformations) to produce structures like the carpet segment in \cref{fig:carpet}.

\tikzsetnextfilename{chain}
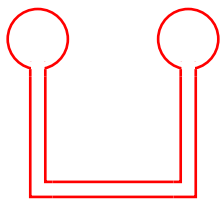
\begin{figure}
\centring
\begin{tikzpicture}
\colorlet{chain}{ring1}
\chain
\end{tikzpicture}
\caption{The Basic Component}
\autolabel
\end{figure}

\tikzsetnextfilename{carpet}
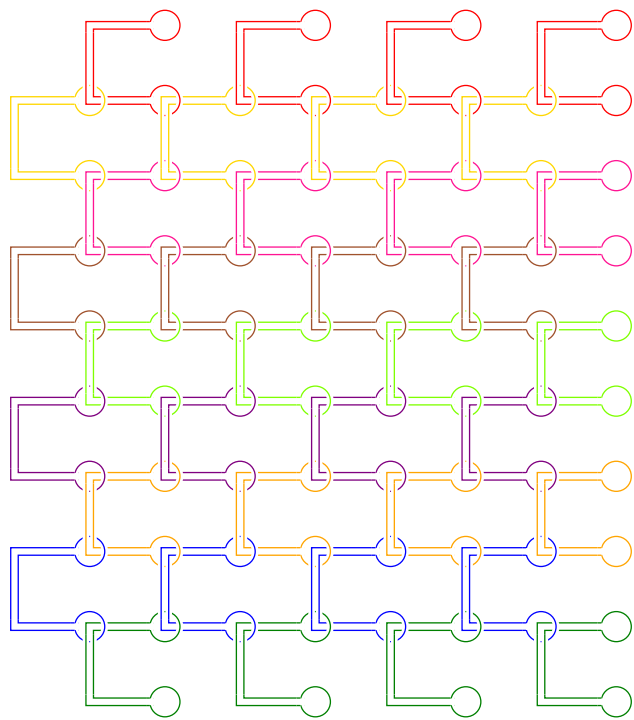
\begin{figure}
\centering
\scalebox{.5}{\begin{tikzpicture}
\foreach \k in {1,...,9} {
\foreach \l in {1,...,4} {
\pgfmathsetmacro{\m}{2*\l + 1*Mod(\k,2)}
\pgfmathparse{int(Mod(\k,9)+1)}
\edef\l{\pgfmathresult}
\begin{scope}[xshift=2*\m cm,yshift=2*\k cm,rotate=-90]
\colorlet{chain}{ring\l}
\chain
\end{scope}
}
}
\end{tikzpicture}}
\caption{Segment of a Brunnian carpet.}
\autolabel
\end{figure}

As drawn, the carpet is unlinked.
To form structures that hold together then the ends need to be joined together.
There are a variety of ways of doing this.
One of the simplest to draw is the carpet in \cref{fig:brcarpet}.

\tikzsetnextfilename{brcarpet}
\begin{figure}
\centering
\scalebox{.5}{%
\begin{tikzpicture}
\pgfmathsetmacro{\m}{0}
\pgfmathsetmacro{\last}{5}
\foreach \k in {1,...,\last} {
\pgfmathparse{int(\k < 5 ? \k : 9 - \k)}
\colorlet{chain}{ring1\pgfmathresult}
\pgfmathparse{\m + 4}
\xdef\m{\pgfmathresult}
\pgfmathsetmacro{\i}{-\k-1}
\pgfmathsetmacro{\j}{-\k}
\foreach \l in {1,...,\m} {
\pgfmathsetmacro{\side}{\l <= \m/4 ? 1 : (\l <= \m/2 ? 2 : (\l <= 3*\m/4 ? 3: 4))}
\pgfmathparse{\l <= \m/4 ? \i + 2 : (\l <= \m/2 ? \k : (\l == \m/2 + 1 ? \k - 1 : (\l <= 3*\m/4 ? \i - 2: -\k)))}
\xdef\i{\pgfmathresult}
\pgfmathparse{\l <= \m/4 ? -\k : (\l == \m/4 + 1 ? -\k + 1 : (\l <= \m/2 ? \j + 2 : (\l <= 3*\m/4 ? \k: (\l == 3*\m/4 + 1 ? \k - 1 : \j - 2))))}
\xdef\j{\pgfmathresult}
\begin{scope}[xshift=2*\i cm,yshift=2*\j cm,rotate=\side * 90]
\pgfmathparse{\k == \last ? "\noexpand\leftchain" : "\noexpand\chain"}
\let\tempchain=\pgfmathresult
\pgfmathparse{\k == \last && Mod(\l-1,\m/4) == 0 ? "\noexpand\leftcornerchain" : "\noexpand\tempchain"}
\pgfmathresult
\end{scope}
}
}
\end{tikzpicture}}
\caption{A Brunnian carpet.}
\autolabel
\end{figure}
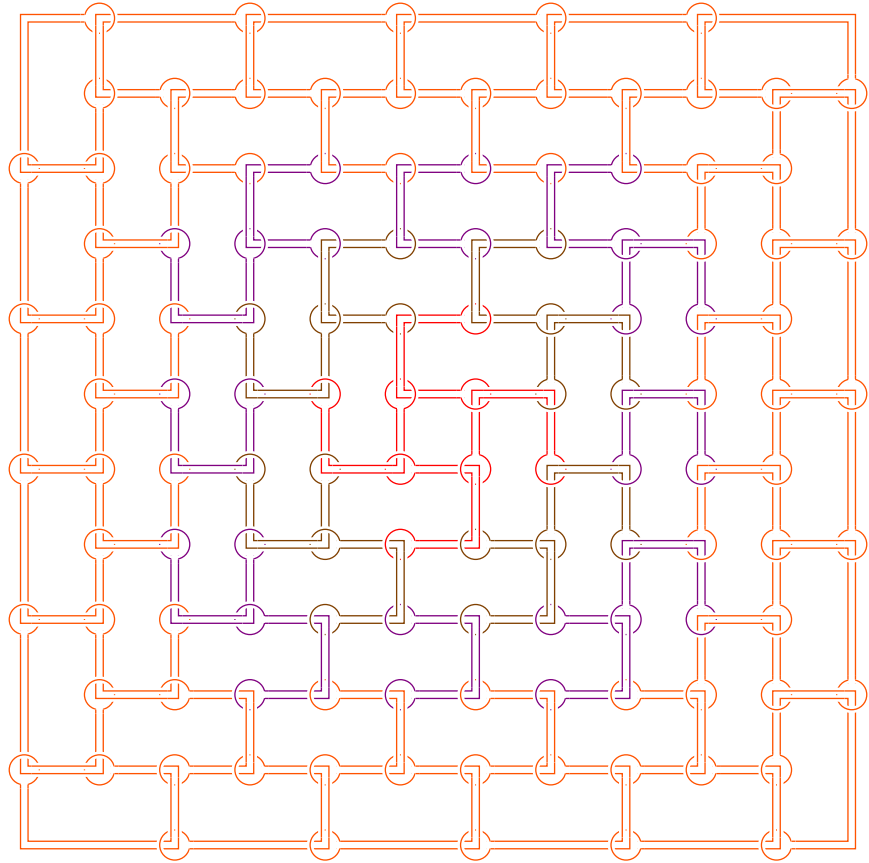

It is not necessarily the case that a structure built in this manner will have the Brunnian property wherein the removal of a single component causes it to fall apart.
One goal of this article is to introduce a framework whereby that question can be studied.
Within this coarse question of whether or not the structure has the Brunnian property are many finer questions relating to how easy it is to disconnect one of these structures.
Our framework will also address these.

\bigskip

To describe the framework we need to take a step back from our pictures.
From far off, the doubled lines merge into one and the loops look like blobs.
By colouring the components we can ensure that they remain visible.
Thus the simple carpet of \cref{fig:scarpet} looks like \cref{fig:scarpetgraph}.

\tikzsetnextfilename{scarpet}
\begin{figure}
\centring
\begin{tikzpicture}
\pgfmathsetmacro{\m}{0}
\pgfmathsetmacro{\last}{1}
\foreach \k in {1,...,\last} {
\pgfmathparse{\m + 4}
\xdef\m{\pgfmathresult}
\pgfmathsetmacro{\i}{-\k-1}
\pgfmathsetmacro{\j}{-\k}
\foreach \l in {1,...,\m} {
\pgfmathsetmacro{\side}{\l <= \m/4 ? 1 : (\l <= \m/2 ? 2 : (\l <= 3*\m/4 ? 3: 4))}
\pgfmathparse{\l <= \m/4 ? \i + 2 : (\l <= \m/2 ? \k : (\l == \m/2 + 1 ? \k - 1 : (\l <= 3*\m/4 ? \i - 2: -\k)))}
\xdef\i{\pgfmathresult}
\pgfmathparse{\l <= \m/4 ? -\k : (\l == \m/4 + 1 ? -\k + 1 : (\l <= \m/2 ? \j + 2 : (\l <= 3*\m/4 ? \k: (\l == 3*\m/4 + 1 ? \k - 1 : \j - 2))))}
\xdef\j{\pgfmathresult}
\begin{scope}[xshift=2*\i cm,yshift=2*\j cm,rotate=\side * 90]
\colorlet{chain}{ring\l}
\pgfmathparse{\k == \last ? "\noexpand\leftchain" : "\noexpand\chain"}
\let\tempchain=\pgfmathresult
\pgfmathparse{\k == \last && Mod(\l-1,\m/4) == 0 ? "\noexpand\leftcornerchain" : "\noexpand\tempchain"}
\pgfmathresult
\end{scope}
}
}
\end{tikzpicture}
\caption{Simple Carpet}
\autolabel
\end{figure}
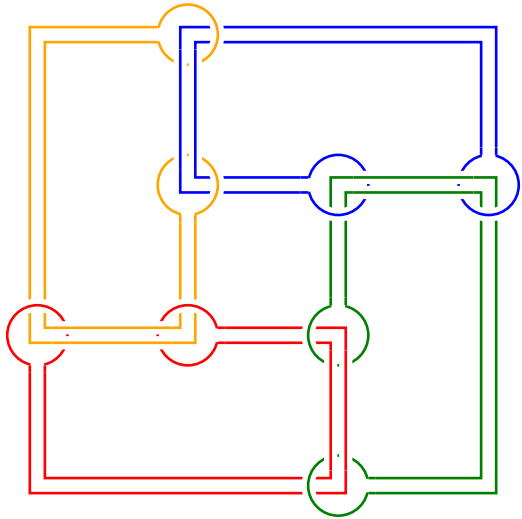

\tikzsetnextfilename{scarpetgraph}
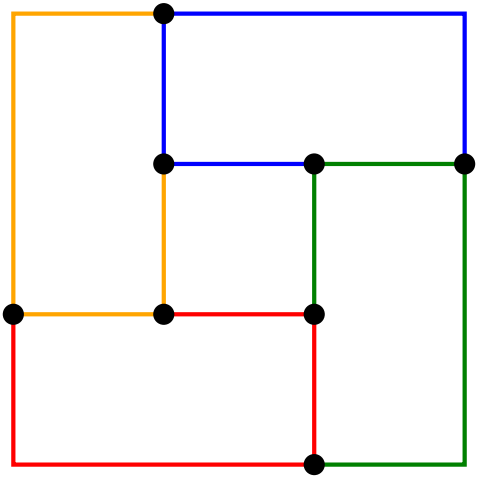
\begin{figure}
\centring
\begin{tikzpicture}
\foreach \k in {1,...,4}
{
  \draw[rotate=90*\k-90,ultra thick,ring\k] (-3,-1) -- ++(0,-2) -- ++(4,0) -- ++(0,2) -- ++(-2,0);
}
\foreach \k in {1,...,4}
{
  \fill[rotate=90*\k-90] (-3,-1) circle[radius=4pt] ++(2,0) circle[radius=4pt];
}
\end{tikzpicture}
\caption{Simple Carpet From a Distance}
\autolabel
\end{figure}

This looks like a graph in which the edges are coloured.
This will be our model for one of these structures.
We shall show that this captures the disconnection structure of the original link with one small exception.

To determine whether the link disconnects, from the graph we construct a category.
The objects of this category are subsets of the vertices of the graph.
The morphisms are, roughly, the disconnection implications.
Thus there is a morphism from one set of vertices to another if disconnecting the first set implies that the second set is also disconnected.

Once we have this category, our key question is to identify its inital objects.
These are collections of vertices with the property that disconnecting such a collection disconnects the entire graph.

There are various special circumstances that deserve names.
We shall say that the graph is \emph{Brunnian} if for each colour, the set of endpoint vertices forms an initial object.
We shall say that it is \emph{strongly Brunnian} if every vertex is initial.
We shall say that a component (or colour) is \emph{Brunnian} if its endpoint vertices are initial.

\section{The Link Graph}

In this section we shall define the graph corresponding to a link and show that the properties of disconnection for the two are the same, modulo one small exception.

We begin with the abstract definition.

\begin{defn}
A \emph{link graph} is a graph together with a colouring of its edges with the following properties:
\begin{enumerate}
\item The maximum valency of any vertex is \(3\).
\item The maximal monochrome subgraphs are trees.
\item Binary vertices are monochrome.
\item Vertices have at most two colours.
\end{enumerate}
\end{defn}

By a colouring of the edges we mean that we have a set of colours, say \(C\), and a function from the edges of the graph to \(C\).
When we talk of a vertex in chromatic terms then we mean the colours of the edges that end at that vertex.

We want to define a notion of equivalence for such graphs that models the process of unravelling.
The goal is to eliminate as many leaves as possible.

\begin{defn}
Two link graphs are said to be unravelled equivalent if they are joined by a series of the following moves (or their inverses):
\begin{enumerate}
\item Removal of binary vertices.
Since binary vertices are monochrome, the remaining edge has a well\hyp{}defined choice of colour.
\item Removal of a unary vertex where the adjacent vertex is monochrome.

That is to say, 
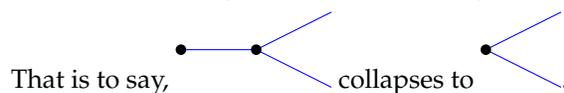
\tikz[baseline=-.5ex] { \draw[ring3] (0,0) -- (1,0) -- (2,.5) (1,0) -- (2,-.5); \fill (0,0) circle[radius=2pt] (1,0) circle[radius=2pt]; } collapses to \tikz[baseline=-.5ex] { \draw[ring3]  (1,0) -- (2,.5) (1,0) -- (2,-.5); \fill (1,0) circle[radius=2pt]; }.

\item Replacement of an adjacent unary vertex and non\hyp{}monochrome ternary vertex by two unary vertices.

That is to say, 
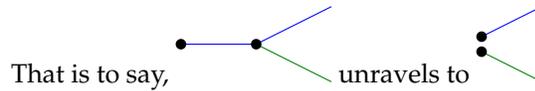
\tikz[baseline=-.5ex] { \draw[ring3] (0,0) -- (1,0) -- (2,.5); \draw[ring2] (1,0) -- (2,-.5); \fill (0,0) circle[radius=2pt] (1,0) circle[radius=2pt];} unravels to  \tikz[baseline=-.5ex] { \draw[ring3] (1.2,.1) -- (2,.5); \draw[ring2] (1.2,-.1) -- (2,-.5); \fill (1.2,.1) circle[radius=2pt] (1.2,-.1) circle[radius=2pt];}
\end{enumerate}
\end{defn}

It is obvious that applying any of these operations to a link graph produces a link graph.
Moreover, any unary vertex unless adjacent to another unary vertex can be collapsed along its branch until either that branch disappears or all that remains of that colour is two unary vertices joined by an edge.
Thus any link graph is equivalent to one of the following form.

\begin{defn}
A \emph{reduced} link graph is one in which each component either is trivalent or is of the form \tikz[baseline=-.5ex] { \draw[ring3] (0,0) -- (1,0); \fill (0,0) circle[radius=2pt] (1,0) circle[radius=2pt];}.

A link graph is \emph{unlinked} if it is equivalent to a reduced link graph with only monochrome components.
\end{defn}

In a reduced graph, a monochrome component must be of the form \tikz[baseline=-.5ex] { \draw[ring3] (0,0) -- (1,0); \fill (0,0) circle[radius=2pt] (1,0) circle[radius=2pt];}.
Thus a disconnected graph is equivalent to a disjoint union of such graphs.

\bigskip

Let us conclude this section by remarking on how to go from a link to a link graph and vice versa.
Recall that our links are built from basic components as in \cref{fig:chain}.
We attach such components by looping a circle of one component around a corner of another.
We can generalise these components by allowing more attaching points of both types, as in \cref{fig:genchain}.
Note that the receiving points are the genuine corners, not the branching points (as these would create a four valent vertex in the graph).

\tikzsetnextfilename{genchain}
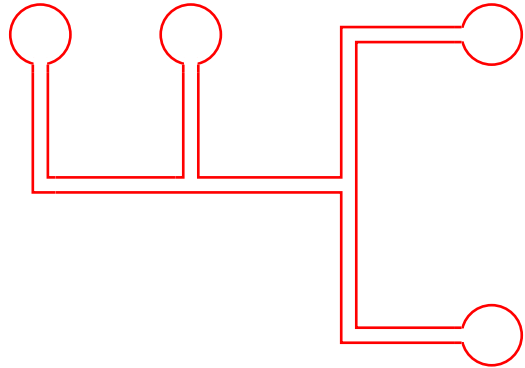
\begin{figure}
\centring
\begin{tikzpicture}
\colourlet{chain}{ring1}
\draw[knot,double=chain] (-1.1,0.7) -- ++(0,-1.8) -- ++(.3,0);
\draw[knot,double=chain] (1.1,0.7) -- ++(0,-1.6) -- ++(.3,0) -- ++(1.6,0) -- ++(0,2) -- ++(1.6,0);
\draw[knot,double=chain] (-.9,0.7) -- ++(0,-1.6) -- ++(.1,0);
\draw[knot,double=chain] (.9,0.7) -- ++(0,-1.6) -- ++(-.1,0);
\draw[knot,double=chain] (-.8,-.9) -- ++(1.6,0);
\draw[knot,double=chain] (-.8,-1.1) -- ++(3.8,0) -- ++(0,-2) -- ++(1.6,0);
\draw[knot,double=chain] (4.5,.9) -- ++(-1.3,0) -- ++(0,-3.8) -- ++(1.3,0);
\draw[knot,double=chain] (-1,1) circle (.4);
\draw[knotbg,line width=1.8mm] (-1,0.5) -- ++(0,.2);
\draw[knot,double=chain] (-1.1,0.5) -- ++(0,.1);
\draw[knot,double=chain] (-.9,0.5) -- ++(0,.1);
\draw[knot,double=chain] (5,1) circle (.4);
\draw[knotbg,line width=1.8mm] (4.5,1) -- ++(.2,0);
\draw[knot,double=chain] (4.5,1.1) -- ++(.1,0);
\draw[knot,double=chain] (4.5,.9) -- ++(.1,0);
\draw[knot,double=chain] (5,-3) circle (.4);
\draw[knotbg,line width=1.8mm] (4.5,-3) -- ++(.2,0);
\draw[knot,double=chain] (4.5,-3.1) -- ++(.1,0);
\draw[knot,double=chain] (4.5,-2.9) -- ++(.1,0);
\draw[knot,double=chain] (1,1) circle (.4);
\draw[knotbg,line width=1.8mm] (1,0.5) -- ++(0,.2);
\draw[knot,double=chain] (1.1,0.5) -- ++(0,.1);
\draw[knot,double=chain] (.9,0.5) -- ++(0,.1);
\end{tikzpicture}
\caption{Generalised Building Block}
\autolabel
\end{figure}

The process of going from the link to the graph is the obvious one: replace each double strand by an edge, or series of edges, and each joint by a vertex.
To go the other way we need to fatten the edges to double strands and then replace trivalent and univalent vertices by either branches or loops around strands.
This, however, is not uniquely defined because at each vertex we have a choice as to how to splice the strands coming in and going out.
Moreover, even with a particular choice of splicing we can introduce twists before splicing which will change the link type.
However, with a single exception, all the links so constructed will have the same properties with regard to disconnection and how they behave under removal of components.

The exception is straightforward to illustrate.
Consider the carpet from \cref{fig:scarpet} with one component removed, as in \cref{fig:cardis}.
In our scheme, this starts disconnecting by sliding out the green component from the blue, then the blue from the yellow.
However, it can also disconnect from the other direction since the yellow can be removed directly from the blue.
However, this involves separating the two strands of the yellow so that one can pass under the blue and one over.
Our methods cannot detect this move.

\tikzsetnextfilename{cardis}
\begin{figure}
\centring
\begin{tikzpicture}
\pgfmathsetmacro{\m}{0}
\pgfmathsetmacro{\last}{1}
\foreach \k in {1,...,\last} {
\pgfmathparse{\m + 4}
\xdef\m{\pgfmathresult}
\pgfmathsetmacro{\i}{-\k-1}
\pgfmathsetmacro{\j}{-\k}
\foreach \l in {2,...,\m} {
\pgfmathsetmacro{\side}{\l <= \m/4 ? 1 : (\l <= \m/2 ? 2 : (\l <= 3*\m/4 ? 3: 4))}
\pgfmathparse{\l <= \m/4 ? \i + 2 : (\l <= \m/2 ? \k : (\l == \m/2 + 1 ? \k - 1 : (\l <= 3*\m/4 ? \i - 2: -\k)))}
\xdef\i{\pgfmathresult}
\pgfmathparse{\l <= \m/4 ? -\k : (\l == \m/4 + 1 ? -\k + 1 : (\l <= \m/2 ? \j + 2 : (\l <= 3*\m/4 ? \k: (\l == 3*\m/4 + 1 ? \k - 1 : \j - 2))))}
\xdef\j{\pgfmathresult}
\begin{scope}[xshift=2*\i cm,yshift=2*\j cm,rotate=\side * 90]
\colorlet{chain}{ring\l}
\pgfmathparse{\k == \last ? "\noexpand\leftchain" : "\noexpand\chain"}
\let\tempchain=\pgfmathresult
\pgfmathparse{\k == \last && Mod(\l-1,\m/4) == 0 ? "\noexpand\leftcornerchain" : "\noexpand\tempchain"}
\pgfmathresult
\end{scope}
}
}
\end{tikzpicture}
\caption{Simple Carpet with One Component Removed}
\autolabel
\end{figure}
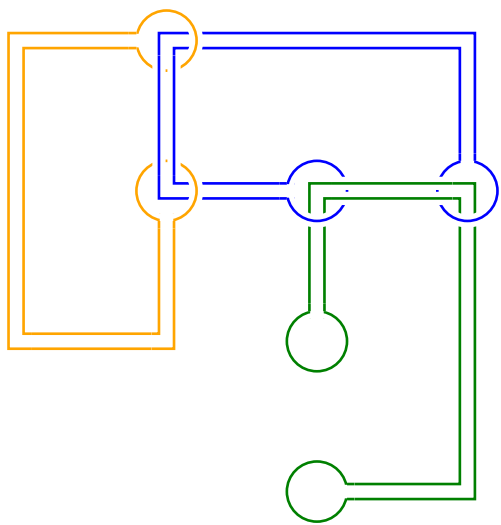

Notice that inserting a twist in the yellow strands disables this unlinking.
So our scheme describes an unlinking that works even if one does not look too closely at individual strands to see whether or not they are twisted around each other.

\bigskip

There is one other setup that is not covered in our scheme.
This is where a link attaches to itself.
To take this into account we would need to distinguish between two types of monochrome trivalent vertex: branching vertices and self\hyp{}attachment vertices.
It would also then be the case that a monochrome subgraph might not be equivalent to a two\hyp{}vertex graph.
A simple example of the sort of thing that we could get is in \cref{fig:selfconngraph}.
Here, the central vertex is a branching vertex and the outer vertices are connecting vertices.
At each connecting vertex, as we approach along the edge from the centre then we view the left\hyp{}hand edge as being the continuation and the right\hyp{}hand edge as being the connecting edge.

\tikzsetnextfilename{selfconngraph}
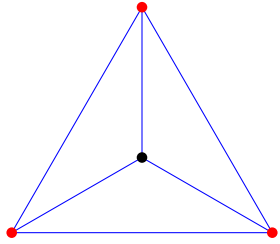
\begin{figure}
\centring
\begin{tikzpicture}
\draw[ring3] (0,0) -- (90:2) -- (210:2) (0,0) -- (-30:2) -- (90:2) (0,0) -- (210:2) -- (-30:2);
\fill (0,0) circle[radius=2pt];
\fill[ring1] (90:2) circle[radius=2pt] (-30:2) circle[radius=2pt] (210:2) circle[radius=2pt];
\end{tikzpicture}
\caption{A Link Graph With Self\hyp{}Connections}
\autolabel
\end{figure}

When we fatten this back to a link, we actually get a knot (as it has a single component).
After some simple manipulations, it looks like the alternating knot in \cref{fig:selfconnknot}.
This has Kauffman polynomial:
\[
A^{17} - 2 A^{13} + 2 A^{9} - 2 A^{5} + A + 3 A^{-11} - 2 A^{-15} + 2 A^{-19} - A^{-23}
\]

\tikzsetnextfilename{selfconnknot}
\begin{figure}
\centring
\begin{tikzpicture}[every path/.style={ring3}]
\foreach \k in {0,1,2}
{
  \begin{scope}[rotate=\k*120]
  \path (85:2) coordinate (a-0-\k) ++(30:.5) coordinate (a-1-\k)  ++(30:.5) coordinate (a-2-\k)
  (-25:2) coordinate (b-0-\k) ++(30:.5) coordinate (b-1-\k)  ++(30:.5) coordinate (b-2-\k);
\foreach \l in {0,1,2} 
  \draw (a-\l-\k) -- (b-\l-\k);
\draw (a-0-\k) -- ++(120:.35) arc[start angle=-150,delta angle=-47,radius=.5];
\draw (a-2-\k) -- ++(120:.35) arc[start angle=30,delta angle=110,radius=.5];
\draw (b-2-\k) -- ++(-60:.64) arc[start angle=30,delta angle=-120,radius=.5];
\draw (a-1-\k) -- ++(120:.74);
\draw (b-1-\k) -- ++(-60:.5) arc[start angle=30, delta angle=-180, radius=.25];
\draw (b-0-\k) -- ++(-60:.24);
  \end{scope}
}
\end{tikzpicture}
\caption{The Knot Corresponding to \Cref{fig:selfconngraph}}
\autolabel
\end{figure}

\section{The Disconnection Category}

In this section we shall introduce a category that will help us study the Brunnian properties of the link graph.
The goal is to determine conditions whereby we can identify which vertices need to be disconnected to completely unlink the graph.

Let \(G\) be a link graph.
We define a category, \(\m{D}_G\), as follows.

Let \(M_G\) be the set of vertices of \(G\) that are either univalent or non\hyp{}monochrome trivalent.
The objects of \(\m{D}_G\) are the subsets of \(M_G\).
For convenience of notation, we identify a vertex with the corresponding singleton subset.

We shall need some more notation.
Recall that for a given colour, the corresponding maximal monochrome subgraph of \(G\) is a tree.
For a colour \(c\) let us write \(G_c\) for this subgraph.
Let \(v \in M_G\) be a trivalent vertex.
As \(v\) is in \(M_G\), it cannot be monochrome and so must have exactly two colours.
One of these colour occurs twice, let us call this \(c\).
Consider \(G_c\).
In this graph, \(v\) is a bivalent vertex.
Its removal therefore splits \(G_c\) into two connected pieces.
We therefore can partition the leaves of \(G_c\) into two sets according to which piece of \(G_c \ssetminus \{v\}\) they end up in.
Let us write these two sets of vertices as \(L_+(v)\) and \(L_-(v)\)\footnote{We choose \(+\) and \(-\) merely as two distinguished symbols and so that we can use \(\pm\) to refer to both, there is no implication of a way to select the parity}.
Note that these are subsets of the leaves of \(G_c\), not of \(G_c \ssetminus \{v\}\).
Each vertex in \(L_\pm(v)\) has a single edge of colour \(c\).
By the rules for being a link graph, each vertex must therefore be either univalent or non\hyp{}monochrome trivalent, hence is in \(M_G\).
Thus \(L_\pm(v) \subseteq M_G\) and so are objects in \(\m{D}_G\).

The morphisms of \(\m{D}_G\) are generated by the following rules:
\begin{enumerate}
\item \(\m{D}_G\) is a \emph{thin} category; that is to say, there is at most one morphism between any two objects.
\item \(\{v_1, \dotsc, v_k\}\) is the categorical product of \(v_1, \dotsc, v_k\).
\item Let \(v\) be a univalent vertex.
Then there is a morphism \(\emptyset \to v\).
\item Let \(v \in M_G\) be a trivalent vertex.
Then there are morphisms \(L_\pm(v) \to v\).
\end{enumerate}

The main property of this category that we want to use is the following.

\begin{theorem}
Let \(G\) be a link graph and \(\m{D}_G\) its disconnection category.
Then disconnecting a set of non\hyp{}monochrome vertices in \(G\) unlinks it if and only if that set of vertices is initial in \(\m{D}_G\).
\end{theorem}

There is one very obvious initial object: \(M_G\) itself.
The question is as to the existence of others.
As the category is a thin category, an object is initial if and only if it has a morphism to \(M_G\).
Therefore, a na\"ive algorithm to find all initial objects is to produce a combinatorial description of the category and then find all objects that admit a morphism to \(M_G\).

We are particularly interested in graphs with particular properties.
The first properties relate to singleton subsets: we are interested in finding out which, if any, singleton subsets are initial and whether or not all singleton subsets are initial.

This latter case is equivalent to  the skeleton of \(\m{D}_G\) being a single point.

For a particular colour, say \(C\), we let \(\m{L}(C)\) be the subset of vertices of \(G\) with the property that one and only one of the incoming edges has colour \(C\).
Then we are interested in whether or not \(\m{L}(C)\) is initial, and for which colours this holds.

In \cref{ap:disprog} we describe an implementation in Perl of an algorithm to find these initial families.

Finally, we hope that the ideas presented here will facilitate the synthesis of Brunnian type surfaces.

\appendix

\section{The Disconnection Program}
\label{ap:disprog}
\lstset{
  language=Perl,
  breaklines=true,
  commentstyle={\upshape\small},
  tabsize=2,
  texcl}
\begin{lstlisting}
#! /usr/bin/env perl -w

use strict;
use Getopt::Long qw(:config auto_help bundling);

sub debug;
sub error;
sub full;
sub vertices;
sub colours;

my $debugging;
my @fns;

GetOptions (
  "d|debug!" => \$debugging,
  "f|full" => sub {push @fns, \&full},
    "v|vertices" => sub {push @fns, \&vertices},
    "c|colours" => sub {push @fns, \&colours},
    );

my @graph;

# The input to this program is a link graph, given either via \Verb+STDIN+ or as a file on the command line.  The graph should be specified as follows.	 First, enumerate all of the vertices consecutively starting from \(0\).  Each edge is then determined by the vertices at its ends and its colour.  This information is given to the program by listing the edges in the format: \Verb+<vertex> <vertex> <colour>+.  A colour is a string with no spaces (or hashes).  To enable comments, anything after a \Verb+#+ is ignored.  Lines not matching the format are also ignored.  Indentation is allowed.

# Internally, a graph is stored as an array.  The indexing set is the vertices (represented by their enumeration).  Each element of the graph is a hash containing the information relating to that vertex.  Specifically, it contains an array of information about its adjacent edges and its valency.  The information stored about an edge consists of an array containing (in order) the index of the vertex at the other end and the colour of the edge.	Note that each edge therefore appears twice in the structure.

while (<>)
{
# Strip out comments
  s/#.*//;
# Look for a match for an edge.
  if (/^\s*(\d+)\s+(\d+)\s+(\S+)/)
  {
# \Verb+$1+ and \Verb+$2+ are the vertices, \Verb+$3+ is the colour.

# Add the edge to the array of edges associated to vertex \Verb+$1+ and increment its valency.	Also check that the valency doesn't exceed \(3\).
    push @{$graph[$1]{"edges"}}, [$2,$3];
    $graph[$1]{"valency"}++;
    if ($graph[$1]{"valency"} > 3) {
      error ("Vertex $1 has valency greater than 3.");
    }

# Add the edge to the array of edges associated to vertex \Verb+$2+ and increment its valency.	Also check that the valency doesn't exceed \(3\).
    push @{$graph[$2]{"edges"}}, [$1,$3];
    $graph[$2]{"valency"}++;
    if ($graph[$1]{"valency"} > 3) {
      error ("Vertex $2 has valency greater than 3.");
    }
  }
}

# Not all vertices contribute to the objects of the category.  The array \Verb+@objects+ contains only those that do.  Its elements are the indices in \Verb+@graph+ of the contributing vertices.  We'll also need a reverse look\hyp{}up for \Verb+@objects+.  We can add this to the vertex hashes in \Verb+@graph+ for simplicity using the key \Verb+catindex+.

my @objects = ();

# Iterate through the graph, looking for suitable vertics.
for (my $i = 0; $i < @graph; $i++) {
  if ($graph[$i]{"valency"} == 1) {
# We include all univalent vertices.
    $graph[$i]{"catindex"} = @objects;
    push @objects,$i;
  } elsif ($graph[$i]{"valency"} == 3) {
# We also include all non\hyp{}monochrome trivalent vertices.  To find these, we look at the colour of the first edge and compare it with the other two.  If one of these differs, it cannot be monochrome. 
    my $nc = $graph[$i]{"edges"}[0][1];
    if (($graph[$i]{"edges"}[1][1] ne $nc) || ($graph[$i]{"edges"}[2][1] ne $nc)) {
      $graph[$i]{"catindex"} = @objects;
      push @objects,$i;
    }
  }
}

# Now that we have our list of contributing vertices we can build the category data structure.  The objects of our category are subsets of \Verb+@objects+ with each subset being the product of its elements (we conflate contributing vertices with the corresponding singleton subsets).  Since a morphism into a product is determined by morphisms into its factors, we only need to keep track of morphisms to singleton subsets, i.e.\ contributing vertices.  Thus the morphisms in our category are completely determined by the set of morphisms to the singleton subsets, and as our category is ``thin'', for each of these there is either a morphism or there isn't.  Thus the morphisms from an object are determined by a subset of the singleton sets (those it has a morphism to), which is again an object in the category.  We can encode our objects as integers in the range \(0\) to \(2^n - 1\) (where \(n\) is the number of objects) using the encoding \(S \mapsto \sum_{i \in S} 2^i\).  Thus our category data structure is an array of length \(2^n\) whose entries are integers in the range \([0, 2^n - 1]\).
my @category;
my $nobjs = 2**@objects;

# Initially, all we know is that the object at index \(i\) has a morphism to each of the singleton subsets that it contains, and so its morphism set is encoded again by \(i\).  So we initialise the category setting the value of \Verb+$category[$i]+ to \Verb+$i+.
for (my $i = 0; $i < $nobjs; $i++) {
  $category[$i] = $i;
}

# We now add in the morphims coming from the graph.

# Our definition of the category says that there is a morphism from the empty set to any univalent vertex.  So we iterate through the contributing vertex looking for univalent vertices.  If the \(i\)th such is univalent, we add \(2^i\) to the morphism encoding of the empty set.  Note that the empty set is encoded as \(0\), and the initial encoding of its morphisms is \(0\).
my $temp = 0;
for (my $i = 0; $i < @objects; $i++) {
  if ($graph[$objects[$i]]{"valency"} == 1) {
    $temp += 2**$i;
  }
}
if ($temp != 0) {
    $category[0] = $temp;
}

# Now we add in the morphisms from the contributing trivalent vertices.  For each, we traverse the corresponding monochrome graph to find the two sets of leaves.  These will be subsets of contributing vertices and so represent objects in the category.  We add in a morphism from each of these objects to the original vertex.

# Iterate over the contributing vertices.
for (my $i = 0; $i < @objects; $i++) {
# Find those of valency \(3\).
  if ($graph[$objects[$i]]{"valency"} == 3) {
# The first task is to determine the dominant colour of this vertex and find the two edges of that colour so that we can follow them to find the leaves.
    my $col;
    my @cedges;
# Do the first two edges have the same colour?
    if ($graph[$objects[$i]]{"edges"}[0][1] eq  $graph[$objects[$i]]{"edges"}[1][1]) {
# Yes.  So that must be the dominant colour, and the first two edges are the ones we need to follow.
      $col = $graph[$objects[$i]]{"edges"}[0][1];
      @cedges = (0,1);
    } else {
# No.  So the dominant colour is the colour of the last edge.
      $col = $graph[$objects[$i]]{"edges"}[2][1];
# And the two edges to follow are the last edge and the one of the first two that is of the same colour.
      if ($graph[$objects[$i]]{"edges"}[0][1] eq $col) {
	@cedges = (0,2);
      } else {
	@cedges = (1,2);
      }
    }
# Now that we have our initial edges, we follow each in turn along the original graph to find the leaves that are of the same colour.
    my $index;
# We do this for each edge in turn.
    foreach my $edge (@cedges) {
# The variable \Verb+$index+ will hold the integer corresponding to the subset of the contributing vertices that we end up finding.  So we start with the empty set.
      $index = 0;
# The array \Verb+@nedges+ consists of edges that we need to follow.  We add to this as we traverse the tree since it may have branches.  Each edge is represented by an array containing the starting vertex and the ending one.  We need to be careful to use the indices of the vertices in the original graph.
      my @nedges = ([$objects[$i],$graph[$objects[$i]]{"edges"}[$edge][0]]);
# Now we iterate through this array of edges.
      while (@nedges) {
# We get the next edge for consideration.
	my $e = shift @nedges;
# The variable \Verb+$add+ is a boolean for whether this edge leads to a leaf.  At the moment, we don't know that it doesn't so we set it to \Verb+1+ (true).
	my $add = 1;
# Now we look at each edge incident to the end vertex of the edge that we are considering.
	foreach my $te (@{$graph[$e->[1]]{"edges"}}) {
# We test to see if it is of the right colour and check that it isn't the edge we just came along.
	  if ($te->[1] eq $col && $te->[0] != $e->[0]) {
# If we have a match, we add this edge to the list of those to be considered and set our boolean \Verb+$add+ to \(0\) (false) as this vertex is not a leaf.
	    push @nedges, [$e->[1],$te->[0]];
	    $add = 0;
	  }
	}
# Once we've tested the incident edges we know whether or not the vertex at the end of this edge was a leaf.  If it was, we add its location to \Verb+$index+ to build up the object.
	if ($add) {
	  $index += 2**$graph[$e->[1]]{"catindex"};
	}
      }
# Once we've finished traversing the tree, we know the encoding of our object which has a morphism to the object \Verb+$i+.  To add that morphism, we \Verb+OR+ it with the existing encoding of the morphisms.
      debug("Adding morphism from $index to $i");
      $category[$index] |= 2**$i;
    }
  }
}

# The other type of initial object that we might be particularly interested in is the families of vertices determined by the leaves of a particular colour.  First, we need to build up a hash of these families.
my %colours;

# Each contributing vertex has a minority colour and is a leaf on the corresponding monochrome subtree.
for (my $i = 0; $i < @objects; $i++) {
# Check the valency.
  if ($graph[$objects[$i]]{"valency"} == 3) {
# It is \(3\), so we need to look for the minority colour.  Once we've found it, we add the index of the vertex to the hash of that colour (using the index in the \Verb+@objects+ array).
    if ($graph[$objects[$i]]{"edges"}[0][1] eq $graph[$objects[$i]]{"edges"}[1][1]) {
      $colours{$graph[$objects[$i]]{"edges"}[2][1]} |= 2**$i;
    } elsif ($graph[$objects[$i]]{"edges"}[1][1] eq $graph[$objects[$i]]{"edges"}[2][1]) {
      $colours{$graph[$objects[$i]]{"edges"}[0][1]} |= 2**$i;
    } else {
      $colours{$graph[$objects[$i]]{"edges"}[1][1]} |= 2**$i;
    }
  } else {
    $colours{$graph[$objects[$i]]{"edges"}[0][1]} |= 2**$i;
  }
}

foreach my $fn (@fns) {
  &$fn(\@category,\@objects,\%colours);
}

exit;

# These are two simple debugging and error messaging routines.
sub debug() {
    my $m = shift;
    if ($debugging) {
	print STDERR "DEBUG: " . $m . "\n";
    }
}

sub error() {
    my $m = shift;
    die "ERROR: " . $m . "\n";
}

sub full() {
  my $c = shift;
# Now that we have our initial data for the category, we need to complete it.  This works by looking for morphisms that are compositions of the ones that we already have.  Recall that the value of \Verb+$category[$i]+ points to the object that is the product of all of the singleton subsets that \Verb+$i+ has morphisms to.  It therefore has a morphism to any product of those subsets, which are represented by integers whose binary expansion is contained in that of \Verb+$category[$i]+.  To test this, we look at those objects \Verb+$j+ for which \Verb+$category[$i] & $j+ is again \Verb+$j+.  For such an object, any morphism out of \Verb+$j+ gives, by composition, a morphism out of \Verb+$i+.  Using our encoding, we can update the morphisms from \Verb+$i+ by \Verb+OR+ing it with the morphisms from \Verb+$j+.  That is, we replace \Verb+$category[$i]+ by \Verb+$category[$i] | $category[$j]+.  Before doing the update, we test to see if this would produce anything new since this is an iterative process and we want to know when we can safely stop.

  my $repeat = 1;
  my $nobjs = @$c;
  my ($i,$j,$t,$temp);
  my @output;
  my @s;

# While \Verb+$repeat+ is true, we loop through the objects.
  while ($repeat) {
# Set \Verb+$repeat+ to false initially.
    $repeat = 0;
# Iterate through the objects, collecting a list of the objects that we need to focus on.  Essentially, the value of \Verb+$c->[$i]+ is the maximal object that it is known to be equivalent to.
    for ($i = 0; $i < $nobjs; $i++) {
# We iterate over these objects of morphisms from \Verb+$i+.  If it is already everything then there is nothing further to do.
      next if ($c->[$i] == $nobjs - 1);
      $t = $c->[$i];
      $j = $t;
      $temp = 0;
      while (1) {
# This loop iterates over the subsets of \Verb+$t+, being those objects that have a morphism from \Verb+$t+, and gathers all the outgoing morphisms from those objects.
	$temp |= $c->[$j];
      } continue {
	last if $j == 0;
	last if $temp == $nobjs - 1;
	($j -= 1) &= $t;
      }
# Now we test to see if there are any new morphisms in that lot.  We are guaranteed to have \Verb+$temp & $t = $t+ since \Verb+$t+ represents morphisms from \Verb+$i+ and we already know about the identity morphism on \Verb+$i+.
      if ($t != $temp) {
# There are new morphisms, so add them in and flag for a repeat cycle.
	$repeat = 1;
	&debug("Updating $i from " . $c->[$i] . " to " . $temp);
# So now \Verb+$i+ has a morphism to \Verb+$temp+.  But we already had a morphism from \Verb+$temp+ to \Verb+$i+ by the way that \Verb+$temp+ was constructed.  So \Verb+$i+ and \Verb+$temp+ are equivalent, and so is everything in between.
	$t = $temp ^ $i;
	$j = $t;
	while (1) {
	    $c->[$j | $i] = $temp;
	} continue {
	    last if $j == 0;
	    ($j -= 1) &= $t;
	}
      }
    }
  }

  my @iobs;

  for ($i = 0; $i < $nobjs; $i++) {
    if ($c->[$i] == $nobjs - 1) {
      push @iobs, [$i,0];
    }
  }

  foreach my $iob (@iobs) {
    if (!$iob->[1]) {
      foreach my $ioc (@iobs) {
	if (!$ioc->[1] && ($iob->[0] != $ioc->[0])
	    && (($iob->[0] & $ioc->[0]) == $iob->[0])) {
	  $ioc->[1] = 1;
	}
      }
    }
  }

  foreach my $iob (@iobs) {
    if (!$iob->[1]) {
      @s = ();
      $t = $iob->[0];
      $j = 0;
      while ($t > 0) {
	if ($t & 1) {
	  push @s, $j;
	}
      } continue {
	$t >>= 1;
	$j ++;
      }
      push @output, "{" . join(",", @s) . "}";
    }
  }

  print "Initial objects: " . join(", ", @output) . " (and all supersets thereof)\n"; 

  return $c;
}

sub vertices () {
  my $c = shift;
  my $o = shift;
# This routine only looks for initial vertices.
  my @output;
  my ($i,$j,$t,$repeat);
  for ($i = 0; $i < @$o; $i++) {
    $repeat = 1;
    while ($repeat) {
      $repeat = 0;
# We start with vertex \Verb+$i+.  It has morphisms to \Verb+$c->[2**$i]+.  We iterate through the subsets of this, looking for morphisms.  If we already have morphisms to every object then there's nothing to do.
      $t = $c->[2**$i];
      next if ($t == $nobjs - 1);
      $j = $t;
      $temp = 0;
      while (1) {
	$temp |= $c->[$j];
      } continue {
	last if $j == 0;
	last if $temp == $nobjs - 1;
	($j -= 1) &= $t;
      }
      if ($temp != $t) {
	$repeat = 1;
	$c->[2**$i] = $temp;
# Add them to \Verb+$c->[$t]+ (the original target) as well.
	$c->[$t] = $temp;
      }
    }
  }

  for ($i = 0; $i < @$o; $i++) {
    if ($c->[2**$i] == $nobjs - 1) {
      push @output, $i;
    }
  }

  print "Initial vertices: " . join(", ", @output) . "\n"; 

  return $c;
}

sub colours () {
# Once we have our hash of colours, each pointing to an array of vertices, we need to convert each array of vertices to the index of an object in the category.  Then we test to see if that object is initial.  If so, we add that colour to our list to display.
  my $c = shift;
  my $o = shift;
  my $cols = shift;

  my @output = ();
  my ($i,$j,$t,$repeat);
  my $nobjs = @$c;

  foreach my $col (keys %$cols) {
    $i = $cols->{$col};
    $repeat = 1;
    while ($repeat) {
      $repeat = 0;
# We start with vertex \Verb+$i+.  It has morphisms to \Verb+$c->[$i]+.  We iterate through the subsets of this, looking for morphisms.  If we already have morphisms to every object then there's nothing to do.
      $t = $c->[$i];
      next if ($t == $nobjs - 1);
      $j = $t;
      $temp = 0;
      while (1) {
	$temp |= $c->[$j];
      } continue {
	last if $j == 0;
	last if $temp == $nobjs - 1;
	($j -= 1) &= $t;
      }
      if ($temp != $t) {
	$repeat = 1;
	$c->[$i] = $temp;
# Add them to \Verb+$c->[$t]+ (the original target) as well.
	$c->[$t] = $temp;
      }
    }
  }

  for my $col (keys %$cols) {
    $i = $cols->{$col};
    if ($c->[$i] == $nobjs - 1) {
      push @output, $col;
    }
  }
  print "Initial colours: " . join(", ", @output). "\n";

  return $c;
}
\end{lstlisting}

\bibliography{synthlinks}

\endgroup

\chapter{Programs}

\begingroup
\def\bibliography#1{}

\appendices
\endgroup
\bibliographystyle{plain}

\nocite{*}
%\bibliography{synthlinks}

\end{document}